\documentclass[12pt]{article}
\usepackage{amscd,amsfonts,amsmath,amssymb,amstext,amsthm}
\newtheorem{theorem}{Theorem}[section]
\newtheorem{lemma}{Lemma}[section]
\newtheorem{proposition}{Proposition}[section]
\newtheorem{corollary}{Corollary}[section]
\theoremstyle{remark}
\def\lie#1{\mathfrak{ #1}}
\def\CC{{\mathbb C}}
\def\RR{{\mathbb R}}
\def\ZZ{{\mathbb Z}}
\def\PP{{\mathbb P}}

\def\GL{\operatorname{GL}}
\def\PGL{\operatorname{PGL}}
\def\U{\operatorname{U}}
\def\SU{\operatorname{SU}}
\def\SL{\operatorname{SL}}
\def\SO{\operatorname{SO}}
\def\Sp{\operatorname{Sp}}
\def\USp{\operatorname{USp}}

\def\Spin{\operatorname{Spin}}

\def\Id{\operatorname{Id}}

\def\dim{\operatorname{dim}}
\def\codim{\operatorname{codim}}
\def\U{\operatorname{U}}

\def\Det{\operatorname{Det}}
\def\Fix{\operatorname{Fix}}

\def\bd{\operatorname{bd}}

\def\Id{\operatorname{Id}}

\def\Spin{\operatorname{Spin}}

\def\PSU{\operatorname{PSU}}
\def\Span{\operatorname{Span}}  
\title{Fibrations and globalizations of compact homogeneous
CR-manifolds\footnote{ This project was partially supported by an 
NSERC Discovery Grant and by SFB/TR 12 of the Deutsche 
Forschungsgemeinschaft.}}
\author{B. Gilligan and A. Huckleberry}
\date{\today}
\begin{document}
\maketitle
\begin{center}
{In memory of A.G. Vitushkin}  
\end{center}  
\tableofcontents
\begin {abstract}
\noindent
Fibrations methods which were previously used for complex homogeneous
spaces and CR-homogenous spaces of special types (\cite{HO1},
\cite{AHR}, \cite {HR,R}) are developed in a general framework.
These include the $\lie g$--anticanonical fibration in
the CR--setting which reduces certain considerations to
the compact projective algebraic case where a Borel--Remmert type
splitting theorem is proved.  This allows a reduction to
spaces homogeneous under actions of compact Lie groups.
General globalization theorems are proved which allow one
to regard the homogeneous CR--manifold as the orbit of a 
real Lie group in a complex homogeneous space of a complex
Lie group. In the special case of CR--codimension at most two
precise classification results are proved and are 
applied to show that in most cases there exists such a globalization.
\end {abstract}
\section{Introduction}  

In this paper we consider homogeneous CR-manifolds (CR=Cauchy-Riemann).
They are assumed to be of the form $M=G/H$ with 
$G$ being a connected Lie group acting on $M$ by CR-automorphisms.  
We present a number of Lie theoretic methods involving 
equivariant CR-fibrations and globalizations for studying these
manifolds.  These are applied to prove structure and classification
results in projective and low-codimensional settings.
It should be emphasized that,  although we only consider the setting where
a Lie group is acting transitively on our CR-manifold, its full group of
CR-automorphisms could be infinite-dimensional.

\bigskip\noindent 
An equivariant CR-fibration is just a CR-mapping $G/H\to G/I$ of
homogeneous CR-manifolds.  In our context the notion of globalization
is defined as follows. Since $M$ is real-analytic, it is embedded
(uniquely at the level of germs) as a generic submanifold of
a complex manifold $X$.  If $X$ is chosen sufficiently small,
then the Lie algebra $\widehat {\mathfrak g}:=\mathfrak g+i\mathfrak g$
acts locally and locally transitively on $X$. We say that the
$G$-action on $M$ can be globalized if there is 
a connected complex Lie group $\widehat G$ with Lie algebra  
$\widehat {\mathfrak g}$
such that $X$ can be taken to be a global complex homogeneous manifold
$\widehat G/\widehat H$, i.e., $M$ is just a $G$-orbit in such a manifold. 

\bigskip\noindent  
If $\widehat G$ is a complex Lie group and $\widehat H$ is a 
closed complex subgroup of $\widehat G$, then  
one lets $\widehat N := N_{\widehat G}(\widehat H^{\circ})$ 
be the normalizer in $\widehat G$ of the connected 
component of the identity of $\widehat H$ 
and then obtains an important tool, the normalizer fibration 
$\widehat G/\widehat H\to \widehat G/\widehat N$.  
This has proved very fruitful, since the analytic problems 
caused by the discrete isotropy of its fiber and the methods of 
algebraic groups that can be applied to its base 
(this is equivariantly embedded 
into some projective space) are now separated and can be analyzed 
somewhat independently.  

\bigskip\noindent  
In the case of a homogeneous CR-manifold $M=G/H$   
an analogue of the normalizer fibration, the $\mathfrak g$-anticanonical 
fibration, is known to exist. One may attempt to use this fibration,  
\[  
         G/H  \stackrel{F}{\longrightarrow} G/J  ,
\]   
to build a globalization of the given $M$. The base of this fibration
is globalizable, because it is a $G$-orbit in a projective space. 
Thus one must determine whether or not the parallelizable fiber
$F=J/H$ is globalizable and, if so, if its globalization can
be fit together with the globalization of the base.   

\bigskip\noindent
In the second section of this paper we recall the necessary basics 
about CR-manifolds,  
the construction and properties of the $\mathfrak g$-anticanonical 
fibration, and show the existence of the globalization $\widehat F$       
if $F$ has codimension at most two.  

\bigskip\noindent 
As indicated above, the base $G/J$ also has a globalization 
$\widehat G/\widehat J$ which is a $\widehat G$-orbit in projective space.   
The main point is to understand whether $\widehat J$ acts 
holomorphically on $\widehat F$.   
There is a surjectivity condition on the induced map of fundamental groups, 
which we call Condition (C),   
needed to deal with the possibility of ineffectivity of the 
$\widehat G$-action on the space $\widehat G/\widehat J$.  
If this condition is satisfied and if  
the radical $\widehat R$ of $\widehat G$ is acting as an Abelian 
group on the base $\widehat G/\widehat J$ (a necessary condition), 
then a globalization does exist.       
This result is proved in the general setting in Theorem 3.1 
and in the particular setting of the  $\mathfrak g$-anticanonical 
fibration in Theorem 3.2.  

\bigskip\noindent  
In the fourth section, using methods involving real and complex 
algebraic groups, we prove structure theorems for projectively 
embedded, compact, homogeneous CR-manifolds. 
We show that the complexified group 
$\widehat G$ is reductive, in particular 
that its radical $\widehat R$ is Abelian, and prove a splitting theorem
which can be regarded as the analogue of the Borel-Remmert theorem
for complex homogeneous spaces.
The results in this setting are general and should prove useful 
for other problems.    

\bigskip\noindent  
In the fifth section, methods of algebraic groups are applied
to give a detailed description of projective homogeneous
CR--manifolds of codimension at most two. In particular, 
the desired globalization results are obtained whenever the
2--dimensional affine quadric $\SL_2(\CC)/\CC^*$ is not
involved.  

\bigskip\noindent   
A fine classification of compact homogeneous CR-hypersurfaces 
with nondegenerate Levi forms is known.
In the last section we note that our classification 
results in section five give a classification of compact 
homogeneous CR-manifolds $M$ with nondegenerate 
Levi form in the codimension two setting.  
We also present remarks on the situation where $\widehat G/\widehat H$ 
is K\"{a}hler.

\bigskip\noindent
This work is dedicated to the memory of Anatolii Georgievich Vitushkin
with whom we had a continuous and extremely fruitful cooperation starting
in the mid-1980's.  Tolya was not only a mathematician of the highest
quality but also a wonderful teacher whose students are working on
a wide range of topics in the general area of complex analysis.
Complex and CR-manifolds having a high degree of symmetry appear
in many of their works. In our opinion a fusion of the methods 
which we present here and those developed by 
members of the \emph{Vitushkin School}, which are complementary 
to ours, will lead to a deeper understanding of the role of symmetry 
in complex analysis.

\bigskip\noindent
Let us close this introduction by giving a brief guideline to
the work in this direction of Vitushkin's students. 
A. Loboda (see, e.g., \cite {L}) has proved
classification results for locally homogeneous CR-manifolds 
in low-dimensions. V. Ezhov and V. Beloshapka and their
coworkers (in particular G. Schmalz who was also in the Moscow
seminar as a student of Shabat) have carried out computations
of Chern-Moser invariants of CR-manifolds which in fact have
a high degree of symmetry.  The classification theorems
of A. Isaev and N. Kruzhilin in situations where big
unitary groups are acting (see \cite{IK}) and the classification
work of Isaev et al for domains and hyperbolic manifolds with 
large automorphism groups relative to the dimension of the 
domain (see \cite{I} and its references) all involve settings 
where orbits of the real Lie group of interest are low-codimensional 
CR-manifolds which fit into the context of the present paper.

\bigskip\noindent
For other recent results in our context we refer to the 
work of G. Fels and W. Kaup (see e.g.\cite{FK}). For 
classification results in the hypersurface case see, e.g.,
\cite {MN}, \cite{N}, \cite {Ro}, \cite {AHR}, \cite {R}, \cite {HR} and \cite {AS}.

\bigskip\noindent
Finally, we wish to thank the referee for his critical remarks.
In our opinion they led to important improvements in this work.
\section{Basics on CR-manifolds}  

\subsection{CR-Structures}  

A CR-manifold is a real manifold that has some ``partial complex structure''.  
These manifolds have always been of interest,   
since they can occur in a number of natural ways.  
For example, any real-valued smooth function on a complex manifold defines 
a CR-submanifold at those points where its gradient does not vanish.  
Another important class of examples arises as the orbits of real subgroups 
in the homogeneous spaces of globalizations of those groups.   
 
\vskip 2ex\noindent{\bf Definitions:}\  
Suppose $M$ is an $n$-dimensional (differentiable) manifold.  
A {\bf CR-structure of type $(n,l)$} on $M$ is a subbundle $H$ of rank $l$ 
of the complexified tangent bundle $TM \otimes \mathbb C$ that satisfies:  
\begin{description}
\item  1) $H \cap \overline{H} = \{ 0 \}$, the zero section of  
$TM \otimes \mathbb C$  
\item 2) $H$ is involutive 
\end{description}  
A {\bf CR-manifold of type $(n,l)$} is a pair $(M,H)$ consisting 
of an $n$-dimensional manifold $M$ and a CR-structure 
$H$ of type $(n,l)$ on $M$. 

\vskip 2ex\noindent 
If $(M,H)$ is a CR-manifold of type $(n,l)$, then $H$ defines 
a subbundle $T^{CR}M$ of $TM$ of (real) rank $2l$, called 
the {\bf holomorphic tangent bundle} to $M$, that satisfies:
\begin{description}
\item 1') There is a bundle isomorphism  
${\cal J}: T^{CR}M \to T^{CR}M$ such 
that ${\cal J}^{2} = - {\rm Id}_{T^{CR}M}$ 
\item 2') For all $\widehat{\xi},\widehat{\zeta}\in\Gamma(M,T^{CR}M)$
one has $[\widehat{\xi},\widehat{\zeta}] \; - \; 
[{\cal J}\widehat{\xi},{\cal J}\widehat{\zeta}] \; \in \; \Gamma(M,T^{CR}M)$ and  
$N(\widehat{\xi},\widehat{\zeta}) := {\cal J}([\widehat{\xi},\widehat{\zeta}] -
[{\cal J}\widehat{\xi},{\cal J}\widehat{\zeta}])-
[\widehat{\xi},{\cal J}\widehat{\zeta}]-[{\cal J}\widehat{\xi},\widehat{\zeta}] \; = \; 0$ .  
\end{description}  
One has the equivalences 1) $\Longleftrightarrow$ 1') and 
2) $\Longleftrightarrow$ 2').  
As well, there is the bundle map $\alpha: T^{CR}M \to H$ 
defined by $\alpha(\widehat{\xi}) := \frac{1}{2} (\widehat{\xi} - i{\cal J}\widehat{\xi})$ and 
$\alpha$ is an isomorphism  
that satisfies $\alpha({\cal J}\widehat{\xi})=i\alpha(\widehat{\xi})$.    
Thus there is a one-to-one correspondence between 
CR-structures $H$ defined on $M$ and holomorphic 
tangent subbundles $T^{CR}M$ of $TM$.  


\vskip 2ex\noindent 
CR-manifolds can arise as submanifolds embedded in complex manifolds.  
Suppose $(X,\widetilde{{\cal J}})$ is a complex manifold. 
Then $\widetilde{{\cal J}}$ can be uniquely extended to a complex linear 
endomorphism of the complexified tangent space $TX \otimes \mathbb C$.  
Because $\widetilde{{\cal J}}^2 = - {\rm Id}_{ TX \otimes \mathbb C}$, 
the eigenvalues of $\widetilde{{\cal J}}$ are $+i$ and $-i$ and 
the corresponding eigenspaces of $TX\otimes\mathbb C$ are denoted, 
as is customary, by $TX^{1,0}$ and $TX^{0,1}$, respectively.  
Now assume $M \subset X$ is 
a (real) submanifold that satisfies the condition that    
\[  
   H_{x} \; := \; (T_{x}M \otimes \mathbb C) \; \cap \; T_{x}X^{1,0} 
\]  
{\em has constant rank for all}   $x\in M$.    
Then $(M,H)$ is a CR-manifold.  
One should note that the CR-structure on $M$ is also given by 
\begin{eqnarray*} 
     T^{CR}M & := & TM \; \cap \; \widetilde{{\cal J}} \ TM \\ 
     {\cal J} & := & \widetilde{{\cal J}}|_{T^{CR}M }  
\end{eqnarray*}  
with $T_{x}^{CR}M$ having constant rank for all $x\in M$.  
Note that for $x\in M$ we have the following decomposition of the 
tangent space  
\[  
    T_{x}(X) \; = \; T^{CR}_{x}(M) \; \oplus \; N_{x} \; \oplus \; {\cal J}(N_{x}) ,  
\]  
where 
\[  
     T_{x}(M) \; = \; T^{CR}_{x}(M) \; \oplus \; N_{x} ,   
\]   
and $N$ is a subbundle of $TM$ complementary to $T^{CR}(M)$.  

\noindent 
We say that $M$ has a {\bf generic CR-structure} if 
$TM + \widetilde{{\cal {\cal J}}}TM = TX|_{M}$.  
(The sum is not direct, in general.)


\vskip 2ex\noindent{\bf Remark:}\  
If the CR-structure on $M$ has type $(n,l)$ 
and $M$ is embedded in a complex manifold $X$ as a generic submanifold, 
then the (real) codimension of this 
structure is given by 
\begin{eqnarray*}  
                                 k & := & \dim_{\mathbb R}\ X \; - \; \dim_{\mathbb R} M \\  
                                 & = &  \dim_{\mathbb R} {\cal J}(N) \\  
                                 & = & \dim_{\mathbb R} N \\
                                 & = & \dim_{\mathbb R} M \; - \; \dim_{\mathbb R} T^{CR}M \\
                                 & = & n \; - \; 2l  ;  
\end{eqnarray*}  
this last formula also holds when $M$ is not embedded.  
It follows from this that 
\[  
            \dim_{\mathbb R} X \; = \; n \; + \; k \; = \; 2(n \; - \; l) ; \; \mbox{ i.e., } \; 
            \dim_{\mathbb C} X \; = \; n \; - \; l .  
\]              

\vskip 2ex\noindent{\bf Definition:}\  
Suppose $(M_{1}, T^{CR}(M_{1}),{\cal J}_{1})$ and 
$(M_{2}, T^{CR}(M_{2}),{\cal J}_{2})$ are two CR-manifolds.  
A smooth map $f:M_{1} \to M_{2}$ is called a CR-{\bf map} if 
\[  
    f_{*}(T_{x}^{CR}(M_{1}) )\; \subset \;  T_{f(x)}^{CR}(M_{2})
\]  
for all $x\in M_{1}$ and 
\[  
      {\cal J}_{2} \circ f_{*} \; = \; f_{*} \circ {\cal J}_{1} 
\] 
Concepts such as CR-submersion, CR-isomorphism, and CR-function 
are obvious.  

\vskip 2ex\noindent  {\bf Definition:}\  
A vector field $\widehat{\xi}\in\Gamma(M,TM)$ is called a {\bf CR-vector field} if 
the local one-parameter group of transformation of $M$ that $\widehat{\xi}$ induces 
consists of CR-transformations.   
We will denote the set of CR-vector fields on $M$ by $\Gamma_{CR}(M,TM)$.  

\vskip 2ex\noindent 
A CR-manifold $(M,H)$ is called an {\bf analytic CR-manifold} 
if $M$ is an analytic manifold and $H$ is an analytic subbundle of 
$TM \otimes \mathbb C$; i.e., $H$ is locally generated by analytic 
local sections of $TM \otimes \mathbb C$.  

\vskip 2ex\noindent 
Given a CR-manifold, one would be in an ideal situation if one could realize this manifold 
as an embedded manifold in a complex manifold from which it inherits its CR-structure.  

\vskip 2ex\noindent 
An embedding $\sigma: M \to X$ of a CR-manifold into a complex manifold 
$X$ is called a CR-embedding if $\sigma(M) \subset X$ is a CR-submanifold 
and the map $\sigma$ is a CR-isomorphism of $M$ onto its image.  
Such an embedding is called {\bf generic} if its image $\sigma(M)$ is a generic 
CR-submanifold of $X$.  
We will then say that $(X,\sigma)$ is a {\bf complexification} of $M$. 

\begin{theorem}[Andreotti-Friedrichs \cite{AF}] 
Every analytic CR-manifold $M$ has a complexification $(X,\sigma)$ that is unique 
up to some biholomorphic map.  
\end{theorem}

\begin{theorem} \label{extension}  
\begin{enumerate}  
\item Suppose $f:M_{1} \to M_{2}$ is an analytic CR-map between 
two generic analytic  
CR-submanifolds $M_{i} \subset X_{i}, i=1,2$, where 
the $X_{i}$ are complex manifolds.  
Then there exist open neighborhoods $U_{i} \subset X_{i}$ of $M_{i}$ 
for $i=1,2$ and a holomorphic map $\widehat {f}:U_{1} \to U_{2}$ such that 
$\widehat {f}|_{M_{1}}= f$.  
\item A holomorphic map $\widehat {f}$ defined on an open, connected neighborhood 
of a generically embedded CR-manifold is constant if and only if 
$\widehat {f}|_{M}$ is constant.  
\item Suppose $M \subset X$ is a generic CR-submanifold of a complex 
manifold $X$ and $\widehat{\xi}\in\Gamma_{CR}(M,TM)$ is a CR-vector field on $M$.  
Then there exist an open neighborhood $U$ of $M$ in $X$ and a holomorphic 
vector field $\widehat{\zeta}$ on $U$ such that 
\[  
             \widehat{\xi}(x) \; = \; \widehat{\zeta}(x) \; + \; 
             \widehat{\zeta}^{\dagger} (x)  
\]  
for every $x\in M$, where  $\widehat{\zeta}^{\dagger}$ denotes the 
complex conjugate of the vector field $\widehat{\zeta}$.  
\end{enumerate}  
\end{theorem}   
 
\vskip 2ex\noindent  {\bf Definition:}\  
Let $M$ and $B$ be CR-manifolds.   
Suppose $(M,\pi,B,G)$ is a principal bundle, 
where  $\pi:M \to B$ is the bundle projection and $G$ is 
the structure group.  
If the right action of $G$ on $M$ is by CR-transformations, 
then we call  $(M,\pi,B,G)$ a {\bf CR-principal bundle}.  


\subsection{The Levi Form}  

Suppose $(M,H)$ is a CR-manifold and let $\pi: TM \otimes {\mathbb C}  \to 
TM\otimes{\mathbb C}/H\oplus\overline{H}$ denote the projection.  
Then $\pi$ induces an ${\mathbb R}$-linear bundle map 
\[  
    L(M,H): H \times H \; \longrightarrow \; TM \otimes{\mathbb C}/H\oplus\overline{H} 
\]  
that is defined in the following way:  
For $x\in M$ and $(a,b)\in H_{x}\times H_{x}$ 
choose local sections $\widehat{\xi}$ and $\widehat{\zeta}$ in $H$ with 
$\widehat{\xi}(x)=a$ and $\widehat{\zeta}(x)=b$. 
Set 
\[  
     L(M,H)(a,b) \; := \; \pi_{x}[\widehat{\xi},\overline{\widehat{\zeta}}](x).   
\]   
The map $L(M,H)$ is well defined and is called the {\bf Levi form}  
of $(M,H)$.  
Some pertinent facts about the Levi form are:
\begin{enumerate}
\item The map $L$ is antihermitian. 
\item The set $L^{\circ} := \{ \widehat{\xi}\in H \ | \ L(\widehat{\xi},\widehat{\zeta}) =0 \; \forall \widehat{\zeta}\in H \}$ 
is called the {\bf Levi kernel} of $M$.  
We say that $L$ is {\bf non-degenerate} if $L^{\circ} = \{ 0 \}$.  
\item If $M$ is a real hypersurface in a complex manifold $X$, 
then $iL$ has the same signature as the restriction of the 
complex Hessian matrix of a defining function of $M$ to 
$H=(TM\otimes{\mathbb C})\cap TX^{1,0}$.  
The latter is the classical Levi form; see 11. Theorem, p. 290 ff in \cite{Gr}.
\item A CR-manifold is called {\bf Levi flat} if 
$L^{\circ}=H$.  
\end{enumerate}  


\subsection{Homogeneous CR-structures}  

\vskip 2ex\noindent{\bf Definition:}\  
A CR-{\bf automorphism} of a given CR-manifold is a CR-map 
of $M$ to itself that is also a diffeomorphism.  
If the group of CR-automorphisms acts transitively on $M$, then 
$M$ is called a {\bf homogeneous} CR-{\bf manifold}.  

\begin{theorem}[\cite{R}, Satz 1.3.1.2]\label{Satz1.3.1.2}
Let $M=G/H$ be a homogeneous CR-manifold.  
The $G$-invariant CR-structures $(R,{\cal J})$ on $M$ 
are in one-to-one correspondence with the pairs  
$(\widetilde{R},\widetilde{{\cal J}})$, where $\widetilde{R}$ is a vector subspace of $\mathfrak g$ 
with ${\mathfrak h} \subset \widetilde{R} \subset {\mathfrak g}$ 
and $\widetilde{{\cal J}}: \widetilde{R} \to \widetilde{R}$ is an endomorphism,     
that satisfy the following:  
\begin{description}  
\item 1) $\widetilde{{\cal J}}\xi=0$ if and only if $\xi\in{\mathfrak h}$  
\item 2) $\widetilde{{\cal J}}^{2} \xi + \xi \in {\mathfrak h}$  
for all $\xi\in\widetilde{R}$  
\item 3) ${\rm Ad}_{g}\xi\in\widetilde{R}$ and 
${\rm Ad}_{g}\widetilde{{\cal J}}\xi-\widetilde{{\cal J}}{\rm Ad}_{g}\xi\in {\mathfrak h}$  
for all $g\in H$ and $\xi\in\widetilde{R}$  
\item 4) $[\xi,\zeta] - [\widetilde{{\cal J}}\xi,\widetilde{{\cal J}}\zeta]\in\widetilde{R}$ and   
$\widetilde{{\cal J}}([\xi,\zeta] - [\widetilde{{\cal J}}\xi,\widetilde{{\cal J}}\zeta]) - 
[\widetilde{{\cal {\cal J}}}\xi,\zeta] - 
[\xi,\widetilde{{\cal J}}\zeta] \in {\mathfrak h}$ for all $\xi,\zeta\in\widetilde{R}$.  
\end{description}  
Two pairs $(\widetilde{R},\widetilde{{\cal J}})$ and $(\widetilde{R}',\widetilde{{\cal J}}')$  
are equivalent if and only if  $\widetilde{R}=\widetilde{R}'$ and  
$\widetilde{{\cal J}}\xi-\widetilde{{\cal J}}'\xi\in{\mathfrak h}$ f
or all $\xi\in\widetilde{R}=\widetilde{R}'$.  
If $H$ is connected, then 3) may be replaced by  
\[  
  3')\; [\xi,\zeta]\in\widetilde{R} \mbox{  and  }   
  \widetilde{{\cal J}}[\xi,\zeta] \ - \ [\xi,\widetilde{{\cal J}}\zeta] \in {\mathfrak h}  
\]  
for all $\xi\in{\mathfrak h}$ and $\zeta\in\widetilde{R}$.  
\end{theorem} 
 
\noindent  
The proof is given in \cite{R}, pp. 18 - 21, basically following  
\cite{KN2} except $\widetilde{R}\not={\mathfrak g}$.  
  
\begin{corollary}  [Richthofer \cite{R}, Zusatz to Satz 2, p. 21]  
The $G$-invariant CR-structures on $M$ are analytic.  
\end{corollary}

\vskip 2ex\noindent  {\bf Consequence: \  
Every homogeneous CR-manifold $M=G/H$ has a complexification.}    


\subsection{Complexification of the Lie algebra} \label{cplxfn} 

Let $M=G/H$ be a homogeneous CR-manifold.  
Since the invariant CR-structure $(R,{\cal J})$ on $G/H$ is analytic, 
a complexification $(\widehat {M},\sigma)$ of $G/H$ exists.  
We will think of $M$ as already embedded in $\widehat {M}$, i.e., $M=\sigma(G/H)$.  
The action map $\lambda:G\times M \to M, (g, \sigma(aH)) \mapsto \sigma(gaH)$ 
induces a homomorphism $\phi:G \to {\rm Aut}_{CR}(M)$, the group of 
CR-automorphisms of $M$ and a homomorphism 
\[  
    \alpha:{\mathfrak g} \; \longrightarrow \; \Gamma_{CR}(M,TM) .  
\]   
(Modulo ineffectivity, we are identifying $\mathfrak g$, the Lie algebra of $G$, 
with a subalgebra of the algebra of CR-vector fields on $M$.)  
For any $\xi\in\mathfrak g$, let $(\widehat{\xi}_{t})_{t\in\mathbb R}$ be the 
one parameter group of the vector field $\alpha(X)$.  
Obviously, $\widehat{\xi}_{t}=\phi(\exp \ t \xi)$.  
Therefore, there exist a neighborhood $U$ of $M$ in $\widehat {M}$ and a 
holomorphic vector field $\widehat{\zeta}$ on $U$ such that 
$\widehat{\zeta}+\overline{\widehat{\zeta}}|_{M} 
= \widehat{\xi} $, see  Theorem \ref{extension}.  
Since $\alpha(\mathfrak g)$ is finite dimensional, one sees that by 
choosing a basis for this Lie algebra one can find an open neighborhood $U$ 
of $M$ such that every vector field $\widehat{\xi}\in\alpha(\mathfrak g)$ can be 
so ``extended'' to a holomorphic vector field $\widetilde{\xi}$ on $U$.  
Moreover, the map $\xi \mapsto \widetilde{\xi}$ is a Lie algebra homomorphism.     
This construction allows us to define the complexification $\widehat {\mathfrak g}$ 
of the Lie algebra $\mathfrak g$ to be 
\[   
     \widehat {\mathfrak g} \; := \; \{ \ \widetilde{\xi} \; + \; 
     {\cal J}\widetilde{\zeta} \ | \ \xi,\zeta \in {\mathfrak g} \ \} ;    
\]  
namely, $\widehat {\mathfrak g}$ is the complex Lie subalgebra of 
$\Gamma(U, T^{1,0}U)$ that is generated by the image of $\mathfrak g$.  
One can easily check that $\mathfrak m := \mathfrak g \cap {\cal J} \mathfrak g$ 
is an ideal in the Lie algebra $\mathfrak g$.  

\subsection{The $G$-action on $\widehat {\mathfrak g}$}\label{gaction}   

\vskip 2ex\noindent  
We would now like to define a $G$-action on 
the complex Lie algebra $\widehat {\mathfrak g}$.    
In order to do this, we consider the action 
\[ 
      G \times \alpha(\mathfrak g) \; \longrightarrow \; \alpha(\mathfrak g)  
\]  
given by 
\[  
      g . \alpha(\xi) \; = \; \alpha({\rm Ad}_{g} \ \xi) 
\]  
or 
\[  
     g. \alpha(\xi)(p) \; = \; d\lambda_{g}(\alpha(\xi)(g^{-1}.p)) , \mbox{ where } 
p\in M.  
\]  
If $\xi\in{\mathfrak m}$, then ${\cal J}\alpha(\xi)$ is tangential to $M$ and thus 
$g. {\cal J} \alpha(\xi) = {\cal J} g.\alpha(\xi)$, since $G$ is acting as a group of 
CR-transformations on $M$.  
One extends this action to $\widehat {\mathfrak g}$ in the following way.  
Let $\eta = \widetilde{\xi} + {\cal J}\widetilde{\zeta}\in\widehat {\mathfrak g}$.  
Define $g\eta = \widetilde{g\xi} + {\cal J}\widetilde{g\zeta}$, where 
$\widetilde{g\xi}$ is the extension 
of $g.\alpha(\xi)$ and $\widetilde{g\zeta}$ is the extension of $g.\alpha(\zeta)$.  
Because of the complex linear $G$-action on $\alpha(\mathfrak m)$, this 
extension is well-defined and the map 
$\eta \mapsto g.\eta$ of $\widehat {\mathfrak g} \to \widehat {\mathfrak g}$ is 
{\bf complex linear}, i.e., $g.{\cal J}\eta = {\cal J}g.\eta$.  

\vskip 2ex\noindent  
A naturally occurring situation is the following.  
Suppose $\widehat {G}$ is a complex Lie group and $\widehat {H}$ is a closed 
complex subgroup of $\widehat {G}$.  
Further, let $G$ be a real subgroup of $\widehat {G}$ such that  
$G/H \hookrightarrow \widehat {G}/\widehat {H}$ is a generic CR-submanifold, 
where $H =\widehat {H}\cap G$.   
Throughout this paper we will call the complex manifold 
$\widehat {M} := \widehat {G}/\widehat {H}$ a 
{\bf globalization} of $M$, since the local action of $\widehat G$ 
in a neighborhood of $M$ has been globalized.   
Under the assumption that the actions are effective it follows that the Lie algebra 
$\widehat {\mathfrak g}$ of $\widehat {G}$ is isomorphic to the complexification 
of the Lie algebra $\mathfrak g$ of $G$ that was discussed above.  

\subsection{The ${\mathfrak g}$-anticanonical fibration}  

There are two main approaches to the definition of the  
${\mathfrak g}$-anticanonical fibration.  
We will present both, beginning with a more algebraic 
method and then later using the complex tools at hand 
due to the embedding.  
Most of what follows can be found in \cite{R}, \cite{HR}, 
and \cite{AHR}.  

\subsubsection{Equivariant CR-Fibrations}\label{defncr}  

If there exists a closed subgroup $I \subset G$ with $H \subset I$ 
such that $G/I$ has a  $G$-invariant CR-structure 
and the projection $\pi: G/H \to G/I$ is a CR-mapping, then we will say 
that we have a CR-fibration.  
A fundamental question is whether such fibrations exist and what 
information they tell us about the homogeneous CR-manifold.  
Suppose the CR-structure on $M=G/H$ is given by a pair $(R,{\cal J})$ and   
the corresponding structure on $\mathfrak g$ is given 
by $(\widetilde{R},\widetilde{{\cal J}})$, see Theorem \ref{Satz1.3.1.2}.  
Set 
\[   
      N(H^{\circ}) := \{ \ g\in G \ | \   gH^{\circ}g^{-1} \subset H^{\circ} \ \}  .    
\]   
Then it is clear that $H\subset N(H^{\circ})$ and at the Lie algebra 
level one has $N(H^{\circ}) := \{ \ g\in G \ | \ {\rm Ad}_{g}{\mathfrak h} 
\subset {\mathfrak h} \ \}$.  
In order to take into account the CR-structure on $G/H$ as well, we consider   
the CR-{\bf normalizer} $N_{CR}(H)$ of $H$ which is defined as follows:   
\[    
   N_{CR}(H) \; := \; \{ \ g\in G \ | \ {\rm Ad}_{g} \xi\in\widetilde{R} \mbox{ and }     
                   {\rm Ad}_{g} \widetilde{{\cal J}}\xi  -\widetilde{{\cal J}}{\rm Ad}_{g}   
                 \xi\in{\mathfrak h}, \quad \forall \xi\in\widetilde{R} \ \}.  
\]   
Note that $N_{CR}(H) = N_{CR}(H^{\circ})$ and the definition of    
$N_{CR}(H)$ only depends on the equivalence class of   
$(\widetilde{R},\widetilde{{\cal J}})$, see Theorem \ref{Satz1.3.1.2}.  
By the same result $(\widetilde{R},\widetilde{{\cal J}})$ also defines a CR-structure 
on $G/H^{\circ}$.  
Now one can define a real analytic right action  
\[  
     r: G/H^{\circ} \times N(H^{\circ}) \; \longrightarrow \; G/H^{\circ}  
\]  
by $r(gH^{\circ},n) := r_{n}(gH^{\circ}) := gnH^{\circ}$  for   
$g\in G$ and $n\in N(H^{\circ})$.  
In general, the mapping $r_{n}$ is not a CR-mapping for every   
$n\in N(H^{\circ})$.   
However, one has the following characterization of the subgroup $N_{CR}(H)$.  

\begin{proposition}[\cite{R} p. 23] \label{CRnorm}  
One has $H \subset N_{CR}(H) \subset N(H^{\circ})$.  
The group $N_{CR}(H)$ consists of all $n\in N_{G}(H^{\circ})$ such that
the right translations $r_{n}$ given by 
\[  
      r_{n}: G/H^{\circ} \; \longrightarrow \; G/H^{\circ} \quad\quad gH^{\circ} 
      \; \mapsto \; gnH^{\circ}   
\] 
are CR-mappings.  
\end{proposition}

\vskip 2ex\noindent{\bf Definition:} 
The map $\varphi_{\widehat {\mathfrak g}} :G/H\to G/N_{CR}(H)$ is called the 
{\bf ${\mathfrak g}$-anticanonical fibration} of $G/H$.   

\medskip\noindent{\bf Remark:}\ 
We note below that this map really is a CR-map.  
Also it is important to keep in mind that another construction of 
the {\bf ${\mathfrak g}$-anticanonical fibration} is given in the 
next subsection, along with further details about some of its properties.  
The notation $\varphi_{\widehat {\mathfrak g}}$ is also justified by 
that construction, see \S \ref{construction}.    

\vskip 2ex\noindent 
An important consequence of the above that will be extremely useful later on 
is the following.  

\begin{corollary}
If the ${\mathfrak g}$-anticanonical fibration is degenerate, i.e., 
$N_{CR}(H) = G$, and the $G$-action is assume to be almost effective, 
then $H$ is discrete.  
\end{corollary}  


\vskip 2ex\noindent
Suppose that $M$ is embedded in a complexification $\widehat {M}$.
Let $\widehat {\mathfrak h} := \{ \zeta\in\widehat {\mathfrak g} \ | \ \zeta(o)=0\ \}$,
where $o=\pi(e)$ and $\pi:G\to G/H$.
Recall
\[
     \widehat {\mathfrak m} \; = \; {\mathfrak g} \ + \ {\cal J}{\mathfrak g}
\]
where ${\cal J}$ is the complex structure on $\widehat {M}$ and $\mathfrak g$
is regarded as a real subalgebra of $\Gamma_{\mathcal O}(\widehat {M},T\widehat {M})$.
It turns out that
\[
     N_{CR}(H) \; = \; \{ \ g\in G \ | \ {\rm Ad}_{g} \widehat {\mathfrak h} =
                       \widehat {\mathfrak h} \ \} ,
\]
where ${\rm Ad}:G\to {\rm Aut}(\widehat {\mathfrak g})$ is the adjoint action
of $G$ on $\widehat {\mathfrak g}$ that is induced by the usual adjoint action
of $G$ on $\mathfrak g$.
This shows that $G/N_{CR}(H)$ inherits a $G$-invariant CR-structure from
the Grassmannian defined by the complex vector subspaces of 
$\widehat {\mathfrak g}$
that have the same dimension as $\widehat {\mathfrak h}$ and $N_{CR}(H)$ is
the isotropy of the $G$-adjoint action computed at the point that corresponds
to the subspace $\widehat {\mathfrak h}$.
The Pl\"{u}cker embedding then gives us a $G$-equivariant map
\[
     G/H \; \longrightarrow \; G/N_{CR}(H) \; \longrightarrow \; 
       {\mathbb P}_{k} .
\]
This map is the same as the one given by the holomorphic sections of the
anticanonical bundle of $\widehat {M}$ that are generated by the
$\widehat {\mathfrak g}$-sections.
This fact is the reason for the name ${\mathfrak g}$-anticanonical
fibration.  

\vskip 2ex\noindent  
We have a representation
of $G$ into the group $GL_{k+1}(\mathbb C)$.  
But we only use the action on the projective space, so we may 
regard $G$ as being mapped into $PGL_{k+1}(\mathbb C)$.  
Let $\widetilde{G}$ denote the smallest complex Lie subgroup that contains
the image of $G$ and let $\widetilde{N}_{CR}(H)$ denote the isotropy subgroup of
the $\widetilde{G}$-action computed at the point corresponding to $N_{CR}(H)$.

\vskip 2ex\noindent
Since $N_{CR}(H)\subset N_{G}(H^{\circ})$, the fiber of the
$\mathfrak g$-anticanonical fibration may be written as
a quotient $A/\Gamma$, where $A := N_{CR}(H)/H^{\circ}$ is a Lie
group and $\Gamma := H/H^{\circ}$ is a discrete subgroup of $A$.
We will make some observations about this particular situation in 2.7.

\vskip 2ex\noindent
The question whether $G/N_{CR}(H)$ admits a  CR-structure such that
the map $G/H \to G/N_{CR}(H)$ is a CR-submersion is answered by the 
following.
Note that this result is slightly more general.

\begin{theorem}[\cite{R}, Satz 1.3.2.2]\label{Satz2}
Let $L\subset N_{CR}(H)$ be a closed subgroup such that $H\subset L$ and
$\widetilde{{\cal J}}({\mathfrak l}\cap\widetilde{R})\subset({\mathfrak l}\cap\widetilde{R})$,
where ${\mathfrak l}$ denotes the Lie algebra of $L$.
Then $G/L$ has a $G$-invariant CR-structure so that
the map $G/H \to G/L$ is a CR-submersion.
\end{theorem}  

\begin{corollary}
The manifold $G/N_{CR}(H)$ has a $G$-invariant CR-structure
so that the map $G/H\to G/N_{CR}(H)$ is a CR-submersion.
\end{corollary}


\vskip 2ex\noindent  
In order to give a further description of the fibration 
$G/H \rightarrow G/N_{CR}(H)$ 
we now present some observations about the Levi form of invariant 
CR-structures $(R,{\cal J})$ on $M=G/H$.  

\vskip 2ex\noindent 
Again let $\pi:G \rightarrow G/H$ be the quotient map and set $o=\pi(e)$.  
Let $\psi:{\mathfrak g} \rightarrow {\mathfrak g}/\widetilde{R}$ be the 
projection.  
Set $\widetilde{L}(\xi,\zeta) := \psi[\xi,\zeta]$ for $\xi,\zeta\in\widetilde{R}$.  
Then 
\[  
    \widetilde{L} \; = \; \pi^{*} L_{R}|_{\widetilde{R}} .  
\]   
This follows from the fact that $P(\widetilde{R})$ generates $\widetilde{R}_{g}$ at 
every point $g\in G$ and because $\widetilde{R}$ is a left invariant subbundle of $TG$.  
(For the notation, see Theorem \ref{Satz1.3.1.2}.)  

\begin{theorem}[\cite{R}, Satz 7]
If $H$ is connected, then one has a CR-principal bundle 
\[  
    G/H \; \longrightarrow \;  G/N_{CR}(H) .  
\]      
\end{theorem}

\begin{lemma}[Lemma 1, \cite{AHR}]
If $\zeta\in\widehat {\mathfrak g}$ vanishes at a point $p\in M$, then $\zeta$ 
vanishes identically on $F_{p}$.
\end{lemma}  

\begin{corollary}[\cite{AHR}]  
The fibers of $\varphi_{\widehat {\mathfrak g}}$ are Levi flat CR-submanifolds of 
$M$.  
If the image of $\varphi_{\widehat {\mathfrak g}}$ is a point, then 
\[  
    {\mathfrak m} \; = \; \{ \ \xi\in{\mathfrak g} \ | \ {\cal J}\widehat{\xi}\in \Gamma(M,TM) \ \} 
\]  
generates $L^{\circ}(M)$.  
If $M$ is a CR-hypersurface, then $\dim_{\mathbb C}{\mathfrak m} = n-1$, 
where $n := \dim_{\mathbb C}\widehat {M}$.
\end{corollary}

\vskip 2ex\noindent 
If the image of $\varphi_{\widehat {\mathfrak g}}$ is a point, then $F_{p}=M$.  
Thus ${\mathfrak m}_{p}F :=  {\mathfrak m}_{p}= {\mathfrak m}$.  
Since the $G$-action is almost effective, this observation suffices 
to complete the proof.  
     
\begin{corollary}[See \cite{AHR}]
Let $(M,G)$ be a homogeneous CR-manifold.   
Set $k := {\rm rank}_{\mathbb C} L^{\circ}(M)$.  
Let $F$ be a fiber of the ${\mathfrak g}$-anticanonical fibration of $M$.  
Then $F$ is a Levi flat CR-submanifold of $M$ of type $(\alpha,\beta)$, 
where $\beta \le k$. 
The ideal ${\mathfrak m}_{{\mathfrak j}/{\mathfrak h}}$ coming from the right 
action 
of $J/H^{\circ}$ on $G/H^{\circ}$ generates a complex subbundle of 
$L^{\circ}(M)$ 
of rank $\beta$.  
\end{corollary}       


\subsubsection{Another Construction and Further Properties}\label{construction} 

Throughout we assume $G$ is connected and is acting almost effectively 
on $M$.  
Then the map $\alpha:{\mathfrak g}\to\Gamma(M,TM)$ 
that is induced by the action $\lambda:G\times M\to M$ 
is injective.     
For $g\in G$, we write $\lambda_{g}:M\to M$ for the 
map $\lambda_{g}(x) := \lambda(g,x)$ and let $g\cdot x$ 
denote the point $\lambda(g,x)$.  

\vskip 2ex\noindent 
Let $(\widehat {M},\sigma)$ be a ${\mathfrak g}$-complexification 
of $M$ (see \S \ref{cplxfn}) and $n = \dim_{\mathbb C} \widehat {M}$.  
In the following we will think of $M$ as already imbedded 
in $\widehat {M}$, i.e., $\sigma|_{M}={\rm id}_{M}$.  
Set 
\[  
    V_{\mathfrak g} \; := \; \bigwedge_{n} \; \widehat {\mathfrak g} .  
\]  
This is a $(k+1)$-dimensional complex vector subspace of the 
space of holomorphic sections of the anticanonical bundle 
of $\widehat {M}$:  
\[  
       V_{\mathfrak g} \; \subset \;  \Gamma_{\mathcal O} 
      (\widehat {M},\bigwedge_{n} T\widehat {M}^{1,0}) .  
\]  
The natural $G$-action on $\widehat {\mathfrak g}$ (see \S \ref{gaction}) defines 
a $G$-action on ${\mathbb P}(V_{\mathfrak g}^{*}) \cong {\mathbb P}_{k}$ 
given by $g\cdot f(\sigma) = f(g^{-1}\cdot\sigma)$ for 
$g\in G, f\in {\mathbb P}(V^{*}_{\mathfrak g}) , \sigma\in V_{\mathfrak g}$.  
The map ${\mathbb P}(V^{*}_{\mathfrak g})\to 
{\mathbb P}(V^{*}_{\mathfrak g}), 
f \mapsto g\cdot f$ is projective linear for every $g\in G$; see \S \ref{gaction}.
Now let $U\subset\widehat {M}$ be an open neighborhood of $M$, so that for every 
$x\in U$ there exists a $\sigma\in V_{\mathfrak g}$ with $\sigma(x)\not= 0$, 
and let 
\[  
     \widehat {\varphi_{\widehat {\mathfrak g}} } : U \; \longrightarrow\; 
      {\mathbb P}(V^{*}_{\mathfrak g}) 
\]  
be given by $\widehat {\varphi_{\widehat {\mathfrak g}} }(x)(\sigma)=\sigma(x)$.  
Since $\bigwedge_{n} T\widehat {M}^{1,0}$ is a line bundle and for every $x\in U$ 
there exists a $\sigma\in V_{\mathfrak g}$ with $\sigma(x)\not= 0$, the 
map $\widehat {\varphi_{\widehat {\mathfrak g}} }$ is well defined.  
In order to show that $\widehat {\varphi_{\widehat {\mathfrak g}} }$ is also holomorphic, 
we pick a basis $\{ \sigma_{0},\ldots,\sigma_{k} \}$ of $V_{\mathfrak g}$ and set 
$\sigma = \sum_{i=0}^{k} a_{i}(\sigma) \sigma_{i}$.  
Then 
\[  
    \widehat {\varphi_{\widehat {\mathfrak g}} }(x)(\sigma) \; = \; \sigma(x) \; = \; 
    \sum_{i=0}^{k} a_{i}(\sigma) \sigma_{i}(x) .   
\]   
Let $\{ \sigma_{0}^{*},\ldots,\sigma_{k}^{*} \}$ be the dual basis of 
$V^{*}_{\mathfrak g}$ and let $[z_{0}: \ldots :z_{k}]$ be the homogeneous 
coordinates on ${\mathbb P}(V^{*}_{\mathfrak g})$ that are defined by  
$\{ \sigma_{0}^{*},\ldots,\sigma_{k}^{*} \}$.  
Then any $f\in{\mathbb P}(V^{*}_{\mathfrak g})$ can be written as 
$f =  [z_{0}(f): \ldots :z_{k}(f)]$, where $z_{j}(f) = f(\sigma_{j})$, 
for $j=0, \ldots, k$.  
Then 
\[  
    z_{j}(\widehat {\varphi_{\widehat {\mathfrak g}} }(x)) \; = \; 
    \widehat {\varphi_{\widehat {\mathfrak g}} }(x)(\sigma_{j}) \; = \; 
              \sigma_{j}(x)  
\]   
and thus 
\[  
    \widehat {\varphi_{\widehat {\mathfrak g}} }(x) \; = \;   [\sigma_{0}(x): \ldots :\sigma_{k}(x)]  .  
\]  
Since the $\sigma_{j}\in \Gamma_{\mathcal O}(\widehat {M},\bigwedge_{n} 
T\widehat {M}^{1,0})$, 
by the choice of $U$ one sees that $\widehat {\varphi_{\widehat {\mathfrak g}} }$ 
is holomorphic and well defined.  

\vskip 2ex\noindent  
For $x\in\widehat {\varphi_{\widehat {\mathfrak g}} }(M)$, one has 
\[  
     g\cdot \widehat {\varphi_{\widehat {\mathfrak g}} }(x)(\sigma) \; = \; 
     \widehat {\varphi_{\widehat {\mathfrak g}} }(x)(g^{-1}\cdot\sigma) 
             \; = \; [\sigma(g\cdot x)] .  
\]  
Thus in the coordinates chosen above 
\[   
     g\cdot \widehat {\varphi_{\widehat {\mathfrak g}} }(x) \; = \; 
          [ \sigma_{0}(g\cdot x) : \ldots : \sigma_{k}(g\cdot x) ] .  
\]   
This implies that the map $\varphi_{\widehat {\mathfrak g}}  := 
\widehat {\varphi_{\widehat {\mathfrak g}} }|_{M}$ is $G$-equivariant.  
Set 
\[  
    J \; = \; \{\; g\in G \; : \; g\cdot\varphi_{\widehat {\mathfrak g}} \circ\pi(e) \; = \;  
                          \varphi_{\widehat {\mathfrak g}} \circ\pi(e) \; \}  ,  
\]  
where $e$ is the identity in $G$ and $\pi:G\to G/H$ 
is the quotient map.  
Then $G/J$ is a (not necessarily closed) CR-submanifold of 
${\mathbb P}_{k}$.


\begin{theorem}
Let  $G/H \to G/J$ be the fibration constructed above.  
Then $J = N_{CR}(H)$ (see \S \ref{defncr}) and if $G/N_{CR}(H)$ is endowed 
with its $G$-invariant CR-structure, then the map 
$G/N_{CR}(H)\to G/J$ is CR.  
In particular, $G/H\to G/J$ is the ${\mathfrak g}$-anticanonical fibration of $G/H$. 
\end{theorem}  

\begin{proof}  
First we note that 
\[  
        J \; \subset \; N(H^{\circ})  .  
\]  
Otherwise, there would be an element $\xi$ in the Lie algebra 
${\mathfrak h}$ of $H$ and an element $g\in J$ such that 
${\rm Ad}_{g}\xi \not\in{\mathfrak h}$.  
But then 
\[  
   \eta(\pi(e)) \; := \; \frac{1}{2} (\widetilde{\xi} - i{\cal J}\widetilde{\zeta})(\pi(e)) \; = \; 
     0 \; \not= \; \eta(\pi(g)) , 
\]  
which contradicts the previous Lemma.  

\vskip 2ex\noindent 
Since $(\widetilde{R},\widetilde{{\cal J}})$ (by the definition of $N_{CR}(H)$ in 
\S \ref{defncr}) is defined by ${\mathfrak h}$ and $(R,{\cal J})$ only, one may 
consider the universal covering $G/H^{\circ}$ of $G/H$ with 
its lifted CR-structure and by Proposition \ref{CRnorm}  one has to show that $J$ consists 
of all elements $g\in N(H^{\circ})$ such that the map $r_{g} : G/H^{\circ} 
\to G/H^{\circ}, a\cdot H^{\circ} \mapsto a\cdot gH^{\circ}$ is CR.  

\vskip 2ex\noindent  
Let $g\in J$ and for $x\in G/H^{\circ}$ set 
${\mathfrak m}_{x} = \{ \widehat{\xi}\in\alpha({\mathfrak g}) : \xi(x)\in R_{x} \}$.  
Let $\widehat{\zeta}$ and $\widehat{\zeta'}$ be the holomorphic vector fields 
on $U$ corresponding to $\widehat{\xi}$ and $\widehat{\xi'}$.  
Then $(\widehat{\zeta'} - i\widehat{\zeta})\in\widehat {\mathfrak g}$ has a zero at the point $x$.  
By the preceding Lemma, this vector field vanishes along $F_{x}$.  
This means that ${\cal J}\widehat{\xi}=\widehat{\xi'}$ on $F_{x}$.  
Thus ${\mathfrak m}_{x}$ generates the space $R_{p}$ for every $p\in F_{x}$.  
Moreover, for $g\in J$ one has  
\[  
   dr_{g}({\cal J}\widehat{\xi}(x)) \; = \; dr_{g}\widehat{\xi'}(x) \; = \; \widehat{\xi}'(r_{g}(x)) \; = \;  
      {\cal J}\widehat{\xi}(r_{g}(x)) \; = \; {\cal J}dr_{g}\widehat{\xi}(x))  ,  
\]  
and so $r_{g}$ is a CR-map for every $g\in J$.  

\vskip 2ex\noindent  
Conversely, if $g\in N_{CR}(H^{\circ})$, then there exist neighborhoods 
$V$ and $W$ of $M$ in $U$ and an extension of $r_{g}$ to a biholomorphic 
map $\widetilde{r}_{g}: V \to W$.  
If $\eta := \widetilde{\xi}+{\cal J}\widetilde{\zeta}\in \widehat {\mathfrak g}$, 
then on $W$ one has  
\[  
   \widetilde{r_{g}}\cdot \widehat{\eta}(x) \; = \; 
   d\widetilde{r}_{g}\widehat{\eta}(\widetilde{r}^{-1}_{g}
    (x)) \; = \; d\widetilde{r}_{g}\widetilde{\xi}(\widetilde{r}^{-1}_{g}(x)) \; + \;   
    {\cal J} d\widetilde{r}_{g} \widetilde{\zeta}(\widetilde{r}^{-1}_{g}(x)) .  
\]  
For $x\in M$ one has  $d\widetilde{r}_{g}\widetilde{\xi}(\widetilde{r}^{-1}_{g}(x)) 
= (r_{g}\widehat{\xi})(x) = \widehat{\xi}(x)$ for all $\widehat{\xi}\in\alpha({\mathfrak g})$.  
Since the real holomorphic vector field $(\widetilde{r}_{g}\cdot\widetilde{\xi})(x)=
  d\widetilde{r}_{g}\widetilde{\xi}(\widetilde{r}^{-1}_{g}(x))$ is the continuation of 
$r_{g}\cdot \widehat{\xi}$ (and this also is true for $\widetilde{r}_{g}\cdot 
\widetilde{\zeta}$), it follows that $\widetilde{r}_{g}\cdot \widehat{\eta}=\widehat{\eta}$ on $W$, 
since $G/H^{\circ}$ is a generic CR-submanifold of $W$.   
In particular, one has 
\[  
   \widehat{\eta}(g\cdot \pi(e)) \; = \; \widehat{\eta}(r_{g}\cdot \pi(e)) \; = \; d\widetilde{r}_{g}  
     \widehat{\eta}(\pi(e)) .
\]  
Hence for $\sigma = \widehat{\eta}_{1} \wedge \ldots \wedge \widehat{\eta}_{n} 
\in V_{\mathfrak g}$ one has 
\[  
     \sigma(g\cdot\pi(e)) \; = \; d\widetilde{r}_{g}\widehat{\eta}_{1}\wedge\ldots\wedge  
           d\widetilde{r}_{g}\widehat{\eta}_{n}(\pi(e))    
\]  
and thus 
\[  
         \widehat {\varphi_{\widehat {\mathfrak g}} }(g\cdot\pi(e))(\sigma) \: = \;  
                \widehat {\varphi_{\widehat {\mathfrak g}} }(\pi(e))\sigma  ,  
\] 
and thus $g\in J$.  
This shows $J = N_{CR}(H)$.  
That the map $G/N_{CR}(H)\to G/J$ is CR follows from the fact that 
$G/H\to G/N_{CR}(H)$ is a CR-submersion and, because the map  
$\widehat {\varphi_{\widehat {\mathfrak g}} }$ is holomorphic, the map 
$\varphi_{\widehat {\mathfrak g}} $ is CR.  
(The CR-structure on $G/J$ contains $d\varphi_{\widehat {\mathfrak g}} (R)$ and is thus 
at most bigger than the structure on $G/N_{CR}(H)$.)   
\end{proof}  


\subsection{Globalization of the Fiber}  

\noindent  
The goal of this section is to prove the globalization result for a 
compact, homogeneous CR-manifold $G/H$ of codimension one and two   
under the assumption that the ${\mathfrak g}$-anticanonical fibration 
is degenerate, i.e., $N_{CR}(H)=G$.  
This implies that the isotropy subgroup $H$ is discrete.  
The proof is a straightforward modification of the proof 
of Satz 1.4.2.1 in \cite{R}.

\begin{theorem}[Satz 1.4.2.1 \cite{R}]  \label{hatsinm}  
Suppose $M=G/H$ is a generic homogeneous CR-manifold of codimension less
than or equal to two.
Let $\mathfrak g$ denote the Lie algebra of $G$.
Assume that $G$ acts effectively on $M$ and the 
${\mathfrak g}$-anticanonical
fibration of $M$ is degenerate, i.e., $N_{CR}(H)=G$.
There exists a simply connected, complex Lie group $\widehat {G}$ and a closed
CR-embedding of the universal covering $\widetilde{G}$ of $G$ into $\widehat {G}$
as well as a discrete subgroup $\widetilde{\Gamma}$ of $\widetilde{G}$ such that
$\widehat {G}/\widetilde{\Gamma}$ is a ${\mathfrak g}$-globalization of
$\tilde{G}/\widetilde{\Gamma}=G/H$.
\end{theorem}

\begin{proof}  
We let  ${\mathfrak m} = {\mathfrak g}\cap {\cal J}{\mathfrak g}$ be the
maximal complex ideal in $\mathfrak g$.
Let $\widehat {M}$ be a $\mathfrak g$-complexification of $M$
and $\widehat {\mathfrak g}$ be the $\widehat {M}$-complexification of $\mathfrak g$.
Further, let $\widehat {G}$ be the 
connected, simply connected, complex Lie group with Lie algebra
$\widehat {\mathfrak g}$ that is uniquely determined up to isomorphism.
We identify $\mathfrak g$ with $\{ \ \widetilde{\xi} \ | \ \xi\in{\mathfrak g} \ 
\}$.
Because the $\mathfrak g$-anticanonical fibration of $G/H$ is degenerate
and because of the assumption on the codimension of the CR-structure
on $M$, it follows that
\[
     \dim_{\mathbb C} \widehat {\mathfrak g} \; \le \; \dim_{\mathbb C} {\mathfrak m} \; + \; 2   .
\]

\vskip 2ex\noindent
Let
\[
     \widehat {\mathfrak g} \; = \;  \widehat {\mathfrak r} \; + \;  \widehat {\mathfrak s}
     \quad\quad \mbox{ and }
     \quad\quad   {\mathfrak g} \; = \;  {\mathfrak r} \; + \;  {\mathfrak 
s}
\]
be a Levi decompositions of $\widehat {\mathfrak g}$ and $\mathfrak g$.
Now $\widehat {\mathfrak m}\cap\widehat {\mathfrak s}$ is an ideal in 
$\widehat {\mathfrak s}$ and so by dimension reasons one has 
\[
       \widehat {\mathfrak m}\cap\widehat {\mathfrak s} \; = \; \widehat {\mathfrak s}
\]
and it follows from this that
\[
       \widehat {\mathfrak s} \; \subset \;  \widehat {\mathfrak m} \; = \;  
       {\mathfrak m} \; \subset \;   {\mathfrak g}  .
\]
In particular, $\widehat {\mathfrak s}={\mathfrak s}$.

\vskip 2ex\noindent
The radical ${\mathfrak r}$ of ${\mathfrak g}$ is obviously given by
${\mathfrak r} = \widehat {\mathfrak r} \cap {\mathfrak g}$.
Let $\widehat {R}$ be the simply connected, complex Lie group with Lie algebra
$\widehat {\mathfrak r}$ and $\widetilde{R}$ be the universal covering of the radical
$R$ of $G$.
Then $\widetilde{G} = \widetilde{R} \cdot S$, where $S$ denotes a (maximal) simply 
connected, semisimple, complex Lie group with Lie algebra 
${\mathfrak s}  = \widehat {\mathfrak s}$.
The group $\widehat {G}$ has Levi-Malcev decomposition 
$\widehat {G}=\widehat {R} \rtimes S$.
Since ${\mathfrak r}$ and $\widehat {\mathfrak r}$ are solvable, the inclusion
${\mathfrak r}\hookrightarrow\widehat {\mathfrak r}$ induces an imbedding
$i:\widetilde{R}\hookrightarrow\widehat {R}$, see Chevalley's result \cite{Chev}.  
The homomorphism $G \to \widehat {G}$ induced by the embedding
${\mathfrak g} \hookrightarrow \widehat {\mathfrak g}$ is given by  
\[  
     (i,{\rm id}_{S}): \widetilde{R} \rtimes S = G \to \widehat {G} = \widehat {R}\rtimes S .    
\]  
Since $S$ contains a maximal, compact subgroup of $\widehat {G}$,
the image of $G$ is closed in $\widehat {G}$, see \cite{Goto}.
Now let $\pi:\widetilde{G} \to G$ be the universal covering and set
$\widetilde{\Gamma} := \pi^{-1}(H)$.
Since $G$ acts effectively on $G/H$ and the ${\mathfrak g}$-anticanonical
fibration is degenerate, $H$ is discrete in $G$ and thus $\widetilde{\Gamma}$ is
discrete in $\widetilde{G}$.
Since $  G \subset\widehat {G}$ is closed, $\widetilde{\Gamma}\subset\widehat {G}$ is 
closed and $\widehat {G}/\widetilde{\Gamma}$ is obviously a ${\mathfrak g}$-globalization of
$ \widetilde{G}/\widetilde{\Gamma}=G/H$.
It is clear that the embedding $\widetilde{G}\hookrightarrow\widehat {G}$ is CR.
\end{proof}  
\section {Globalization}\label{globalization}
\noindent
Our goal here is to derive a condition under which the local
action of $\widehat G$ in a neighborhood of the CR--homogeneous
manifold $M$ can be globalized.  By this we mean that there
is a complex homogeneous manifold $X=\widehat G/\widehat H$ with 
$M$ being the CR--equivariantly equivalent to the $G$--orbit
of its neutral point. Here we do not assume that $M$ is 
compact, but otherwise we operate under the usual assumptions
and notation of this paper.  In particular, $G$ is asssumed
to be embedded in the simply-connected complex Lie group
$\widehat G$.   

\bigskip\noindent
As the following example shows, in order to achieve globalization it might 
be necessary to modify $M$ (see, e.g., \cite{HR} for 
more details and \cite{KZ} for much more general considerations). 
In Section \ref{classification in projective case} we will provide
examples where the necessary modifications are much more serious.
On the other hand, we use there the criteria developed in
the present section to show that ``most'' $M$ can be globalized
without significant modifications. 

\bigskip\noindent
{\bf Example.} Consider the 2--dimensional affine quadric
$X:=\widehat G/\widehat J$, where $\widehat G=\SL_2(\CC)$ and $\widehat J$
is the subgroup of diagonal matrices. Note that $\widehat J$
contains the group consisting of $\pm \Id$ so that
$\widehat G$ acts with this small ineffectivity. Fix $x_0\in X$ as a 
neutral point where $\widehat J=\widehat G_{x_0}$. The unipotent
group $\widehat U\cong \CC$ of upper--triangular matrices
realizes $X$ as the total space of the principal 
$\CC$--bundle 
$\widehat G/\widehat J\to \widehat G/\widehat J\widehat U\cong \PP_1(\CC)$.
Identifying the $\widehat U$--orbit $\widehat U.x_0$ with $\CC$, 
define $\Sigma $ to be the subset of $X$ which corresponds
$\RR^{\ge 0}$. 

\bigskip\noindent
Now let $G=\SU_2$ and observe that every $G$--orbit in $X$
intersects the \emph{slice} $\Sigma $ in exactly one point.
The $G$--orbit $M_{x_0}$ of the neutral point is a copy of the
2--sphere which is embedded as a totally real submanifold.
Otherwise for all $x\in \Sigma \setminus \{x_0\}$ the
orbit $M_x=G.x$ is a hypersurface.  Since $G_x$ is just
the ineffectivity mentioned above, $M_x$ is simply the
group $\PSU_2$ equipped with a left--invariant CR--structure.
  
\bigskip\noindent
Consider the universal cover $\tilde Z$ of the complement $Z$
of $M_{x_0}$ in $X$. Here $G$ acts freely as a group
of holomorphic transformations.  The slice $\Sigma $ lifts
to a slice $\tilde \Sigma $ for the $G$--action.  For 
$\tilde x\in \tilde \Sigma$ the  CR--homogeneous space
$\tilde M_x$ is just the group $G$ equipped with 
a left--invariant CR--structure. It is 
an example of a strongly pseudoconvex hypersurface which
can \emph{not} be filled in to a Stein space. It also can not
be globalized in our sense, because if there would be 
a globalization $\widehat G/\widehat H$, then $\widehat H$ would be a subgroup
of $\widehat J$ and this would force $\widehat H=\widehat J$ so that
$\tilde M_x=M_x$.\qed
\subsection {Homogeneous fibrations}
\noindent
Here we provide a criterion for the existence of a globalization
of the total space of a homogeneous fiber bundle where it
is known that the fiber and base are globalizable. As usual
the connected complex Lie group $\widehat G$ is assumed to be
simply connected and $G$ is a connected, real (not necessarily
totally real) subgroup of $\widehat G$ with 
$\widehat {\lie g}=\lie g+i\lie g$. We assume that the homogeneous
CR--manifold $M=G/H$ is the total space of a $G$--homogeneous
fiber bundle $\pi: G/H\to G/J$ the base of which is the
$G$-orbit of the neutral point in a complex homogeneous
manifold $\widehat G/\widehat J$.  The map $\pi $ is assumed to
be holomorphic and locally $\widehat G$--equivariant in some
neighborhood of $M$.  The connected component 
$J^\circ H/H$ of the fiber $J/H$ is denoted by $F$.

\bigskip\noindent
If the fiber $F$ possesses a $\widehat J$--globalization $\widehat F$, 
then one is naturally led to consider the complex $\widehat G$--manifold
$\widehat M:=\widehat G\times _{\widehat J}\widehat F$.  The $G$--orbit of the
neutral point in $\widehat M$ is indeed the CR--homogeneous space
$M=G/H$.  

\bigskip\noindent
In applications one can at most hope that $\widehat F$ is a complex
manifold equipped with a local holomorphic action of $\widehat J^0$.  
Thus we at first assume that $\widehat J$ is connected, an
assumption that can be realized by going to a $G$-equivariant 
covering space of $M$. In the end this turns out to be 
no assumption at all.  Even though $F$ might not be connected,
since $\widehat J$ is connected, we must assume that $\widehat F$ 
is connected.  So, given $F$ and $\widehat F$, we must replace
them by their connected components.  Finally, we assume that
the holomorphic vector fields coming from the $\widehat J$--action
on $\widehat F$ can be integrated so that the universal cover
$\widehat J_1$ of $\widehat J$ acts holomorphically on $\widehat F$.  If
all of these assumptions are satisfied, we say that \emph{the
fiber and base of the CR--homogeneous bundle $G/H\to G/J$ 
are globalizable}.  The following yields a first criterion
for a $\hat G$--globalization.
\begin {proposition}
If the inclusion $J^\circ \hookrightarrow \widehat J$ induces
a surjective map of fundamental groups, then the
$\widehat J_1$--action on $\widehat F$ descends to a $\widehat J$--action.
\end {proposition}
\begin{proof}
If $J_1^\circ $ is the lift of $J^\circ$ into $\widehat J_1$,
then the condition on surjectivity of the fundamental groups
implies that the kernel $\Lambda $ of the map $J_1^\circ \to J^\circ$
is the same as that of $\widehat J_1\to \widehat J$.  Since $\Lambda $
acts trivially on $\widehat F$, the action of $\widehat J_1$ descends
to that of $\widehat J_1/\Lambda =\widehat J$.
\end {proof}

\bigskip\noindent
In our applications we are only able to answer such homotopy 
questions modulo the ineffectivity of the actions on the
base manifold.  Thus a useful criterion should be given at
that level.  For this we let $\widehat I$ be the ineffectivity
of the $\widehat G$--action on $\widehat G/\widehat J$ and 
$I=\widehat I\cap G$ the $G$--ineffectivity. 
For notational simplicity, let us refer to the following as \\

\noindent
{\bf Condition (C):}
\begin {itemize}
\item The inclusion $J^\circ/(I\cap J^\circ)\hookrightarrow \widehat J/\widehat I$
induces a surjective map of fundamental groups.  
\end{itemize}  
\begin {proposition}\label{technical condition}
If
\begin {equation}\label{intersection condition}
J^\circ \cap \widehat I^\circ =I^\circ \,,
\end {equation}
then condition (C) implies
that the $\widehat J_1$ action on $\widehat F$ descends
to a $\widehat J$--action.
\end {proposition}
\begin {proof}
Condition (\ref{intersection condition})
implies that there is a natural
commutative diagram defined by the homotopy
sequences associated to the fibrations 
$$
\widehat I^\circ\to \widehat J\to \widehat J/\widehat I^\circ
\ \text{and} \
I^\circ\to J^\circ \to J^\circ/I^\circ \,.
$$
Now $\widehat I^\circ$ is
simply connected, because it is a normal, connected subgroup
of a simply connected Lie group. Thus, using the above
mentioned homotopy sequences, in order 
to show that the inclusion $J^\circ \hookrightarrow \widehat J$
induces a surjective map of fundamental groups, we
must only show that
$$
\pi_1(J^\circ/I^\circ)\to \pi_1(\widehat J/\widehat I^\circ)
$$
is surjective. We claim that this follows from condition (C).

\bigskip\noindent
For this observe that if $\gamma $ represents a homotopy 
class in $\widehat J/\widehat I^\circ$ with a base point in $J^\circ/I^\circ$,
then the given surjectivity implies that there is a 
homotopy of its image in $\widehat J/\widehat I$ to a closed curve
in $J^\circ/J^\circ\cap I$.  Since this homotopy can be
lifted to a homotopy of $\gamma $ to a curve in $J^\circ /I^\circ$,
the desired result follows.
\end {proof}

\bigskip\noindent
The following presents a situation where 
Proposition \ref{technical condition} can be
applied. 
\begin {proposition}\label{Abelian action}
If the radical $\widehat R$ of $\widehat G$ is acting on $\widehat G/\widehat J$ as 
an Abelian group, i.e., if $\widehat J$ contains the commutator
group $\widehat R'$, then (\ref{intersection condition}) in 
Proposition \ref{technical condition} is satisfied.
\end {proposition}

\bigskip\noindent
The following Lemma is the general fact behind this
result.
\begin {lemma}
If $\widehat N$ is a connected, complex normal subgroup of $\widehat G$
which contains the commutator subgroup $\widehat R'$, then
$\widehat N\cap G$ is connected.
\end {lemma}
\begin {proof}
At the Lie algebra level we have $\widehat {\lie r}'=\lie r'+i\lie r'$.
Thus $\widehat R'\cap G=R'$. Consequently, $\widehat N\supset R'$ and 
it is sufficient to prove the result in the case where $\widehat R$ 
is Abelian.  Now in general $\widehat N\cap \widehat R$ is just the
radical $\widehat R_{\widehat N}$. In the Abelian case this is 
a vector subspace of $\widehat R$.  Since the same is true of
$R$, it follows that $\widehat N\cap R$ is connected.  Thus it is
only necessary to show that $\widehat N/\widehat R$ has connected
intersection with the image of $G$ in $\widehat G/\widehat R$.

\bigskip\noindent
For this it is convenient to consider a Levi-Malcev
decomposition $G=R\rtimes S$ which lines up with
a Levi-Malcev decomposition $\widehat G=\widehat R\rtimes \widehat S$. By this
we simply mean that $S\subset \widehat S$.  Thus the simple factors
of $S$ are either complex simple factors of $\widehat S$, or real
forms of complex simple factors of $\widehat S$ or antiholomorphic
diagonals of products of two isomorphic simple factors of $\widehat S$.
Since we may identify the image of $G$ with $S$ in the quotient      
$\widehat G/\widehat R=\widehat S$ and $\widehat N/\widehat R$ 
with a product of certain of the simple factors of $\widehat S$, 
the desired result follows.
\end {proof}

\bigskip\noindent
{\it Proof of Proposition \ref{Abelian action}.} The lemma
shows that $\widehat I^0\cap G=I^0$ and, since $I^0\subset J^0$,
condition (\ref{intersection condition}) of 
Proposition \ref{technical condition} is fulfilled. \qed 
\subsection {Remarks on coverings}

\noindent
Here we at first continue under the assumptions of the
previous paragraph, in particular that $\widehat J$ is
connected.  The work there shows the following.
\begin {proposition}
If $\widehat R$ is acting as an Abelian group on $\widehat G/\widehat J$
and condition (C) is fulfilled, then $\widehat J$ acts holomorphically
on $\widehat F$.
\end {proposition}

\bigskip\noindent
In this situation we therefore consider the $G$--orbit 
$G/\tilde H$ of the neutral point in the complex 
$\widehat G$-manifold $\widehat G\times _{\widehat J}\widehat F=\widehat G/\widehat H$. 
Now by construction $H$ fixes this neutral point; so 
$\tilde H\supset H$.  Furthermore, $G/\tilde H$ is
an $F$--bundle over $G/J$ and consequently 
$\dim G/H=\dim G/\tilde H$.  
Therefore $M=G/H\to G/\tilde H=\tilde M$ is a covering.  
Hence, in the setting above
the previous paragraph, if condition (C) is fulfilled
and $\widehat R$ acts as an Abelian group on $\widehat G/\widehat J$,
then after replacing $M$ by a discrete $G$--equivariant
quotient, the local $\widehat G$ action near $M$ can be globalized.

\bigskip\noindent
The example in the previous section shows that a discrete
quotient may in fact be necessary. In that case it is just
the quotient that maps the given nonglobalizable manifold
to the hypersurface orbit in the affine quadric. 

\bigskip\noindent
There is a covering which has been
implicitly used above and which is actually not necessary
for globalization.  This occurs as follows.
Let us not assume that $\widehat J$ is connected.  Then, abusing the
notation which was used above, we let $H_1:=H\cap \widehat J^\circ$
and $J_1:=J\cap \widehat J^\circ$.  Having done this, we may apply the
above results. If globalization conditions are satisfied
and $\widehat G/\widehat H_1$ is the resulting complex $\widehat G$-homogeneous
manifold, then it is possible to return to the original 
situation. For this note that the fiber $F$ of the
resulting map $G\times _{H_1}F\to G/J_1$ is just the
connected component of the fiber $J/H$.  In particular
$H$ fixes the neutral point in $F$ and consequently
$H\subset \widehat H_1$. If we then return to the original
situation by replacing $\widehat J^\circ$ by the original
complex group $\widehat J$ we have the desired globalization,
of course with a possibly disconnected fiber.
Let us now formulate this result for future reference.
For simplicity we bundle together the assumptions of the\\

\noindent
{\bf Standard Situation:}
\begin {itemize}
\item
The complex Lie group $\widehat G$ is 
connected and simply connected.
\item
The real (not necessarily totally real) subgroup $G$ is
connected  with $\lie g+i\lie g=\widehat {\lie g}$.
\item
$M=G/H$ is a homogeneous CR-manifold with
$\widehat G$ acting locally holomorphically on a complex manifold
which is a local neighborhood of $M$.  
\item
$G/H\to G/J$ is a CR-homogeneous fiber bundle.
\item
The base $G/J$ is the $G$--orbit of the neutral point in a 
$\widehat G$--homogeneous space $\widehat G/\widehat J$.  
\item
The connected component $F$ of the fiber $J/H$ possesses
a (connected) globalization $\widehat F$ on which the
universal cover $\widehat J_1$ of the  
connected component $\widehat J^\circ$ holomorphically acts. 
\end {itemize}
\begin {theorem}\label{globalization theorem}
If in the standard situation the radical $\widehat R$ acts
as an Abelian group on $\widehat G/\widehat J$ and condition (C)
is fulfilled, then $\widehat J^\circ$ acts holomorphically
on $\widehat F$ and there exists a globalization
$\widehat G/\widehat H\to \widehat G/\widehat J$ of the CR-bundle
$G/\tilde H\to G/J$, where $M:=G/H\to G/\tilde H=\tilde M$ is a covering.
\end {theorem} 
Note that in the case where $M$ is compact the 
covering $M\to \tilde M$ is at most finite--fibered.
This is a small price to pay for a globalization.

\bigskip\noindent
In situations where we wish to apply Theorem \ref{globalization theorem}
it is quite often only possible to verify a weaker
version of condition (C), namely that
the inclusion $J^0/(I\cap J^0)\hookrightarrow \widehat J/\widehat I$
induces a map of fundamental groups with the property that
the image of $\pi_1(J^0/(I\cap J^0))$ has finite index 
in $\pi_1(\widehat J/\widehat I)$.  We now replace condition (C)
by this weaker version and note the following

\bigskip\noindent
{\bf Zusatz.} \emph{Theorem \ref{globalization theorem} holds under the weakened
version of condition (C)}.
\begin {proof}
If only the weakened version of conditon (C) holds, then the above
shows that nevertheless a finite covering space $\widehat J_1$ of 
$\widehat J$ acts (transitively) on $\widehat F$. Let us write
$\widehat F=\widehat J_1/\widehat I$ and denote by $\Gamma $
the kernel of $\widehat J_1\to \widehat J$.  It is a finite
central subgroup of $\widehat J_1$. If we replace
$\widehat F$ by  
$\widehat F_1:=\widehat J_1/\Gamma \widehat I$, then $\widehat J$ acts on
$\widehat F_1$ and we have the globalization 
$\widehat G\times _{\widehat J}\widehat F_1$.  The $G$-orbit of the
neutral point in this manifold is perhaps a finite quotient
of the original manifold $G/H$, but finite quotients are allowed
in the statement of Theorem \ref{globalization theorem}
\end {proof}

\subsection {The case of the $\lie g$-anticanonical fibration} 

\noindent
Our main application of the globalization criterion is
in the case of the $\lie g$--anticanonical fibration 
$M:=G/H\to G/J\hookrightarrow \widehat G/\widehat J$. Here we
let $\ell $ be the real Lie algebra $\lie j/\lie h$
and $\widehat {\ell}:=\widehat j/\widehat h$. If $L$ and $\widehat L$
are the associated groups, where as in our general setup
$\widehat L$ is taken to be simply connected, then
$F=L/\Gamma $, where $\Gamma $ is discrete.  If
the local $\widehat L$--action can be globalized to
holomorphically act on $\widehat F$ as in the standard 
assumptions, then the globalization criterian can be
applied.

\begin {theorem}
Let $M=G/H\to G/J\hookrightarrow \widehat G/\widehat J$ be
the $\lie g$--anticanonical fibration of the homogeneous
CR-manifold $M$.  Suppose that the connected component
$F=L/\Gamma $ of the fiber possesses an $\widehat L$
globalization $\widehat F$. Then, if $\widehat R$ acts as
an Abelian group on $\widehat G/\widehat J$ and 
condition (C) is fulfilled, $\widehat J^\circ $ acts
holomorphically on $\widehat F$ and there exists a globalization
$\widehat G/\widehat H\to \widehat G/\widehat J$ of the CR-bundle
$G/\tilde H\to G/J$, where $M:=G/H\to G/\tilde H=\tilde M$ is a covering.
\end {theorem}
\section {Structure theorem in the projective case}
In the previous section we gave a criterion for the existence 
of a $\widehat G$--globalization of a CR-homogeneous space $M$.  
To apply this criterion one needs explicit knowledge of certain properties
of the base of the $\lie g$--anticanonical fibration.    
The main purpose of this section is to prove the first of these properties, 
namely, that the radical $\widehat R$ 
of the group $\widehat G$ acts as an Abelian group on the globalization 
$\widehat G/\widehat J$ of the base $G/J$ of the $\lie g$--anticanonical fibration.  
We emphasize that there is no restriction on the codimension 
of the CR-structure on $G/J$ for this to hold.  
The other property one needs in order to apply Theorem 3.2 is condition (C) 
and settings where condition (C) is fulfilled are discussed in the next section.  

\bigskip\noindent  
The base of the $\lie g$--anticanonical fibration 
is itself a CR-homogeneous space which is a $G$--orbit in the
projective space $\PP(V)$ of a $\hat G$--representation space $V$.  
In the notation of the previous section we have the globalization
$$
G/J\hookrightarrow \widehat G/\widehat J \hookrightarrow \PP(V)\,.
$$
We prove here a structure theorem that gives a first description
of this situation.  
It should be underlined that the group $\widehat G$ is only represented on $V$;   
in particular it may be acting with ineffectivity and both $G$ and $\widehat G$ 
are possibly not closed in $\PGL_\CC(V)$.  
Although the globalization criterion
requires information about the $G$-action, we nevertheless replace
$G$ by the closure of its image in $\PGL_\CC(V)$.  
Since we are concerned here with \emph{compact} CR--homogeneous spaces,
this closure stablizes the base of the $\lie g$--anticanonical fibration. 
Using the structure theorem proved here and the 
detailed classification results of the next section, we then 
recapture enough information about the original group in order to 
apply the globalization theory.

\bigskip\noindent  
A complex Lie group is called reductive if it is the complexification of a 
maximal compact subgroup.  
One should recall that reductive Lie groups always carry the structure of  
linear algebraic groups.  

\begin{theorem}
Let $G$ be a connected real
(not necessarily totally real) closed Lie subgroup of $\PGL_\CC(V)$
and $M=G.x_0$ be a compact orbit in $\PP(V)$ with    
$\PP(V)$ assumed to be the projective linear hull of $M$.   
Let $\widehat G$ be the smallest complex Lie group 
containing $G$ in $\PGL_\CC(V)$, i.e., the group corresponding to  
the Lie algebra $\widehat {\lie g}=\lie g +i\lie g$.  
Denote the complex algebraic closure of $G$ by $\bar G$.  
Let $G=R\cdot S$ denote a Levi--Malcev decomposition
of $G$ and $\widehat G=\widehat R\cdot \widehat S$
a Levi--Malcev decomposition of $\widehat G$.    
Then $R$ is central, compact, acts freely on $M$,     
and is totally real with complexification $\widehat R \cong (\CC^*)^k$.  
The group $\widehat G = \widehat S \cdot \widehat R$ is reductive, and hence 
algebraic, and thus $\widehat G = \bar G$.  
Moreover, any maximal compact subgroup of $G$ acts transitively on $M$.  
\end{theorem} 

\noindent  
The result follows from a number of observations that we now give.  

\begin {proposition}
The radical $R$ of $G$ is Abelian.
\end {proposition}
\begin {proof}
The commutator group $R'$ is a real unipotent group and is
in particular a real algebraic subgroup of $\PGL_\CC(V)$.  Thus, 
on the boundary of every $R'$-orbit in $\PP(V)$ there are only $R'$--orbits
of lower--dimension.  Furthermore, every $R'$--orbit is
algebraically diffeomorphic to some $\RR^m$. In particular, the
only possibility for a closed orbit is a fixed point.  Since
$R'$ is a normal subgroup of $G$, it follows that $G$ acts transitively
on the set of $R'$--orbits in $M$.  In particular, the $R'$-orbits
in $M$ are equidimensional
and since $M$ is compact, they must all be closed, i.e. $R'$ fixes
$M$ pointwise. Since $\PP(V)$ is the projective linear hull of
$M$, it follows that $R'=\{\Id\}$.
\end {proof}

\bigskip\noindent
{\bf Remark.} It should be noted that the process of replacing
the original group $G$ by the closure of its representation
on $\PP(V)$ only enlarges the radical.  Thus the above result
guarantees that the radical of the original group is
acting as an Abelian group on $\PP(V)$.\qed

\bigskip\noindent
Since $R$ is Abelian, it has a unique maximal compact subgroup
$T$, i.e., its maximal compact torus. We will now show that
in fact $R=T$. Note that since $T$ is stabilized by conjugation
by $S$ and the group--theoretic automorphism group of $T$ is
discrete, it follows that $T$ is a central subgroup of $G$. 
\begin {proposition} 
The radical $R$ is a central subgroup of $G$.
\end {proposition}
\begin {proof}
The radical $R$ (analytically) decomposes into a $G$--invariant product
$R=T\times V$, where $V$ is the additive group of a vector space.
Since $T$ is central, $R_{G'}\subset V$. By Chevalley's Theorem $G'$ 
is a real algebraic group and therefore so is $R_{G'}$.  Unless it is
trivial, it is noncompact with all orbits noncompact. 
Since $R_{G'}$ is a normal subgroup of $G$, all of its orbits in
$M$ have the same dimension.  Therefore the same argument as
that above which showed that $R'$ is trivial shows here that
$R_{G'}$ is trivial and therefore $R$ is central.
\end {proof} 
\begin {corollary}
The radical of $G$ is a compact torus, i.e., $R=T$.
\end {corollary}
\begin {proof}
Since $R$ is central in $G$, the radical $\bar R$ of the
algebraic closure $\bar G$ is central in $\bar G$.  Thus,
for any two points $x,y\in \bar G.x_0$ the isotropy groups
$\bar R_x$ and $\bar R_y$ agree.  Since the linear hull
of $M$ is the full space $\PP(V)$ and $\bar G.x_0\supset M$,
it follows that $\bar R$ acts freely on $\bar G.x_0$.
Since $\bar R$ is an normal algebraic subgroup of $\bar G$, its orbits
in $\bar G.x_0$ are closed. Now $R$ is a closed
subgroup of $\bar R$.  So all of its orbits are closed
in $\bar G.x_0$ as well; in particular, its orbits in $M$
are closed in $M$.  But $R$ is acting freely on $M$ and 
therefore $R=T$.
\end {proof}

\bigskip\noindent
Since $T$ is a compact torus in a linear group, it is 
totally real and therefore $\widehat T=T^\CC\cong (\CC^*)^m$.
It follows that $\widehat G=\widehat T\widehat S$ is reductive
and therefore $\widehat G=\bar G$.  

\bigskip\noindent
Note that since $G$ is the product $G=T\cdot S$ of a compact torus
and a semisimple group, it is a (real, but not necessarily 
totally real) algebraic subgroup of $\PGL_\CC(V)$.
Thus its isotropy group $G_{x_0}$ has only finitely many
components. 
\begin {corollary} 
Every maximal compact subgroup $K=T\cdot K_S$ of $G$ acts transitively
on $M$.
\end {corollary}
\begin {proof}
Since $M$ is compact and $G_{x_0}$ has only finitely many components,
Montgomery's theorem \cite{Mont} guarantees this.
\end {proof}

\bigskip\noindent
Quotienting out by 
$\widehat T$ leads to the following picture:
\begin {gather}\label{quotient picture}
\begin {matrix}
M & = & G/H & \hookrightarrow & \widehat G/\widehat H &
\hookrightarrow & \PP(V)\\
  &   & \downarrow & &       \downarrow & & & &\\
N & = & G/TH & \hookrightarrow & \widehat G/\widehat T\widehat H &
\hookrightarrow & \PP(W)
\end {matrix}
\end {gather}
Here $G/H\to G/TH$ and 
$\widehat G/\widehat H\to \widehat G/\widehat T\widehat H$ are
principal $T$-- and $\widehat T$--bundles.  Since
$\widehat G/\widehat H\to \widehat G/\widehat T\widehat H$
is a quotient in an algebraic group setting, as the notation
indicates the base is therefore equivariantly embedded as a 
$\hat G$--orbit in a projective space.  

\bigskip\noindent
Now $S$ acts transitively on $N$ and is acting algebraically
on $M$.  Furthermore, it is a normal subgroup of $G$. Hence
by the same argument that has been applied several times above,
the $S$--orbits in $M$ are also compact.  We regard 
$S.x_0$ as a \emph{thick section} for the fibration 
$M=G/H\to G/TH=S/I_S=N$. The connected component of
the fiber of the map $\Sigma :=S.x_0=S/H_S\to S/I_S$ can be regarded
as a subtorus $T_0$ of $T$.  If $T_1$ in
a complementary torus to $T_0$ in $T$, then we have the
following observation.  
\begin {proposition}
The map $T_1 \times \Sigma \to M$ defined by the 
action by $T_1$ and the canonical injection of $\Sigma $
realizes $M$ as a finite quotient $T_1\times _\Gamma \Sigma $.
\end {proposition}

\bigskip\noindent
Thus up to finite quotients the CR--homogeneous space is
the product of the totally real torus $T_1$ and the 
CR-homogeneous space $\Sigma $.  This product stucture
may not be the optimal one from the Cauchy--Riemann viewpoint,
because at the complexified level the bundle 
$\widehat G/\widehat H\to \widehat G/\widehat T_1\widehat H$
might not split accordingly.  For example, $\widehat S$
could act transitively on $\widehat G/\widehat H$! 

\bigskip\noindent
Using the fact that $S$ is a normal subgroup of $G$ we are
able split off the maximal complex subgroup of $G$ so that in
all future considerations of projective CR--homogeneous
spaces we may assume that $G$ is a real form of $\widehat G$.
For this we let $\lie l=\lie g\,\cap \,i\lie g$ be the
ideal which defines this subgroup at the Lie algebra level.
Since $\lie g$ is the Lie algebra direct sum 
$\lie g=\lie t\oplus \lie s$ and $\lie t$ is totally real
in $\widehat {\lie g}$, it follows that $\lie l$ is the sum
of the simple summands of $\lie s$ which are simple summands
of $\widehat {\lie s}$.  Thus the associated
complex group $L$ is a product of the factors of $S$ which are
also factors of $\widehat S$. Since $L=\widehat L$ is an algebraic
normal subgroup of $G$ which is acting algebraically on
the compact manifold $M$, its orbits are all isomorphic
to a fixed compact projective algebraic homogeneous space
$Z$ (an $L$--flag manifold).
\begin {proposition}
The complex algebraic bundle 
$\widehat G/\widehat H\to \widehat G/\widehat H\widehat L$
is $\widehat G$--equivariantly trivial.
\end {proposition}
\begin {proof}
The $\widehat L$--isotropy group $\widehat P$ at $x_0$ has
exactly one fixed point in the flag manifold 
$Z=\widehat L/\widehat P$.  Thus the same is true of every
$\widehat L$--orbit in $\widehat G/\widehat H$. Consequently
the fixed point set $\mathcal F:=\Fix(\widehat P)$ 
of $\widehat P$ in $\widehat G/\widehat H$ is a section of
the fibration, and the natural map 
$\widehat L\times \mathcal F\to \widehat G/\widehat H$ factors
through a $\widehat L$-equivariant isomorphism of the product
of $Z=\widehat L/\widehat P$ and the base 
$\widehat G/\widehat H\widehat L$. Since the product $\widehat G_1$ 
of the remaining factors of $\widehat S$ and the complex torus 
$\widehat T$ centralize $\widehat L$, this product stabilizes
the fixed point set $\mathcal F$ and the manifold 
$\widehat G/\widehat H=Z\times \mathcal F=
\widehat L.x_0\times \widehat G_1.x_0$ splits at the level
of the groups as well. 
\end {proof}

\bigskip\noindent
Now let $G_1$ be the normal totally real subgroup of $G$ which
is defined as the product of $T$ with the factors of $S$ 
which are not in $L$.  
In summary one should note that we have shown a ``Borel-Remmert'' 
type structure result \cite{BR}, as was mentioned in the Introduction.   
\begin {theorem}
The CR-homogeneous space $M$ is the Cauchy-Riemann product
of the compact complex flag manifold $Z=L.x_0$ and the
CR--homogeneous space 
$$
M_1=G_1.x_0=G_1/H_1\hookrightarrow \widehat G_1/\widehat H_1\,.
$$ 
\end {theorem}
 \section {Projective homogeneous spaces of codimension at most two}
\label{classification in projective case}
Here we continue in the setting of the previous section with
$\widehat G$ a connected complex Lie group
acting via a representation on a complex projective space
$\PP(V)$ with a real (connected) subgroup $G$ so that the
orbit $G.x_0=:M$ is the compact CR--homogeneous space of 
interest.  The structure theorems allow us to assume that
$\widehat G=\widehat T\widehat S$ is reductive with radical 
$\widehat T$ and semisimple part $\widehat S$ and that
$G=TS$ is a real form. Since $G$ is a real algebraic group
acting algebraically and $M$ is compact, we know that
every maximal compact subgroup $K$ of $G$ acts transitively
on $M$. Our goal here is to give a detailed description of
this situation under the further assumption that
$M$ is at most 2--codimensional as a CR--manifold, i.e.,
at most 2--codimensional in the complex orbit
$\widehat G.x_0=\widehat G/\widehat H$. 

\subsection {Description of the projective globalization}
The complex homogeneous manifold $X:=\widehat G.x_0=\widehat G/\widehat H$
is the $\widehat G$--globalization in projective space of the 
CR--manifold $M=G.x_0$ which is assumed to be of codimension at most
two. Our goal here is to describe $X$ using
now classical methods from the theory of actions of complex
algebraic groups.   
\subsubsection {The spherical property}
As a first step we show here that $X$ is \emph{spherical}.
This notion, which orginated in classical harmonic analysis,
is naturally translated into the setting of actions of complex 
reductive groups to the condition that a Borel subgroup
$\widehat B$ of $\widehat G$ has an open orbit in $X$, see \cite{VK}. 
By definition a Borel subgroup is a maximal, connected solvable
subgroup of $\widehat G$.  Such are complex algebraic 
subgroups and any two are conjugate (see e.g. \cite{B} for
the basic theory).  Therefore the condition \emph{spherical}
is defined independent of the Borel subgroup in question.
This condition turns out to be quite restrictive and leads
to fine classification results which are of particular use
in our situation.

\bigskip\noindent
There are a number of different ways to verify the spherical
property.  Here we focus on the $K$--action and
use the Hamiltonian viewpoint.  For this
we let $\omega $ be a $K$--invariant K\"ahlerian structure
on $X$.  Since we may assume that $K$ is represented as
a group of unitary transformations, we may take this to
be the restriction of the Fubini-Study form.  In particular
we have the associated $K$--equivariant moment map
$\mu:X\to \lie k^*$. Let $\widehat K$ denote the (reductive)
complexification of $K$ in $\widehat G$.   
\begin {proposition}\label{spherical}
The manifold $X$ is $\widehat K$--spherical.
\end {proposition}
\begin {proof}
It is sufficient to show that the $K$--action on 
$X$ is coisotropic, i.e., that (generically) $\mu $--fibers
are contained in the $K$--orbits (\cite {HW}).   
For this we consider a dimension theoretically generic $K$--orbit
$Y := K.x$, which we know to be of codimension at most two.  
The tangent space of the $\mu $--fiber at $x$ is $T_x Y^{\perp_\omega}$, 
see, e.g., (26.3) in \cite{GS}.  
We must show that this tangent space is contained in $T_x Y$.

\bigskip\noindent 
Now the complex tangent space $T^{CR}_x Y$ to $Y$ in $T_x Y$ is 
a complex subspace of $T_xX$ of 
codimension at most two.   
Since $\omega $ is K\"ahlerian, $T_x Y^{\perp_\omega }$ is contained 
in the orthogonal complement
$P_x$ of $T^{CR}_x Y$ with respect to the induced Hermitian metric.  
Note that $P_x$ has a natural real structure with
$P_x^\RR$ being defined as $P_x\cap T_x Y$.    
Since $\widehat K$ acts transitively on $X$, the orbit $Y$ is not complex.  
So there are two cases to consider.  
First, if $P_x$ is complex 1-dimensional, i.e., $\dim_\RR P_x^\RR=1$, 
then $\codim_{\RR}(Y)=1$ and $Y$ is odd-dimensional.  
But $\mu(Y)$ is a flag manifold and thus is even-dimensional;
so the fiber of $\mu \vert Y$ is positive-dimensional.  
Since the tangent space of the $\mu $-fiber at $x$ is
$T_xY^{\perp_\omega }$, the full $\mu $-fiber is 1-dimensional. Therefore it
must be (locally) contained in $Y$, i.e., $Y$ is coisotropic.

\bigskip\noindent  
In the case where $P_{x}$ is 2-dimensional, it follows that 
$\codim_{\RR}(Y)=2$ and $P^{\RR}_{x}$ is 2-dimensional.  
It would be theoretically possible that
$K_x$ acts with positive dimensional orbits in $P_x$.  
But we have chosen $Y$ to be a generic $K$--orbit; so all nearby
$K$--orbits are also 2--codimensional.  
As a result the connected component $K_x^0$ acts trivially on 
the 2--dimensional complement $P_x^\RR$
to $T^{CR}_x Y$ in $T_xY$.  
Consequently, the orbit $N(K_x^0).x$ of the normalizer of
$K^0_x$ in $K$ is at least 2--dimensional. 
Now since $\mu(Y)$ is even dimensional, we know that either 
$\mu|Y$ has two dimensional fibers, in which case the result follows, 
or $\mu$ maps $Y$ bijectively onto a coadjoint orbit.  
In the latter case, in particular due to the fact that
coadjoint orbits are simply-connected, one would have $K_{\mu(x)} = K^0_{x}$.  
But since coadjoint orbits are flag manifolds, $N(K_{\mu(x)}).x$ is finite 
and this is a contradiction.  
Thus $\mu|Y$ has 2-dimensional fibers and it follows that $Y$ is coisotropic.  
\end {proof}
\subsubsection {Affine-Rational fibrations} 
As a result of Proposition \ref{spherical} we now know that
$X=\widehat G/\widehat H$ is $\widehat G$--spherical.  Since
$K$ acts transitively on $M$, its generic orbits in $X$ are
at most 2--codimensional.  Letting $\widehat G_u$ to be 
a maximal compact subgroup of $\widehat G$ which contains
$K$, the same is true of it.  One therefore says that
$X$ is a $\widehat G$--homogeneous spherical variety of
\emph{rank} at most two. A great deal is known about
spherical varieties (see e.g. \cite{Ak1}, \cite{Ak2}, \cite{BLV}, \cite{Br}, \cite{LV}) 
so that it would be possible to give a detailed list of the manifolds
$X$ which occur in our setting.  Our goal here is to
give sufficient detail so that for any given application
the reader can work out whatever fine point is needed.

\bigskip\noindent
The ``affine--rational'' fibration, which is in a certain 
sense canonical, is the first method which we apply.
As the name indicates, the basic building blocks of this
fibration $\widehat G/\widehat H\to \widehat G/\widehat Q$
are affine and projective rational homogeneous 
varieties, the fiber being affine and the base being rational.
This condition for the base is equivalent to $\widehat Q$
containing a Borel subgroup, i.e., the group $\widehat Q$
is (complex) parabolic.  Using the most elementary aspects
of \emph{root theory}, the parabolic subgroups of $\widehat G$
can be described in complete detail (see, e.g., \cite{B}).    

\bigskip\noindent
Since the possibility of having a nontrivial
finite group $\widehat H/\widehat H^0$ causes only
notational difficulties,  in our discussion here we assume 
that $\widehat H$ is connected. Instead of applying 
the Levi--Malcev decomposition which is applicable in complete
generality, we make use of the Levi--decomposition 
$\widehat H=\widehat H_u\rtimes \widehat H_r$ which is only valid
for algebraic groups. Here $\widehat H_u$ is the unipotent radical
of $\widehat H$ which consists of the unipotent elements of its
radical, and $\widehat H_r$ is a maximal reductive subgroup.  At
this point the ``hat notation'' only indicates that we are dealing
with complex algebraic groups and has nothing particular to do
with the real group $G$.

\bigskip\noindent
If $\widehat H$ is not reductive, then using what we call
the \emph{Weisfeiler method} (see \cite{Hum} $\S 30$) there
is a systematic constructive method for determining a parabolic
group $\widehat Q=\widehat Q_u\rtimes\widehat Q_r$ with 
$\widehat H_u\subset\widehat Q_u$ and 
$\widehat H_r\subset\widehat Q_r$. 
The fibration 
$\widehat G/\widehat H\to \widehat G/\widehat Q$ has a
projective rational manifold as its base. The
fiber $\widehat Q/\widehat H$ is itself a homogeneous
fiber bundle 
$$A:=\widehat Q/\widehat H\to\widehat Q/ \widehat Q_u\widehat H=
\widehat Q_r/\widehat H_r$$  
with fiber the affine homogeneous space 
$\widehat Q_u/\widehat H_u=\CC^n$ and base which is the
quotient of a reductive group by a reductive subgroup.  
The latter is affine algebraic and, since we are dealing with
algebraic homogeneous bundles, it follows that the total
space is also affine.  Thus we refer to 
$\widehat G/\widehat H\to \widehat G/\widehat Q$ as
an \emph{affine--rational} fibration. It is in general
not unique, but we force it a bit in the direction of
unicity by assuming that the fiber is minimal in the
sense that it can not be $\widehat Q$--equivariantly
fibered over a projective rational manifold.

\subsubsection {The fiber of an affine--rational fibration}
In our case of interest where $X$ is at most of rank
two, although the group $\widehat Q$ may not be acting
as a reductive group on the fiber $A$, since the generic
orbits of its maximal compact subgroups, e.g., 
$U:=\widehat G_u\cap \widehat Q$, are at most 2--codimensional,
we view $A$ as a spherical variety of rank at most two in a 
slightly more general sense. We now describe all possible
cases which can occur for $A=\widehat Q/\widehat H$ (see
the table below).  Our detail is sufficient so that by using
elementary root and representation theory precise 
combinatorial descriptions can be determined.

\bigskip\noindent
Recall that the center of $\widehat Q$ is a complex
torus $\widehat T\cong (\CC^*)^r$ and consider the
fibration $\widehat Q/\widehat H\to\widehat Q/\widehat T\widehat H$.  The
fiber is $(\CC^*)^n$.  Since the base is affine and
for each $\CC^*$ the codimension of the generic $U$--orbit
in the base decreases by one, it follows that $0\le n\le 2$.
If $n=2$, it follows that $\widehat H_u=\widehat Q_u$. In that
case the semisimple parts $\widehat H_r^{ss}$ and 
$\widehat Q_r^{ss}$ also agree and in fact $\widehat Q$ is
just acting as $(\CC^*)^2$.  In other words, 
$\widehat G/\widehat H\to \widehat G/\widehat Q$ is 
a $(\CC^*)^2$--principal bundle.

\bigskip\noindent  
If $n=1$, then the generic $U$--orbits in $\widehat Q/\widehat T\widehat H$
are 1--codimensional.  This is a situation which is described
in detail in (\cite {AHR}). There are two cases, i.e., either
$\widehat Q_u$--acts or it doesn't.  If it does, then, again since
the codimension of the $U$ goes down every time one fibers by
a noncompact fiber, it follows that $\widehat Q_u$ acts transitively on
$\widehat Q/\widehat T\widehat H=\CC^m$ and $U$ is acting linearly
there as either the compact symplectic group, $\SU_m$ or $\U_m$.
If $\widehat Q_u$ acts trivially on this manifold, i.e., if 
$\widehat Q_u=\widehat H_u$, then
it is just the quotient $\widehat Q_r^{ss}/\widehat H_r^{ss}$ of the
semisimple parts of the Levi--factors.  Since the compact group
$U$ has 1--codimensional orbits, we know that this manifold is
just a semisimple symmetric space of rank 1, i.e., the tangent bundle
of either the sphere $S^n$ or its 2:1 quotient $\PP_2(\RR)$,
complex projective space $\PP_n(\CC)$, hyperbolic projective
space $\PP_n(\mathbb H)$ or the Cayley plane $\PP_2(\mathbb O)$.
These last results are due to Morimoto and Nagano
(\cite {MN}).

\bigskip\noindent
Up to this point we have only handled the case where the
$\widehat T$--action on $A$ is nontrivial.  If it is trivial,
then we consider the fibration 
$\widehat Q/\widehat H\to \widehat Q/\widehat Q_u\widehat H$
which has $\CC^n$ as its fiber.  As usual there are two cases:
The base is either nontrivial or it isn't!  If it is trivial,
this means that the reductive parts $\widehat Q_r$ and $\widehat H_r$
are the same and that $\widehat H$ is constructed from $\widehat Q$
by removing two root groups.  This is a situation that is easily
classified.  If the base $\widehat Q/\widehat Q_u\widehat H$
is nontrivial, then, since the center of $\widehat Q$ acts trivially
on $A$, this is again the quotient of the semisimple parts
of $\widehat Q_r$ and $\widehat H_r$, i.e., a complex semisimple
symmetric space of rank one as above.  Once one understands
all of the possibilities for the pairs $(\widehat Q,\widehat H)$
with this property, one only needs to sort out those where 
a root group of $\widehat Q_u$ can be removed.  This type of
combinatorial discussion has also been carried out in (\cite {AHR}).

\bigskip\noindent
Finally we come to the case where both the center $\widehat T$
and the unipotent radical $\widehat Q_u$ both act trivially
on the fiber $A$.  
In other words, after moding out ineffectivity $A$ is 
a semisimple affine spherical space of rank two.  
These have been classified by M.Kr\"amer (\cite {Kr})
with very useful additional remarks being given by
D.~Akhiezer (\cite {Ak2}).  
Here is one way of thinking about this classification. 

\bigskip\noindent
With one exception where 
$\widehat Q$ is acting as $\SO_9$ with isotropy $\Spin_7$, all examples occur 
either as affine symmetric spaces of rank two
or are naturally defined bundles over symmetric spaces of rank one
which were described as above. The two infinite series of bundles are 
the following:
\begin {equation} \label{complex}
\SL_{m+1}/SL_m\to \SL_{m+1}/\SL_m\cdot \CC^*
\end {equation}
and
\begin {equation}\label{quaternions}
\Sp_{2n}/(\Sp_{2n-2}\times \CC^*)\to \Sp_{2n}/\Sp_{2n-2}\times \Sp_2\,.
\end {equation}
Let us describe these examples in further detail. In the
case of (\ref{complex}) the base is the tangent bundle
of $\PP_m$ and the fiber is $\CC^*$.  The total space
can be regarded as the tangent bundle of the total space
of the unit circle bundle of the hyperplane bundle over
$\PP_m$.  In the case of (\ref{quaternions}) the base
is the tangent bundle of $\PP_n(\mathbb H)$ and the fiber
is the 2--dimensional affine quadric.  For this latter
point it is important to recall that $\Sp_2=\SL_2$.  Note
that this manifold is the tangent bundle
of a very natural $S^2$--bundle over $\PP_n(\mathbb H)$.

\bigskip\noindent
The above descriptions of the fiber $A$ of the affine-rational
fibration $\widehat G/\widehat H\to \widehat G/\widehat Q$
are summarized in the following table.
Note that the classification of symmetric spaces can be found in \cite{Helg}.   

\begin{table}[tbh]
\renewcommand{\arraystretch}{1.4}  
\begin{tabular}{c c c} \hline 
fiber & base & remarks \\ \hline  \hline
 \multicolumn{3}{ c } { 
$\widehat T$-fibration }  \\ \hline
  $(\CC^{*})^{2}$ & point & $\widehat G/\widehat H 
\stackrel{(\CC^{*})^{2}}{\to} \widehat G/\widehat Q$ is principal \\ 
  $\CC^{*}$ & $\CC^{m}$ & linear $U$-action on 
 $\widehat Q/\widehat T \widehat H$ \\ \hline  \hline 
  \multicolumn{3}{ c } { $\widehat Q_{u}$-fibration } \\   \hline 
  $\CC^{n}$ & point &   $\widehat H = \widehat Q$ - 2 root groups \\  
  $\CC^{n}$ & symmetric space of rank 1 & 
 $\widehat H_{u} = \widehat Q_{u}$ - 1 root group \\ \hline \hline 
  $\CC^{*}$ & $T(\PP_{m})$ & $A=\SL_{m+1}(\CC)/\SL_{m}(\CC)$  \\ 
  $\SL_{2}(\CC)/\CC^{*}$ & $T(\PP_{n}(\mathbb H))$ & 
   $A=\Sp_{2n}(\CC)/(\Sp_{2n-2}(\CC)\times\CC^{*})$  \\ 
      \multicolumn{3}{c }{affine spherical spaces of rank two } \\ 
  \multicolumn{3}{ c } { 
 one exceptional case: $\SO_{9}(\CC)/\Spin_{7}(\CC)$ } \\ \hline  \hline 
\end{tabular}       
\caption{Fiber $A=\widehat Q/\widehat H$ of the affine-rational fibration} 
 \end{table}  

\subsection {The case of a nontrivial radical}
Here we return to the study of the $G$--action on the
$\widehat G$--homogeneous space $\widehat G/\widehat H$.
As always in the section we assume that $M=G.x_0$ is
at most 2--codimensional.  Here we describe how to
handle the situation where the radical $R$ of $G$ is
nontrivial. Of course we do this by a case-by-case
analysis of the sequence of fibrations
$$
\widehat G/\widehat H\to \widehat G/\widehat R\widehat H
\to \widehat G/\widehat Q\,.
$$
The fiber is $(\CC^*)^n$, and, assuming that
$M$ is at most 2--codimensional and that $R$ is nontrivial,
one has $n=1,2$.  
Note that in the case where $n=2$ the group $G$
acts transitively on the base $\widehat G/\widehat Q$.  If
$n=1$, there are two cases: either $\widehat R\widehat H=\widehat Q$,
and $M$ is mapped to a real hypersurface orbit of $G$ in the base,
or $G$ acts transitively on $\widehat G/\widehat Q$ and $M$ is mapped
to a $G$--hypersurface orbit in the noncompact manifold 
$\widehat G/\widehat R\widehat H$.

\bigskip\noindent
In the case where $n=2$ the manifold $M$ is just $(S^1)^2$--principal
bundle over the compact base $\widehat G/\widehat H$ where
$G$ is acting transitively.  If $n=1$ and $G$ is again acting transitively
on the base, then $M$ is just an $S^1$--principal bundle over
the hypersurface orbit of the semisimple part of $G$ in 
$\widehat G/\widehat R\widehat H$.   A precise description
of this hypersurface setting can be found in (\cite {AHR}).
The remaining situation is where $n=1$ and 
$\widehat G/\widehat R\widehat H=\widehat G/\widehat Q$
and $M$ is mapped to a $G$--hypersurface orbit in the
base. Below we give a complete description in the base and
then elementary considerations show which $S^1$--principal
bundles arise over the $G$--hypersurface orbit. 

\bigskip\noindent
Thus, except for this last case which is handled below,
the existence of a nontrivial radical allows us to reduce
the classification to known results.  
\subsection {Transitive action on the base}
The following is a major simplifying step for our description
of the situations where $G$ is noncompact.
\begin {theorem}\label{transitive on base}
If $X=\widehat G/\widehat H$ is noncompact and the
real form $G$ acts transitively on the base $Y=\widehat G/\widehat Q$
of an affine--rational fibration, then $G$ is compact.
\end {theorem} 
The proof will be given in detail for the case where $G$
is semisimple, which we now assume, and at the end
we will note the necessary adjustments to handle
the case where $G$ has a nontrivial radical. 
We  require several intermediate steps where it is assumed
to the contrary that $G$ does act transitively on the base
of the given affine--rational fibration.  Let us begin by
recalling that, although $G$ is a real form of $\widehat G$
in its representation in the projective linear group, it 
is possible that this is not the case for its action
on $Y$.  Fortunately, this can happen
in only one way, namely when one or more simple factors $S$ of $G$
are themselves complex Lie groups which are embedded as a real
form of their complexifications $S\times S$ as \emph{antiholomorphic
diagonals}. This means that there is an antiholomorphic 
automorphism  $\varphi :S\to S$ with $S$ embedded by 
$s\mapsto (s,\varphi (s))$. If either factor of the
complexification $S\times S$ acts trivially on 
$Y$, then the group $S$ which was originally
a real form acts holomorphically as a complex Lie group.

\bigskip\noindent
Now $Y$ splits into a product of homogeneous rational factors
of the simple factors of $\widehat G$ (modulo ineffectivity).
We write it as
$$
Y=\Pi_iY^1_i\times \Pi_jY^{2,1}_j\times \Pi_kY^{2,2}_k\,
$$
where the simple noncomplex factor $S_i$ acts transitively on $Y^1_i$ as
a real form of its simple complexification $\widehat S_i$,
the simple complex factor $S_j$ acts as a real form of
its complexification $S_j\times S_j$ on $Y^{2,1}_j$ where
the latter acts almost effectively, and the simple
complex factor $S_k$ acts as a complex Lie group on $Y^{2,2}_k$
where its complexification also acts as $S_k$. We will immediately
see that most of these possibilities can not occur if $G$ is
to act transitively on $Y$. 

\bigskip\noindent
For example, if we write $Y^{2,1}_j$
as $S/Q_1\times S/Q_2$, then in order for the
diagonally embedded copy of $S$ to act transitively on this
product, we must have a situation where a parabolic subgroup
$Q_2$ of $S$ would necessarily act transitively on the
rational manifold $S/Q_1$. Taking conjugates so that $Q_1$ and
$Q_2$ contain the same Borel subgroup, elementary root considerations
show that this is not possible. Hence, factors of the
type $Y^{2,2}_j$ don't occur.  

\bigskip\noindent
Our strategy for proving Theorem \ref{transitive on base}
is to first handle the case where $G$ and $\widehat G$
are simple (Proposition \ref{simple case}), then the case where $G$ is simple
and complex with its complexification being $\widehat G=G\times G$
(Proposition \ref{complex case}), and finally to put things together in 
Theorem \ref{transitive on base}.
\subsubsection {Transitive action of a simple real form}
\begin {proposition}\label{simple case}
If $\widehat G$ is simple, $\widehat Q$ is 
a proper subgroup of $\widehat G$ and $G$ acts
transitively on $\widehat G/\widehat Q$, then
$G$ is compact.
\end {proposition}
\noindent
The proof requires a bit of preparation.  First, the 
classification work of A.~Onishchik (\cite{O1,O2}) shows that
there are only two series of examples where the situation
in the proposition could occur, i.e., where a simple real form $G$ acts
transitively on a homogeneous rational manifold 
$Y=\widehat G/\widehat Q$ of its simple 
complexification.  They can be described as follows:
\begin {enumerate}
\item
The odd dimensional complex projective space $\PP_{2n-1}$
\item
The space $\mathcal C_n$ of complex structures on $\RR^{2n}$.
\end {enumerate}
These are compact Hermitian symmetric spaces which have respective
isometry groups $\SU_{2n}$ and $\SO_{2n}(\RR)$. As coset spaces
they are then described as $\PP_{2n}=\SU_{2n}/\U(2n-1)$
and $\mathcal C_n=\SO_{2n}(\RR)/\U(n)$.  In the first case
the real form of $\SL_{2n}(\CC)$ which acts transitvely is
the group of quaternionic linear transformations $\SL_{2n}(\mathbb H)$
and of course its maximal compact subgroup $\USp_{2n}$ acts transitively
as well.  In the second case the noncompact real form is $\SO(1,2n-1)$
with maximal compact subgroup $\SO_{\RR}(2n-1)$ also acting transitively
on $\mathcal C_n$.

\bigskip\noindent
To prove Proposition \ref{simple case} amounts to showing that
none of these cases can occur in our situation. One of the
situations which arises (and must be eliminated) is where the
unipotent radical $\widehat Q_u$ acts trivially on the fiber 
$\widehat Q/\widehat H$. Before coming to the proof of
Proposition \ref{simple case}, we note a fact that is
useful in handling that case.
\begin {proposition}\label{transitive action}
Let $\widehat G=K_G^\CC$ be a complex reductive group with reductive subgroups
$\widehat H$ and $\widehat L$ with the property that $\widehat L=K_L^\CC$
has an open orbit $\Omega $ in the affine homogeneous space 
$\widehat G/\widehat H$. 
Assume that the maximal compact group $K_G$ 
contains $K_L$ and that all groups under consideration
are connected. Then the (dimension theoretically) minimal $K_G$--orbits 
in $\widehat G/\widehat H$ are contained in $\Omega $ and are
$K_L$--orbits. In particular, $\Omega $ is also affine.
\end {proposition}
\begin {proof}
The real dimension of a dimension theoretically minimal 
$K_G$--orbit $N$ agrees with the complex dimension of 
$\widehat G/\widehat H$, i.e., that of $\Omega $.  All $K_L$--orbits in 
$\Omega$ have at least this dimension.  Thus, if a $K_G$--orbit
has nonempty intersection with $\Omega $ and is minimial, then
it is contained in $\Omega $.  On the other hand, no minimal
$K_G$--orbit is contained in a proper complex analytic subset
of $\widehat G/\widehat H$, because it totally real with
its real dimension being the same as the complex dimension 
of $\widehat G/\widehat H$. Applying this to the complement
of $\Omega $, the desired result follows.
\end {proof}
\bigskip\noindent
{\bf Example.} It should be pointed out that in the above
setting $\widehat L$ doesn't necessarily 
act transitively on $\widehat G/\widehat H$.
For example, if $\widehat G=\SL_3(\CC)$ and $\widehat H\cong \GL_2(\CC)$
is the stabilizer in $\widehat G$ of the decomposition 
$\CC^3=\Span(e_1)\oplus \Span(e_2,e_3)$, then in fact $\widehat H$ has
an open orbit in $\widehat G/\widehat H$.\qed

\bigskip\noindent
Suppose that the unipotent radical $\widehat Q_u$ acts trivially on the 
fiber $\widehat Q/\widehat H$.
Then, making an intermediate fibration if necessary, we may assume
that this fiber is either $(\CC^*)^n$ or an affine symmetric space of 
the group $\widehat Q_r$.  Since the center of $\widehat Q_r$ is 
1--dimensional in
both cases, if the first case occurs, then $n=1$.  Thus, to 
prove Proposition \ref{simple case} we must eliminate the 
following three cases: 1.) The unipotent radical $\widehat Q_u$ acts 
transitively on $\widehat Q/\widehat H$, 2.) $\widehat Q/\widehat H=\CC^*$, and
3.) The radical of $\widehat Q_r$ acts trivially on $\widehat Q/\widehat H$ and
with $\widehat Q/\widehat H$ being a symmetric space of the
semisimple part $Q_r^{ss}$.

\bigskip\noindent
{\it Proof of Proposition \ref{simple case} in cases 1.) and 2.):}
Note that the representation of $\widehat Q_r$ on $\widehat Q_u$ is 
irreducible and $\widehat Q_u$ can be identified with the 
tangent space of the Hermitian
symmetric space $\widehat G/\widehat Q$.  Thus, if $\widehat Q_u$ acts
transitively on the fiber, then $\widehat G/\widehat H$ is
just the tangent bundle of $\widehat G/\widehat Q$.  Now
in both cases $\widehat Q_r$ is some $\GL_m(\CC)$ whose
center acts as scalar multiplication by $\CC^*$ on the tangent space at
the neutral point.  As usual denote the basic character by
$\Det :\widehat Q_r\to \CC^*$.  Now, using the standard choices
of a Cartan decomposition $\lie g=\lie k\oplus \lie p$ and a
maximal Abelian subspace $\lie a$ of $\lie p$, one checks that the 
isotropy subgroup of the $G$--action on the base contains a 
Cartan subgroup $H=T\times A$ with $\Det (A)=\RR^{>0}$.  Thus
$A$ can not stabilize a compact subset of the tangent space
and therefore this case doesn't occur.  Case 2.) goes in the
same way, because the $\CC^*$--bundle is defined by some
power $\Det ^k$.\qed

\bigskip\noindent
{\it Proof of Proposition \ref{simple case} in case 3.):} Here
we just note that the fiber $\widehat Q/\widehat H$ is an 
affine symmetric space of some $\SL_m(\CC)$.  If it is of rank one,
then it is the tangent bundle of a projective space and the
same argument as above shows that the $G$--orbit can not
be compact. If it is rank two, then at the level of compact
groups it is either a Grassmannian of 2--planes
$\SU(2+q)/S(\U(2)\times \U(q))$ or one of the isolated
cases $\SU(3)/\SO(3)$ or $\SU(6)/\USp(6)$. Here we apply
Proposition \ref{transitive action}: $K.x_0=M$ fibers
over $\widehat G/\widehat Q=K/L$, where $L.x_0$ is
2-codimensional (over $\RR$) in $\widehat Q/\widehat H$.
Since $L.x_0$ can not be complex (it would then be a compact
complex submanifold of an affine manifold), the orbit
$\widehat L.x_0$ of its complexification is open in the
fiber. Thus Proposition \ref{transitive action} implies
that the smaller group $L$ also acts transitively
on the compact symmetric space of rank 2.  It is known
that no smaller compact group acts transtively on any 
Grassmannian except for projective space (\cite{O1,O2}).  

\bigskip\noindent  
A direct check of the groups in our
cases shows that $\SU(3)/\SO(3)$ does not occur.
Finally, $\SU(5)$ does indeed act transitively on
$\SU(6)/\USp(6)=\SU(5)/\USp(4)$.  However, it is
not a symmetric space and does not fiber over
a symmetric space and is therefore not spherical
of rank two!\qed

\subsubsection {Simple real forms which are complex}

\bigskip\noindent
The result in this case can be stated as follows.
\begin {proposition}\label{complex case}
The group $G$ is not a simple complex group which
is a real form of $\widehat G=G\times G$ and which
acts transitively on a positive--dimensional
base $\widehat G/\widehat Q$.
\end {proposition}
\begin {proof}
As in the previous case, the proof requires a bit
of background. First, we consider the decomposition
$G=KAN$.  Since $G$ is complex, $N$ is the unipotent
radical of a complex Borel subgroup and $A\cong (\RR^{>0})^r$
is a noncompact real form of a maximal complex torus. Since
$M$ is compact and $AN$ is acting algebraically, it follows
that it has a fixed point in $M$.  Therefore we may assume
that the isotropy subgroup $\widehat H$ contains its 
complexification $\widehat {AN}=T\ltimes (N\times N)$
in $\widehat G$.  Here $T\cong (\CC^*)^r$ is 
the maximal complex torus mentioned above which is (holomorphically)
diagonally embedded in $\widehat G=G\times G$. 

\bigskip\noindent
To keep things straight we refer to the first factor of the
complexification as $G_1$ and the second as $G_2$ and consider
the sequence of fibrations
$$
G_1\times G_2/(T\ltimes (N\times N))\to
G_1\times G_2/H\to
G_1\times G_2/(H_1\times H_2)\,.
$$ 
Here $H_1$ and $H_2$ are the respective projections of $H$
into the factors $G_1$ and $G_2$.  Since they contain the Borel
group $T\ltimes N$, they are parabolic and these quotients
are compact. By assumption $G$ acts transitively on the
base $\widehat G/\widehat Q$.  Thus $\widehat Q$ contains
one of the factors, say $G_2$, and it follows that
$\widehat Q=H_1\times G_2$.  So the minimality assumption, 
i.e., that 
the fiber can not be nontrivially fibered onto a 
positive--dimensional compact base, implies that $H_2=G_2$.
However, this implies that a maximal semisimple subgroup $H_r^{ss}$
contains $G_2$ as a simple factor and by our effectivity assumption
it therefore follows that $H$ contains a diagonally embedded copy
of $G$.  Since such is a maximal subgroup and $(G\times G)/G$
is affine, this violates our assumption that the base
is positive--dimensional.
\end {proof}  

\bigskip\noindent
As we observed above, if $G$ is simple, complex and embedded
via an antiholomorphic automorphism as a real form of $G\times G$,
then $G$ can not operate transitively on the product 
$G/P_1\times G/P_2$ of $G$--homogeneous rational manifolds.  However,
interesting 2--codimensional orbits can arise in this way.

\bigskip\noindent
{\bf Example.} Let $V$ be a complex vector space and
equip $V\oplus V^*$ with its standard symmetric complex
bilinear form $b$ which is defined by $b(v,f)=f(v)$ and
is invariant by the diagonal action of $G=\SL(V)$.  This 
defines a $G$--invariant, complex hypersurface 
$$
D=\{([v],[f])\in \PP(V)\times \PP(V^*); b(v,f)=0\}\,.
$$
In fact, $G$ acts transitively on both $D$ and its
complement in $\PP(V)\times \PP(V^*)$, the latter
being the rank one affine symmetric space which is
the cotangent bundle of $\PP(V)$.

\bigskip\noindent  
Choosing a basis of $V$ and letting $z$ and $w$ be the
associated coordinates of $V$ and $V^*$, one writes
$b(z,w)=z^tw$.  At the matrix level the $G$--action is
given by $z\mapsto Az$ and $w\mapsto (A^t)^{-1}w$.  One 
can view this as an action on $\PP_n\times \PP_n$
defined by the holomorphic automorphism $A\mapsto (A^t)^{-1}$.
Changing this slightly, we consider the action given
by $(z,w)\mapsto (Az,\varphi (A)w)$, where 
$\varphi(A)=(\bar A^t)^{-1}=A^\dagger$.  In this way
$\PP_n\times \PP_n$ is a two orbit variety with the
lower--dimensional orbit being the real 2--codimension manifold
$M:=\{([z],[w]);z^t\bar w=0\}$.\qed

\bigskip\noindent
The abstract setup of this example is the following. Let
$S$ be a complex Lie group equipped with an antiholomorphic
automorphism $\varphi :S\to S$.  If $I$ is a closed complex
subgroup, then so is $\varphi (I)$. The mapping $\varphi $
induces an antiholomorphic $S$--equivariant diffeomorphism 
$S/I\to S/\varphi (I)$ where the action on the image
space is given by $t(sI):=\varphi(t)sI$.  If $I_1$ and $I_2$ 
are two closed complex subgroups of $S$,
then the diagonal action on $S/I_1\times S/\varphi(I_2)$ of the real 
form $S$ of $S\times S$ defined by $s\mapsto (s,\varphi (s))$ is 
just the transfer of the standard holomorphic diagonal action on 
$S/I_1\times S/I_2$ by the identity on the first factor and this 
antiholomorphic map on the second. If $I_1=I_2=I$ and $Z=S/I$, 
then we write $Z\times \bar Z$ for the product $S/I\times S/\varphi (I)$  
\begin {proposition}\label{antiholomorphic example}
Let $G$ be a simple complex Lie group which is a real form
of $\widehat G=G\times G$ via an antiholomorphic automorphism $\varphi $.
Suppose that $G$ has a closed orbit $M$ of codimension at
most two in a product 
$Y=\widehat G/\widehat Q=G/Q_1\times G/Q_2=Y_1\times Y_2$.
Then $G=\SL(V)$, $Y_1=Z=\PP(V)=G/Q$, $Y_2=\bar Z=G/\varphi (Q)$ and 
$Y=Z\times \bar Z$. 
In coordinates the manifold $M$ is the transfer
of the complex hypersurface $D$ in the above example by the
map $Z\times Z\to Z\times \bar Z$ which is induced by $\varphi $.
\end {proposition}
\begin {proof}
We may write $Q_2$ as $\varphi (Q_2)$ and consider the standard
(holomorphic) diagonal action of $G$ on $G/Q_1\times G/Q_2$. 
The only possibility for a closed orbit of real codimesion
at most two is to have a complex hypersurface orbit $E$ in this
setup.  Since a 
maximal compact subgroup $K$ of $G$ acts transitively on $E$,
one easily checks that the generic orbits of $K$ in the complement
of $E$ are real hypersurfaces, e.g., because the normal bundle
of $E$ can not be topologically trivial.

\bigskip\noindent  
We regard the open $G$--orbit $Y\setminus E$ as a bundle over
$G/Q_1$ with a fiber which is the open orbit of the
parabolic group $Q_1$ in $G/Q_2$ and see that $Q_2$ is
a maximal parabolic, because otherwise $G/Q_2$ would fiber
over some $G/P$ where the parabolic $Q_1$ would necessarily
act transitively which is impossible.  Consequently, 
$Pic(G/Q_2)\cong \ZZ$ and thus
the lower--dimensional $Q_1$--orbit $E\cap Y_2$ is an ample divisor 
and as a result its complement in $G/Q_2$ is affine.

\bigskip\noindent  
From the classification of the homogeneous affine varieties where
the maximal compact subgoup $L$ at hand has real hypersurfaces 
as its generic orbits (see \cite{AHR}) we know that the only
case where there is a larger group than $L^\CC$ acting is 
the case of $\PP(V)$ where the group is the isotropy subgroup
of the $\SL(V)$--action.  Applying this argument to both factors
we have the desired result.
\end {proof}
\subsubsection {The case of several factors}

\bigskip\noindent
We now turn to the final step in the proof of 
Theorem \ref{transitive on base}, namely to
handle the case where the group $G$ is semisimple but
not simple. At the end we remark how to handle
the case of a nontrivial radical.

\bigskip\noindent
{\it Proof of Theorem \ref{transitive on base}:} Let us
review the situation. Here $G$ is a semisimple real form
of $\widehat G$ and 
$M=G.x_0\hookrightarrow \widehat G.x_0=\widehat G/\widehat H$.
We assume that the given affine--rational fibration
$X=\widehat G/\widehat H\to \widehat G/\widehat Q$ has positive
dimensional fiber and that the real form $G$ acts transitively
on the base.  Also recall that we may assume that $G$ is
(up to connected components) the stabilizer of $M$ in $\widehat G$
and that the action of $G$ on  
$X$ is almost effective. 

\bigskip\noindent  
Under these conditions we show that $G$ is compact.  For this we assume to the
contrary that $G_1$ is a simple noncompact real form of $\widehat G_1$
with a complementary semisimple factor $G_2$ so that 
$Y=\widehat G_1/\widehat Q_1\times \widehat G_2/\widehat Q_2$.
We project the base onto the the second factor $\widehat G_2/\widehat Q_2$ 
and consider the fiber 
$\widehat G_1\times \widehat Q_2/\widehat H\widehat Q_2$ of the
associated fibration of $X$.  Here we replace $M$ by its intersection
with this fiber which is a CR--homogeneous manifold with respect to
$G_1\times (G_2\cap \widehat Q_2)$.  Since this group acts transitively
on the base $\widehat G_1/\widehat Q_1$, we are almost in a position
to apply induction on dimension.  However, the way
we have set up our argument, i.e., assuming that $G$ is semisimple,
the induction assumption does not apply.  Thus we must consider
the action of the radical $Z$ of $(G_2\cap \widehat Q_2)$. It acts
as a compact, central subgroup and we consider the intermediate
fibration
$$
\widehat G_1\times \widehat Q_2/\widehat H\to
\widehat G_1\times \widehat Q_2/\widehat H\widehat Z\,.
$$
The fiber is $(\CC^*)^n$ for $n=0,1,2$ with $M$ being fibered
as an $(S^1)^n$--bundle.

\bigskip\noindent
If $n=2$, then the base is already $\widehat G_1/\widehat Q_1$
where the noncompact real form $G_1$ is assumed to act transitively.
If $\widehat G_1$ is simple, then on the base we have one of 
the exceptional examples given to us by Onishchik (\cite{O1,O2}).  
But, as we have already seen, every character on $\widehat Q_1$ restricts
to a charactor on the $G_1$--isotropy which takes its values
in $\RR^{>0}$.  Since such can not stabilize the compact
fiber $(S^1)^2$, it follows that the $G_1$--isotropy acts
trivially on the fiber.  Thus the $G_1$--orbits in $M$ are
sections and at the complex level the bundle is a topologically trivial 
$(\CC^*)^2$--bundle over $\widehat G_1/\widehat Q_1$.  For
complex geometric reasons (the base is rational) or by a simple
group argument, this implies that the bundle is holomorphically
trivial and that the $\widehat G_1$--orbits are sections which
then agree with the $G_1$--orbits.  Thus $\widehat G_1$ stabilizes
the original $M$ (not just the intersection with this fiber), 
contrary to our assumption that $G$ is a real form and agrees
with the stabilizer of $M$ in $\widehat G$.  Note that
if $G_1$ is complex and $\widehat G_1=G_1\times G_1$, then
it would be theoretically possible for $G_1$ to act transitively
on the base. However, this would mean that one of the factors
of $\widehat G_1$ would act trivially and would therefore act
almost effectively on the fiber $(\CC^*)^2$.  This is of course
not possible and therefore the argument for the case $n=2$
is complete.

\bigskip\noindent
If $n=1$, then we consider the base 
$\widehat G_1\times \widehat Q_2/\widehat H\widehat Z$
in more detail.  Now here the unipotent radical of $\widehat Q_2$ 
could theoretically be acting nontrivially.  On the other
hand the relevant group $G_2\cap \widehat Q_2$ for the
CR--fibration is of the form $Z\cdot S$ where $S$ is semisimple.
So in the discussion we may replace $\widehat Q_2$ by
$\widehat Z\widehat S$ so that the group acting on the base
of this intermediate fibration is the semisimple group
$\widehat G_1\times \widehat S$.  By a similar argument to
that above, $G_1\times S$ is acting here as a real form
of $\widehat G_1\times \widehat S$ and the case where one
of the simple real factors is complex can not occur.  Now
it is theoretically possible that some complex simple
factor of $\widehat G\times S$ stablizes the image
of $M$ here and therefore in order to have the induction
assumption, we must factor it out, thereby splitting off
a complex factor of $M$.  If this factor is not $\widehat G_1$,
then we apply induction to the resulting manifold to
obtain the desired contradiction.  If it is $\widehat G_1$,
then we restrict the $\CC^*$--bundle to the corresponding
compact orbit of $\widehat G_1$.  Due to our minimality assumption
on the affine--rational fibration 
$\widehat G/\widehat H\to \widehat G/\widehat Q=Y$, this
orbit is just the same as the $G_1$--orbit in $Y$,
i.e., $\widehat G_1/\widehat Q_1$ where $G_1$ is also
acting transitively.  Thus we are now in a position
to use the same arguments here as we did in the
case above of the $(\CC^*)^2$--bundle over 
$\widehat G_1/\widehat Q_1$ to complete the proof in
this case.

\bigskip\noindent
Finally, if $n=0$, then we may replace $(G_2\cap \widehat Q_2)$
by a semisimple group $S$ so that we are dealing
with a CR--manifold $M.x_0=G_2S.x_0$ of a semisimple
group where the induction assumption may be applied.

\bigskip\noindent
In order to simplify the above discussion, we assumed
that $G$ is semisimple. Inspection of the proof in the
cases where the radical of $\widehat Q_2$ acts 
($n=1,2$ above) shows that exactly the same arguments
handle the case where the radical of $G$ is nontrivial.
This completes the proof of Theorem \ref{transitive on base}.
\qed
\subsection {Actions of real forms on flag manifolds}
We now turn to the case where $G$ is not acting transitively
on the base of an affine--rational fibration 
$\widehat G/\widehat H\to \widehat G/\widehat Q$ and for the
moment only analyze the $G$--action on the homogeneous 
rational manifold $Y=\widehat G/\widehat Q$ where we 
assume that the $G$--orbit $G.x_0=:M$ of interest is
either 1-- or 2--codimensional.  Recall that
$Y$ splits $Y=Y_1\times \cdots \times Y_m$ according to the splitting
of $\widehat G=\widehat G_1 \cdots \widehat G_m$
into its simple factors.  Thus if one
of the simple factors of $G$ happens to be acting as
a complex simple factor of $\widehat G$ we split off
the corresponding factor (say $Y_1$) of $Y$ which will then 
define a product structure $M=Y_1\times M_1$.  Having done
this we see that $M$ is a product of certain factors of
$Y$ together with a $G$--orbit where $G$ is acting as
a real form.  Thus it is sufficient to consider the
case where $G$ is a real form of $\widehat G$.  
  
\bigskip\noindent  
Furthermore, $Y$ splits into a product $Y=Y_1\times \cdots \times Y_m$
according the splitting of $\widehat G$ into factors
$\widehat G_i$ which are the complexifications of the
simple factors $G_i$ of $G$. Since $M$ splits accordingly,
it is enough to study the actions of the $G_i$ on the $Y_i$.
Note that if one of the $\widehat G_i$ is not simple, then
$G_i$ is complex and embedded in $\widehat G_i=G_i\times G_i$
via an antiholomorphic automorphism.  This situation
has been completely described above (see 
Proposition \ref{antiholomorphic example}).  Thus it remains
to consider the case where $G$ is a simple real form
of the simple complex group $\widehat G$.  

\subsubsection {Background}
For the remainder of this paragraph we assume that $G$ is
a simple real form of a simple complex group $\widehat G$
and consider its action on a $\widehat G$--flag manifold
$Y=\widehat G/\widehat Q$.  
Our particular situation is quite special in that $G$ has 
a closed orbit $M=G.x_0$ of codimension at most two.  
Before proceeding with this
case we summarize here the basic general background which
is needed. For details and more information see (\cite {FHW}).

\bigskip\noindent
Basic for our applications is the fact that $G$ has exactly
one closed orbit in $Y$.  In fact $G$ has only finitely many
orbits in $Y$ and therefore it has open orbits.  Recall that
$K$ denotes a maximal compact subgroup of $G$ and that
$\widehat K$ is its complexification in $\widehat G$.  There
is a important duality between the $G$-- and $\widehat K$--orbits
in $Y$ (called Matsuki--duality, see \cite{Mat1,Mat2}) 
which states that to every 
$G$--orbit there is a unique $\widehat K$--orbit which intersects
the $G$--orbit in a $K$--orbit and vice versa.  In the case
of the closed $G$-orbit $M$ this is the simple fact that 
$K$ acts transitively on $M$ and the $\widehat K$--orbit of 
any point in $M$ is open in $Y$. 

\bigskip\noindent  
There is a very natural $K$--invariant
gradient flow on $Y$ which has its critical points at exactly
the $K$--orbits which arise in the duality.  This flow
is tangent to both the $G$-- and the $\widehat K$--orbits and, for example,
retracts each $\widehat K$--orbit onto the $K$--orbit which intersects
the $G$--orbit as in the duality theorem.  We will use this
in the case where the flow retracts the open $\widehat K$--orbit
onto the closed $G$--orbit.

\bigskip\noindent
If $D$ is an open $G$--orbit in $Y$, then the above duality states
that there is a unique compact (complex) $\widehat K$--orbit
$C_0$ in $D$.  It is often of interest to consider it as a 
point in the cycle space of $D$ or $Y$ and, given a choice
of $K$, we therefore refer to it as the \emph{base cycle}.
It can also be characterized as the unique $K$--orbit in $D$
of minimal dimension.

\bigskip\noindent
It would be desirable to have $K$--invariant exhaustions of the
open orbits which reflect the group theoretic situation.  These
are known to be available for \emph{measurable orbits}. If
$G$ is of Hermitian type, i.e., if the center of $K$ is positive
dimensional (and therefore 1--dimensional), every open
$G$-orbit in every $\widehat G$-flag manifold is measurable.  Such exhaustions
$\rho :D\to \RR^{\ge 0}$ are $q$--convex in the sense that
at every point of $D$ the Levi--form $L(\rho )$ has at least $n-q$ 
positive eigenvalues where $n:=\dim D$.  In fact, near the
base cycle it is exactly of signature $(n-q,q)$.

\bigskip\noindent
In the following two paragraphs we give exact descriptions of
the two cases which are relevant for our classification, i.e.,
where $\codim _YM=1,2$.  There are very few possibilities,
a fact that does not at all reflect the general situation for
actions of real forms on flag manifolds.

\subsubsection {Closed orbits of codimension one}
Here $G$ is a simple real form of $\widehat G$ which itself
is simple and we assume that the closed $G$--orbit on
the flag manifold $Y=\widehat G/\widehat Q$ is 
1--codimensional.  Let us begin with an example which is
in fact the only one which arises.

\bigskip\noindent
{\bf Example.}
For nonnegative integers $p$ and $q$ with $n+1=p+q$ we consider the 
action of the real form $G=\SU(p,q)$ of $\widehat G=\SL_{n+1}(\CC)$ on 
projective space $Y:=\PP_{n}(\CC)$. If we regard a point
in $Y$ as a 1--dimensional subspace $L$ in $\CC^{n+1}$, then, restricting
the mixed signature Hermitian form $\langle \ ,\ \rangle_{p,q}$ to
$L$, we speak of $L$ as being positive, negative or isotropic,
depending on this restriction being positive- or negative-definite or
zero.  The group $G$ has three orbits in $Y$: the open sets 
$D_+$ and $D_-$ of positive (resp. negative) lines and the
real hypersurface $M$ of isotropic lines. In coordinates,
the norm of a vector $z=(z_0,\ldots ,z_n)$ is given by
$$
\Vert z\Vert_{p,q}=\sum_0^{p-1}\vert z_i\vert ^2
-\sum_{p}^n\vert z_i\vert ^2\,.
$$ 
This corresponds to the splitting $\CC^{n+1}=V_+\oplus V_-$,
where $V_+:=\{z_p=\ldots z_n=0\}$ and 
$V_-=\{z_0=\ldots =z_{p-1}=0\}$ which is invariant by the
maximal compact subgroup $K=S(U(V_+)\times U(V_-))$ of $G$.
The base cycles in $D_+$ and $D_-$ are $C_+:=\PP(V_+)$ and
$C_-:=\PP(V_-)$, respectively.  Since $G$ has three orbits,
$\widehat K$ does as well, namely the two base cycles and the
complement $Y\setminus (C_+\cup C_-)$ which is the open
$\widehat K$--orbit containing $M$.

\bigskip\noindent
Notice that if $n=2m-1$ is an odd number, then the symplectic
group $\Sp_{2m}(\CC)$ acts transitively on $\PP_n$ and 
the real form $\Sp(2p,2q)$ which is contained in $\SU(2p,2q)$
is at least a candidate for a smaller real form that acts transitively
on $M$.  We will show below that this is in fact the case. 
\qed

\bigskip\noindent  
Now we prove that these examples are the only
ones which arise.  First we note that, since the minimal $G$--orbit is
a real hypersurface, all other $G$--orbits are open. Furthermore,
the $\widehat K$--orbit $\widehat K.x_0=\widehat K/\widehat H$
of a point in $M$ is open.  Let us consider an affine--rational
fibration $\widehat K/\widehat H\to \widehat K/\widehat P$.  As
usual we $\widehat K$--equivariantly compactify this to 
a $\widehat K$--manifold $Y_1$ by taking
the unique $\widehat P$--compactification of the affine fiber and
going to the associated bundle over the base.  In fact we will
see that $\widehat P/\widehat H$ is just $\CC^*$, i.e., that
$\widehat K/\widehat H$ has two ends corresponding to two
base cycles.  If it didn't, then we see that $M$ is the
level surface of a $K$--invariant exhaustion and therefore 
bounds a relatively compact domain in $\widehat K/\widehat H$ and
consequently bounds two domains in $Y$.

\bigskip\noindent  
Therefore $G$ has exactly two open orbits, $D_+$ and $D_-$, with 
their boundaries containing $M$.
Since $M$ is the unique closed $G$--orbit, these are the only
open $G$--orbits in $Y$ and each contains base cycles $C_+$ and
$C_-$.  Consequently the open $\widehat K$--orbit does indeed
have two ends and $\widehat P/\widehat H=\CC^*$.  The compactification
$Y_1$ therefore has two 1--codimensional $\widehat K$-orbits
$E_+$ and $E_-$.  The natural identification of the open
$\widehat K$--orbits in $Y$ and $Y_1$ extends to an equivariant birational
mapping $\pi:Y_1\to Y$ which, since indeterminacies only exist in
codimension two, is in fact regular.  As the notation indicates,
$E_+$ is mapped to $C_+$ and $E_-$ is mapped to $C_-$. Observe
that the principal $\CC^*$--action extends to an action
on $Y$ which fixes the cycles pointwise.  The $S^1$--action
stabilizes every $K$--orbit and centralizes the
$K$--action on $Y$.  By our maximality assumption we may assume
that $\widehat G$ contains this $\CC^*$ as a subgroup, and
one checks that this forces $K$ to contain the $S^1$.  Therefore
$G$ is of Hermitian type!

\bigskip\noindent  
Now let us define $p-1:=\dim C_+$ and $q-1:=\dim C_-$.  Since
$G$ is of Hermitian type, we have the $K$--invariant exhaustions
$\rho _+$ and $\rho _-$ discussed above of the open $G$--orbits 
$D_+$ and $D_-$.  Near $C_+$ the Levi form $L(\rho _+)$ has
signature $(n-(p-1),p-1)$ and near $C_-$ the Levi--form
$L(\rho _-)$ has signature $(n-(q-1),q-1)$.  The essential point
now is that the $\RR^{>0}$ coming from the central
$\CC^*$--action acts transitively on the $K$--hypersurface orbits.
Thus, from the point of view of complex geometry they are all
the same.  Since the positive eigenvalues of the Levi forms
come from the norm in the direction normal to the cycles,
we see that the restriction of the Levi forms to the complex
tangent spaces of these $K$--orbits has signature
$(n-p,p-1)$ and $(n-q,q-1)$ respectively.  But when discussing
the Levi form near $C_+$ we are discussing it from the point
of view of $C_+$ being inside a domain defined by $\rho _+$.    
The same is true of $C_-$. Since the hypersurfaces are 
complex analytically the same and the signature of the
restricted Levi form is a complex analytic invariant, we 
have the following fact.
\begin {proposition}\label{Levi form}
If a $K$--hypersurface orbit in $Y$ is viewed as the boundary
of a domain which contains $C_+$, then the signature of 
its Levi form is $(p-1,q-1)$, where $p-1=\dim C_+$, $q-1=\dim C_-$,
$(p-1)+(q-1)=n-1$ and $n=\dim Y$.
\end {proposition}
\noindent  
We wish to show that $G=\SU(p,q)$ and $Y=\PP_n$ as in the above
example. For this we consider the fibrations 
$E_+=\widehat K/\widehat P\to \widehat K/\widehat P_1=C_+$
and $E_-=\widehat K/\widehat P\to \widehat K/\widehat P_1=C_-$.
For example, following the fibration 
$\widehat K/\widehat H\to \widehat K/\widehat P$ of the open
$\widehat K$--orbit with the first fibration, we have a fibration
of the open orbit which extends to a fibration of the disjoint
union of the open orbit and $C_+$ onto $C_+$.  The fiber
is a smooth blow down of zero--section (to a point)
of the $\CC$--bundle which arises by
adding the one end of the $\CC^*$--bundle over $\widehat P_1/\widehat P$. 
Thus $\widehat P_1/\widehat P$ is a projective space and a simple
dimension count shows that it is $\PP_{n-p}$.  Analogously we
see that $\widehat P_2/\widehat P$ is $\PP_{n-q}$.  This puts us
in a position of being able to prove the main result of this
paragraph.
\begin {theorem}\label{codim one classification}
Let $G$ be a simple real form of a simple complex Lie group
$\widehat G$ whose closed orbit $M$ in the flag manifold 
$Y=\widehat G/\widehat Q$ is 1--codimensional.  Then
$Y=\PP_n$ with $M$ being the manifold
of isotropic lines in $Y$. If $\widehat G=\SL_{n+1}(\CC)$,
then $G=\SU(p,q)$ where $p+q=n+1$.  If $n=2m-1$ is odd,
then the subgroup $\Sp(2p,2q)$, which is a real form
of $\Sp_{2m}(\CC)$ and is a subgroup of $\SU(2p,2q)$, also
acts transitvely on $M$. 
\end {theorem}
\begin {proof}
We will show that the diagonal $\widehat K$--action
on $\widehat K/\widehat P_1\times \widehat K/\widehat P_2$
coming from the two fibrations of $\widehat K/\widehat P$
defines an isomorphism 
$$\label {key isomorphism}
\widehat K/\widehat P\to 
\widehat K/\widehat P_1\times \widehat K/\widehat P_2\cong
\PP_{p-1}\times \PP_{q-1}\,.
$$ 
For this, observe that the fibers $\widehat P_1/\widehat P=\PP_{q-1}$
and $\widehat P_2/\widehat P=\PP_{p_1}$ are transversal to
one another, because they correspond to the negative and positive
eigenspaces of the Levi form.  So, for example, the closed 
$\widehat P_2$--orbit $\widehat P_2/\widehat P$ in 
$\widehat K/\widehat P$ is mapped biholomorphically onto
the base of $\widehat K/\widehat P\to \widehat K/\widehat P_1$.
Thus, as desired, the bases of the the fibrations are
$\PP_{p-1}$ and $\PP_{q-1}$ respectively.

\bigskip\noindent  
Now, modulo the ineffectivity of the $\widehat K$--action on the base 
$\widehat K/\widehat P_1$, a semisimple part of $\widehat P_1$ is
at most $\SL_{p-2}(\CC)$.  Since we may assume without loss of 
generality that $q\ge p$, this is not enough to be the semisimple
part of a group which acts transtively 
on the fiber $\widehat P_1/\widehat P$.  Thus a semisimple part
of the ineffectivity of the $\widehat K$--action on $\widehat K/\widehat P_1$
must bring up the semisimple part of $P_1$ to a group which can
act transitively on $\PP_{q_1}$.  The only possibility for this
is that this ineffectivity itself acts transitively on the fiber
$\PP_{q-1}$.  This proves (\ref{key isomorphism}).

\bigskip\noindent  
To prove the final statement we must only note that the line bundles
on $\PP_{p-1}\times \PP_{q-1}$ are just the tensor products of
pullbacks of powers of hyperplane section bundles from the 
factors.  The only possibility to blow down the 0-- and $\infty$--section
as we do here (obtaining a \emph{smooth} manifold) is the
bundle with Chern number $(1,1)$.  Blowing this down we obtain
the projective space $\PP_n$.

\bigskip\noindent   
Finally, we must show that the only groups which do the job
are $\SU(p,q)$, as a real from of $\SL_{n+1}(\CC)$ or, 
in the case where $n=2m-1$ is odd, $\Sp (2p,2q)$ as a real
form of $\Sp_{2m}(\CC)$.  Since the manifold $M$ is known
to be the bundle space of the $S^1$--principal bundle of
the $(1,1)$--bundle over $\PP_{p-1}\times \PP_{q-1}$, if
there is to be a subgroup of $\SU(p,q)$ which acts transtively
on $M$, its maximal compact subgroup must act transitively
on this product.  Changing the notation, the
maximal compact subgroup $USp_{2p}\times USp_{2q}$ of
the real form $\Sp(2p,2q)$ of $\Sp_{2m}$ does indeed act
transitively on this product.  Further, the $S^1$--bundle
in question is topologically nontrivial; so this maximal
compact subgroup must act transitively there as well.
Thus, in this case the subgroup $\Sp(2p,2q)$ of $\SU(2p,2q)$
also acts transitively on $M$.
\end {proof} 
\subsubsection {Closed orbits of codimension two}
Again we remind the reader of the situation.  Here
$G$ is a semisimple real form of the semisimple
complex group $\widehat G$.  We assume that $G$ has
a closed orbit $M$ of codimension two in the flag manifold
$Y=\widehat G/\widehat Q$ and split the manifold
$Y=Y_1\times \cdots \times Y_m$ according to the
factors $\widehat G_i$ which are complexifications
of the simple factors $G_i$ of $G$.  If we have
a situation where the simple factor $G_i$ is a complex
group and its complexification is $G_i\times G_i$
and $G_i$ does not act transitively on the flag manifold
$Y_i$, then $M$ splits as a product of the other $Y_j$
with $M_i\hookrightarrow Y_i=\PP_n\times \PP_n$ as in
Proposition \ref{antiholomorphic example}. This case being 
handled, we may assume that all of the $\widehat G_i$ are simple.

\bigskip\noindent  
If we split off all factors $Y_i$ where $G_i$ acts transitively,
then we are left with two possibilities: $Y=Y_1\times Y_2$, where
the closed $G_i$--orbit is a real hypersurface as described in
the previous paragraph, and the case which is to be handled
in this paragraph, where both $G$ and $\widehat G$ are simple.
Thus we assume that we have this latter situation. The following
is then the main result of this paragraph.
\begin {theorem}\label{codimension two}
If $G$ and $\widehat G$ are simple and the closed $G$--orbit
$M$ in the flag manifold $\widehat G/\widehat Q$ is 2--codimensional,
then $G=\SL_3(\RR)$, $\widehat G=\SL_3(\CC)$, the
flag manifold $\widehat G/\widehat Q$ is projective space 
$\PP_2(\CC)$ and the closed $G$--orbit is the set of real
points $M=\PP_2(\RR)$.
\end {theorem} 
\noindent  
We will prove this by induction on the dimension $n=\dim Y$.  This 
begins with $n=2$ where the only possibility for $Y$ is $\PP_2(\CC)$.
The only complex group $\widehat G$ which acts transitively
on $\PP_2(\CC)$ is (up to coverings) $\SL_3(\CC)$ which
only has $\SL_3(\RR)$ and $\SU(2,1)$ as its real forms.  The latter
case has been discussed above and the former is exactly the
case which occurs in the theorem. For the remainder of this
paragraph we will operate under the induction assumption.  In
other words, we will be dealing with a flag manifold of 
dimension $n$ where the assumption of the theorem holds and
will assume that the theorem holds for flag manifolds of
every lower dimension.
\begin {proposition}\label {maximal parabolic}
The group $\widehat Q$ is a maximal parabolic subgroup of $\widehat G$.
\end {proposition}
\begin {proof}
If $\widehat Q$ is not a maximal parabolic, then we have
a nontrivial fibration $\widehat G/\widehat Q\to \widehat G/\widehat Q_1$
and we study the various possible situations for the $G$--action
on the base.  If this action is transitive, our situation
reproduces itself in the fiber.  Here, however, we must be
careful with the induction assumption.  If the fiber group
$G_1=G\cap \widehat Q_1$ is acting as in the theorem, then,
modulo ineffectivity this is $\SL_3(\RR)$ which has a maximal
compact subgroup $\SO_3(\RR)$.  In this case we check the
possibilities given to us by Onishchik (\cite{O1,O2}) and see that a
maximal compact subgroup in the $G_1$--isotropy of the base never
has an $\SO_3(\RR)$--quotient.  The other two possibilities, where
we in fact can not apply the induction assumption,
are handled analogously.  These are where the fiber is a product
of hypersurface orbits or where it is an orbit of a complex
simple group which is acting as a real form.  The groups which
are needed for this are simply not to be found as quotients
of the $G_1$--isotropy given to us by Onishchik (\cite{O1,O2}).

\bigskip\noindent  
If $G$ does not act transitively on the base and the
situation of the theorem is reproduced, then $G=\SL_3(\RR)$,
the base is $\PP_2(\CC)$ and the only possibility for
$Y$ is the manifold of full flags where the closed $G$--orbit
is 3-codimensional. Thus this case also doesn't occur.

\bigskip\noindent   
Finally, if the closed $G$--orbit in $\widehat G/\widehat Q_1=\PP_n(\CC)$
is a real hypersurface, then we directly check that neither of
the possibilities for a real form $G$ acts with 2--codimensional
closed orbit on the bigger flag manifold $\widehat G/\widehat Q$.
\end {proof}
\noindent  
We now turn to an analysis of the open $\widehat K$--orbit
$\Omega $ in $Y$.  For this we let $E:=\bd(\Omega )$ be
its boundary.  For example, we wish to show that
the generic $K$--orbit in $\Omega $ is a real hypersurface.
To see this, observe that since the gradient flow which
realizes duality and retracts $\Omega $ to $M$ is 
$K$--invariant, every $K$--orbit in $\Omega $ fibers
over $M$. 
 In particular, if the generic $K$--orbit in
$\Omega $ is not a hypersurface, then \emph{every} 
$K$--orbit in $\Omega $ is 2--codimensional. Let us show
that this is not the case.

\bigskip\noindent 
In the following we say that points in a $\widehat K$--orbit
$\mathcal O_1$ are only accessible from an orbit
$\mathcal O_2$ if the only orbit which has $\mathcal O_1$
in its closure is $\mathcal O_2$.
\begin {lemma}
Let $\mathcal O$ be a $\widehat K$--orbit in $E$ and let
$y_0\in \mathcal O$ only be accessible by
the open orbit $\Omega $. Assume that
all $K$--orbits in $\Omega $ are 2-codimensional over 
$\RR $. Then the codimension of $K.y_0$ in $\mathcal O$
is at most one.
\end {lemma}
\begin {proof}
This follows from the fact that the isotropy group $K_{y_0}$
stabilizes a polydisk transversal to $\mathcal O$ where
its orbits are at most 1-codimensional.
\end {proof}
\bigskip\noindent  
Now if $y_0\in \mathcal O$ is only accessible from the
open orbit, then every point $y\in \mathcal O$ has this
property. Thus we have the following remark.
\begin {corollary}
If $\mathcal O$ contains a point which is only accessible
from $\Omega $, then $\mathcal O$ is not affine.
\end {corollary}
\begin {proof}
If $\mathcal O$ is affine and the semisimple part of 
$\widehat K$ acts nontrivially on $\mathcal O$, then
the minimal $K$--orbit in $\mathcal O$ is not 1--codimensional.
Recall that unless $G$ is of Hermitian type, then $K$
is semisimple and even in the Hermitian case it has only
a 1--dimensional center.  So if the semisimple part 
acted trivially on $\mathcal O$, then $\mathcal O\cong \CC^*$
and $\widehat K$ would have a fixed point on the boundary.
The $G$--orbit of this cycle would be the bounded symmetric
domain in $Y$. But since $Y\not =\PP_1$,  we know that the
minimal $K$--orbit in the bounded symmetric domain is
larger than 1--dimensional.
\end {proof}
\noindent   
Now suppose $E$ contains irreducible components which are 
1--codimensional and let $E_0$ be the union of these components.
Since $\widehat Q$ is maximal, this divisor is ample and the
complement is Stein.  Let $\mathcal O$ be a maximal dimensional
orbit in $E$ in the complement of $E_0$.  It is clearly only accessible
in the above sense and therefore the generic $K$--orbits
in $\mathcal O$ are real hypersurfaces in $\mathcal O$.  Since
the maximal dimensional orbit on the boundary of $\mathcal O$
is only accessible from $\mathcal O$, it follows by the
same argument that the $K$--orbits there are 0--codimensional,
i.e., the $\widehat K$--orbits on the boundary of $\mathcal O$
are even closed in $Y$. Since $E_0$ is an ample divisor,
such orbits must have nonempty intersection with $E_0$ and 
consequently they are contained in $Y$.  Thus we have the
following situation.
\begin {proposition}
Either $E_0=\emptyset$, i.e., the boundary $E$ of the open
$\widehat K$--orbit in $Y$ is everywhere at least 2--codimensional,
or $E_0=E$ is an ample divisor.
\end {proposition}
\begin {proof}
We have just seen that if $E_0$ is nonempty, but $E_0\not =E$,
then every orbit $\mathcal O$ in the complement of $E_0$
in $E$ is closed in the complement of $E_0$ in $Y$.  Since
$E_0$ is ample, this means that every such $\mathcal O$
is affine.  But this was ruled out above (recall we are
still operating under the assumption that every $K$--orbit
in $\Omega $ is 2--codimensional).
\end {proof}
\bigskip\noindent  
Finally, we come to the key point!  If $E_0=\emptyset$,
then the real codimension of $E$ in $Y$ is at least four.
Hence the homology groups $H_2(\Omega )$ and $H_2(Y)$ agree.
Since $\widehat Q$ is maximal, $b_2(Y)=1$.  This information
allows us to prove the following result. 
\begin {theorem}
The boundary $E$ is an ample divisor in $Y$ and 
the generic $K$--orbit in $\Omega $ is a real hypersurface.
\end {theorem}
\begin {proof}
We have seen above that if $E$ is not an ample divisor, then  
$H_2(\Omega )=H_2(Y)=\ZZ$. Now the fiber $\widehat P/\widehat H$
is an affine homogeneous space which is either
a $\CC^*$-bundle over an affine symmetric space or is itself
$(\CC^*)^2$.  In the former cases there are always $K$--orbits
which are of codimension more than two.  So if every 
$K$--orbit is 2--codimensional, then this fiber is $(\CC^*)^2$.
In this case $\widehat P$ is represented on the fiber as
$(\CC^*)^2$ and is therefore not a maximal parabolic.  Thus
$b_2(Y_1)\ge 2$.  

\bigskip\noindent  
Finally, $M$ is a strong deformation retract of $\Omega $
and, even in the Hermitian case where there can be a contribution
from the center of $K$, $\pi_1(M)$ is either $\ZZ$ or finite.
Thus, if all $K$--orbits are 2--codimensional and $E_0\not =E$,
then an application of the homotopy sequence to the fibration
$Y\to Y_1$ shows that $b_2(Y)\ge 2$!  This contradiction means
that either $E$ is an ample divisor or the generic
$K$--orbits in $\Omega $ are hypersurfaces or both.  We claim
both.

\bigskip\noindent  
To see this, observe that if the generic $K$--orbits are hypersurfaces,
then the complement of $\Omega $ consists of the union of the
closed $\widehat K$--orbits, i.e., the cycles. Nagano's theorem
states that in this situation there are at most two $K$--orbits
in $Y$ which are lower--dimensional than hypersurfaces.  Since $M$
is already one of them, it follows that there is exactly 
one cycle in this situation.  Furthermore, $M$ is the minimal
$K$--orbit in $\Omega $.  Since $M$ is 2--codimensional,
the fiber of the affine--rational fibration is therefore also
2--dimensional, i.e., it is either the affine 2--dimensional
quadric or its 2:1--quotient which is $\PP_2$ with a quadric
curve removed.  In any case this fiber is uniquely
equivariantly compactified by adding a copy of $\PP_1$ to   
obtain a two orbit manifold which is mapped regularly
to $Y$ with the lower--dimensional orbit which was constructed
by compactifying the fiber being mapped to $E$.  But this map blows down
the curve that was added to the affine fiber unless the
cycle itself was already a divisor!  Thus, if the generic
$K$--orbit in $\Omega $ is a hypersurface, then $E$ is indeed
an ample divisor.

\bigskip\noindent  
Conversely, if $E$ is an ample divisor, then $\Omega $ is
already affine and consequently $M$ is totally real and therefore
is 2-dimensional.  Since this situation is completely understood,
i.e., we are at the beginning step of the induction,
we actually prove the complete classification result in this
case!
\end {proof}
\noindent  
As the reader has undoubtedly noticed, the proof of 
Theorem \ref{codimension two} is now complete.  Since
$E$ is now known to be an ample divisor, just as in the
last step above, we see that $Y=\PP_2(\CC)$ and $G=\SL_3(\RR)$
is acting on $Y$ as usual.\qed   
\subsection {Classification theorem}
Above we have proved classification results in all situations
which can arise where $M=G/H=G.x_0$ is a compact Cauchy--Riemann 
homogeneous space which is realized as a $G$--orbit in some
projective space.  Here we put this classification together
in a systematic way.

\bigskip\noindent   
In our particular situation we automatically have
the globalization $X=\widehat G/\widehat H=\widehat G.x_0$
on the ambient projective space.
Recall that the radical $T$ of $G$ is a compact torus and
that $G=T\cdot S$ is a product (with possible 
finite intersection) of $T$ with its semisimple part $S$.
The complex group $\widehat G=\widehat T\cdot \widehat S$
splits accordingly and $T$ is a real form of its radical
$\widehat T=(\CC^*)^n$.  However, it is quite possible
that one or more of the simple factors of $S$ is a complex
simple factor of $\widehat S$.  On the other hand, we
have shown that the orbits of such a factor are compact and  
that $M$ is a $G$--equivariant CR--product $M=M'\times Z$,
where $Z=\widehat G_1/\widehat Q_1$ is the compact complex
homogeneous rational manifold which is the orbit of
the product of such factors.  Since we are dealing with
algebraic groups, the situation reproduces
itself in the sense that $M'$ is such an orbit, but with
the advantage that now $G$ is a real form of $\widehat G$.
At various stages of the classification we will see that
there is such a compact complex factor $Z$ and therefore
it is convenient to refer to $M$ as having a certain
form \emph{up to compact complex factors}.  This means
that we are referring to the manifold $M'$.

\bigskip\noindent
Thus, after eliminating such compact factors we begin with
an affine--rational fibration 
$$
X=\widehat G/\widehat H\to \widehat G/\widehat Q=Y
$$ 
which is assumed to be minimal in the sense that there is no
intermediate fibration where the base is compact.  This
restricts to a CR--fibration of $M=G/H\to G/Q=N$.  At the
level of the complex groups the base $\widehat G/\widehat Q$
is a classical object where for example the structure of
$\widehat Q$ is well--understood.  Above we have explained
all possibilities for the fiber $\widehat Q/\widehat H$ and
given $\widehat Q$ in root--theoretic terms it is not difficult
to determine the possibilities for $\widehat H$. Let us
now consider the possibilities for how $M$ is embedded
in $X$.

\bigskip\noindent
{\bf Compact spherical type}

\bigskip\noindent
Above it was shown that if $G$ acts transitively on the
base $\widehat G/\widehat Q$, i.e., 
$$
Y=\widehat G/\widehat Q=G/Q=N\,,
$$  
then $G$ is compact.  We refer to the CR--homogeneous 
manifolds $M=G/H\hookrightarrow \widehat G/\widehat H=X$
which arise in this way as being of \emph{compact spherical type}.
They are constructed as follows.

\bigskip\noindent   
Start with a compact group $G$ with complexification $\widehat G$
and consider a parabolic subgroup $\widehat Q$ with
$Y:=\widehat G/\widehat Q$.  Without
loss of generality it may be assumed that the root description
of $\widehat Q$ is based around a maximal complex torus 
$\widehat T$ in $\widehat Q$ which is the complexification of
a maximal compact torus in $G$.  This root description defines
a maximal reductive subgroup $\widehat L$ of $\widehat Q$
(a Levi--factor $\widehat Q_r$) which is the complexification of a   
compact subgroup $L$ of the compact group $K$.  
Now one checks the list
of possible affine--spherical fibers $\widehat F$ (see the above table)
and determines whether or not there is one where $\widehat L$
can be represented (with possible ineffectivity) as the
reductive group acting on $\widehat F$.  

\bigskip\noindent   
Recalling that $\widehat Q=\widehat L\ltimes \widehat Q_u$, one
then has an action of $\widehat Q$ on $\widehat F$ and 
the $\widehat G$--algebraic homogeneous space 
$$
X=\widehat G/\widehat H=\widehat G\times _{\widehat Q}\widehat F\,.
$$ 
Here we choose $x_0\in X$ to be a point in the standard
embedding of the fiber $\widehat F$ in 
$\widehat G\times_{\widehat Q}\widehat F$ where the orbit
$L.x_0$ is of codimension one or two, depending on the
rank of $\widehat F$.  The CR--homogeneous space of interest
is then $M:=G.x_0$.  It should be emphasized that there is
an open dense set $X_{\text{gen}}$ in 
$X$ where any two $G$--orbits are equivariantly, real--analytically
diffeomorphic.  
We refer to such orbits as being \emph{generic}.  
Each of these will inherit invariant
CR--structure from the ambient manifold.  
One should expect,
however, any two of these orbits are \emph{not} CR--equivalent.

\bigskip\noindent
{\bf Basic models with noncompact symmetry groups}

\bigskip\noindent
Now we turn to the basic building blocks of projective
Cauchy--Riemann manifolds $M:=G/H$ of
codimension one or two which are homogeneous under
a noncompact symmetry group.  These are orbits of
simple real forms $G$ of complex semisimple groups
in $\widehat G$--flag manifolds $X=\widehat G/\widehat H$.
As was shown above, in this case there 
are two infinite series of examples and one
exceptional case.  
The following is a quick summary of these basic models.

\bigskip\noindent
{\bf Mixed signature quadrics}

\bigskip\noindent  
Here $X=Y=\PP_n(\CC)$ and $M=M_{p,q}$ is the manifold of isotropic
lines of the standard Hermitian norm of signature $(p,q)$. We
refer to the real hypersurface $M$ as a \emph{mixed signature
quadric}.  Its stabilizer in the automorphism group of $\PP_n$
is the group $\SU(p,q)$ which is a real form of $\SL_{n+1}(\CC)$.  
In the case where $n+1$ is even where we are
considering the form of signature $(2p,2q)$ the real form
$\Sp(2p,2q)$ of $\Sp_{n+1}(\CC)$ also acts transitively
on $M$.  As is proven above, if $G$ is a noncompact Lie group
acting on $M$ under the standard assumptions of this section
and $M$ is 1--codimensional, then, up to compact factors,
$M$ is a mixed signature quadric $M_{p,q}$ and $G$ one of the
groups just mentioned.

\bigskip\noindent
{\bf Twisted diagonal actions on products of projective space}

\bigskip\noindent
If $G$ itself is a simple complex Lie group and $\varphi :G\to G$ is 
an antiholomorphic automorphism, then we embed it in its
complexification $\widehat G=G\times G=G_1\times G_2$ by the 
twisted diagonal map $(g,\varphi (g))$.   
If $Q_1$ is a parabolic
subgroup of $G_1$, then $Q_2=\varphi (Q_1)$ is parabolic in $G_2$ and
$G$ acts by this embedding on the homogeneous rational manifold
$X=Y=\widehat G/\widehat Q=G_1/Q_1\times G_2/Q_2$.  
This is a context which should be studied further. 
However, in our case
there is only one situation of interest which arises when
$G_1=G_2=\SL_{n+1}(\CC)$ and $Q_1$ is the standard isotropy
group for the $G_1$--action on $\PP_n(\CC)$. In the case
of the antiholomorphic automorphism which is defined by complex
conjugation, i.e., $\varphi (A)=\bar {A}$, $Q_2=Q_1$ and
both factors of $X$ are the same.  The closed $G$--orbit in
question is just the antiholomorphic diagonal
$M=\{([z],[w])\in \PP_n\times \PP_n;z=\bar w\}$.  
Thus in our case where $M$ has
been assumed to be at most 2--codimensional there is only
one examples of interest, i.e., for $n=1$.

\bigskip\noindent
Up to conjugation, the only other antiholomorphic automorphism of 
$\SL_{n+1}(\CC)$ is given by composing the map $\varphi $ above
with the holomorphic outer automorphism $A\to (A^{-1})^t$.  Thus
$X=\PP_n\times \PP_n^*$, where by abuse of notation $\PP_n^*$
denotes $\PP((\CC^n)^*)$. The closed orbit in this case can
be described in standard coordinates as
$$
M=\{([z],[w])\in \PP_n\times \PP_n; z^t\bar w=0\}\,.
$$
It is the image by complex conjugation in the second factor
of the complex hypersurface $\{([z],[w]); z^tw=0\}$ and therefore,
independent of $n$ is 2--codimensional. We refer to 
these CR--homogeneous spaces as those which arise by
\emph{twisted diagonal actions on products of projective spaces}.
In the case $n=1$ the projective space $\PP_1$ is 
holomorphically equivariantly isomorphic to its dual.  Thus
the above case of the antiholomorphic diagonal in 
$\PP_1\times \PP_1$ is the same as the case which was just discussed
in the special situation where $n=1$.  Thus there is only
one series of manifolds $M_n$ which arise by twisted diagonal actions.

\bigskip\noindent
{\bf The real points in $\PP_2(\CC)$}

\bigskip\noindent
If we equip $\CC^n$ with its standard real structure, then
the associated projective space $X=\PP_n(\CC)=\PP(\CC^n)$ is defined
over the reals and its set of real points is just the projectivization
of $\PP(\RR^n)=\PP_n(\RR)$.  The real form $G=\SL_{n+1}(\RR)$
of $\widehat G=\SL_{n+1}(\CC)$ has exactly two orbits on
$\PP_n(\CC)$, namely $\PP_n(\RR)$ and its complement. Since 
the case of $n=1$ is a mixed signature quadric for the group
$\SU(1,1)$, we don't consider it here.  So the only case which
occurs in our context of interest is the closed orbit $\PP_2(\RR)$
of $G=\SL_3(\RR)$ in $\PP_2(\CC)$.

\bigskip\noindent
The following summarizes the classification results on basic
models with noncompact symmetry groups in the standard
situation of this section.
\begin {theorem}\label{noncompact symmetry}
If $G$ is simple and noncompact, then $M$ is one of the 
following:
\begin {enumerate}
\item
A mixed signature quadric $M_{p,q}$ of codimension one.
\item
The 2--codimensional closed orbit $M_n$ in $\PP_n\times \PP_n^*$
of $G=\SL_{n+1}(\CC)$ acting by an the antiholomorphically
twisted diagonal embedding $A\to (A,(\bar A^{-1})^t)$
in $\SL_{n+1}(\CC)\times \SL_{n+1}(\CC)$.
\item
The exceptional example of the real projective plane
$\PP_2(\RR)$ embedded as usual as the closed orbit
of $G=\SL_3(\RR)$ on $\PP_2(\CC)$.
\end {enumerate}
\end {theorem}

\bigskip\noindent
{\bf Splitting Theorem}     

\bigskip\noindent
The following result allows us to put together a clean statement 
of our classification result.  In order to put it in perspective
one should recall that if the real form $G$ acts transitively on be base
$\widehat G/\widehat Q$, then it is compact.
\begin {proposition}\label{splitting theorem}
Let $G_1$ be a noncompact simple factor of $G=G_1\cdot G_2$ 
with complexification $\widehat G_1$. Denote the
corresponding splittings of the complexified group
and base by $\widehat G=\widehat G_1\cdot \widehat G_2$
and 
$X=\widehat G/\widehat Q=
\widehat G_1/\widehat Q_1\times \widehat G_2/\widehat Q_2$.
Assume that $G_1$ does not act transitively on $\widehat G_1/\widehat Q_1$.
Then $X$ splits $\widehat G$--equivariantly as a product
$X=\widehat G_1/\widehat Q_1\times \widehat G_2/\widehat H_2$,
where $\widehat Q_1=\widehat H\cap \widehat G_1$ and 
$\widehat H_2=\widehat H\cap \widehat G_2$.  
The CR--homogeneous manifold  
splits $G$--equivariantly $M=M_1\times M_2$ along the same lines.
\end {proposition}
\noindent  
Before turning to the proof, let us formulate our main classification
theorem which follows directly from the splitting theorem
and the classification results proved above. In it we assume that
we have the standard situation of this section where in particular
$M=G/H$ is at most 2--codimensional in $X=\widehat G/\widehat H$
and the latter has an affine--rational fibration 
$X=\widehat G/\widehat H\to \widehat G/\widehat Q=Y$.  The
induced fibration of $M$ is denoted by $M:=G/H\to G/Q=N$ 
\begin {theorem}\label{main classification}
If $N$ is 2--codimensional in the base, then $M=N$, $X=Y$
and $M$ one of the following:
\begin {enumerate}
\item
A product $M=M_{p,q}\times M_{p',q'}$
of mixed signature quadrics with the groups $G=G_1\times G_2$
and $\widehat G=\widehat G_1\times \widehat G_2$ and complex
model $X=\widehat G_1/\widehat Q_1\times \widehat G_2/\widehat Q_2$
splitting accordingly.
\item
One of the series $M_n$ which is the closed $\SL_{n+1}(\CC)$--orbit
in $\PP_n(\CC)\times \PP_n(\CC)^*$, where $G=\SL_{n+1}(\CC)$ is
acting by the antiholomorphically twisted diagonal representation.
\item
The closed orbit $M=\PP_2(\RR)$ of $G=\SL_3(\RR)$ in $\PP_2(\CC)$.
\end {enumerate}
If $N$ is 1--codimensional in the base and $M$ is 1--codimensional
in $X$, then $X=Y=\PP_n(\CC)$ and $M$ is a mixed signature quadric
with $G$ either $\SU(p,q)$ or $Sp(2p,2q)$.

\medskip\noindent
If $N$ is $1$--codimensional in the base, and $M$ is 2--codimensional
in $X$, then $M=M_1\times M_2$ splits as a product of a mixed 
signature quadric $M_1$ and a manifold $M_2$ of compact spherical
type. This splitting is $G$--equivariant and is defined by
a splitting $X=\widehat G/\widehat H=
\widehat G_1/\widehat Q_1\times \widehat G_2/\widehat H_2=X_1\times X_2$,
where $M_1$ is the mixed signature quadric orbit of the
noncompact real form $G_1$ ($\SU(p,q)$ or $\Sp(2p,2q)$)
in $Y_1=\PP_n$. The group $G_2$ is compact with $M_2$ being
a 1--codimensional $G_2$--orbit in  $X_2=\widehat G_2/\widehat H_2$ 

\medskip\noindent
If $G$ acts transitively on the base, i.e., if $N=Y$, then $M$
is of compact spherical type.
\end {theorem}
\noindent
We now turn to the 

\bigskip\noindent
{\it Proof of Proposition \ref{splitting theorem}.} It is only
necessary to handle the case where the image $N=G/Q$
in the base is the product of the closed orbit $N_1=G_1/Q_1$ of
the simple factor $G_1$ in $Y_1=\widehat G_1/\widehat Q_1$ with
the other factor $Y_2=\widehat G_2/\widehat Q_2=G_2/Q_2=N_2$
of the rational homogeneous base.  For this we consider the
fiber $\widehat F=\widehat G_1\times \widehat Q_2/\widehat H$
of the fibration 
$\widehat G_1\times \widehat G_2/\widehat H\to \widehat G_2/\widehat Q_2$
over this factor. We have the intermediate fibration
$$
\widehat G_1\times \widehat Q_2/\widehat H\to
\widehat G_1\times \widehat Q_2/\widehat H\widehat Q_2=
\widehat G_1/\widehat H_1\to
\widehat G_1/\widehat Q_1\,.
$$
From the minimality assumption on the original affine--rational
fibration it follows that $\widehat H_1=\widehat Q_1$.  

\bigskip\noindent  
The only case of interest is where the fiber of 
the induced fibration of the CR--manifold $F$ is a real hypersurface
in its complexification $\widehat H\widehat Q_2$.  If the real
group acting on this fiber were noncompact, then our
classification would imply that $\widehat H\widehat Q_2/\widehat H$ 
would be compact.  This would mean that $X=Y$ is compact and
$M=M_1\times M_2$ would be a product of mixed signature quadrics.
This is of course possible, but, as we remarked  at the outset of 
this proof, it is only necessary to handle the case where only
one noncompact factor of $G$ has a noncompact orbit in the base
$\widehat G/\widehat Q$.

\bigskip\noindent   
In summary, the above shows that we have a fibration
$$
\widehat G_1\times \widehat Q_2/\widehat H \to \widehat G_1/\widehat Q_1
=\PP_n(\CC)\,,
$$
where the real group acting on the base is either
$\SU(p,q)$ or $\SU(2p,2q)$.  If the mixed signature
quadric in the base is the orbit $N_1=G_1.y_0=G_1/Q_1$, then
we know that $Q_1$ acts as a compact group on the fiber.
Since connected solvable groups of compact groups are Abelian,
this implies that the commutator group of every noncompact
solvable subgroup of $Q_1$ acts trivially on the fiber.  In
particular $\widehat Q_1$ contains a normal subgroup $\widehat I$ 
which contains every such group.  An explicit check of the form of $Q_1$
shows that the only possibility for $\widehat I$ is $\widehat Q_1$ itself.
For example, if one considers the stabilizer in $G_1$ of 
a 2--dimesional complex subspace of $\CC^n$ of signature
$(1,1)$, then the isotropy of this copy of $\SU(1,1)$
at the base point of $N_1$ is a real Borel group whose
commutator is contained in no proper normal subgroup
of $\widehat Q_1$.

\bigskip\noindent  
The proof is now complete, because the above shows that
$\widehat Q_1$ acts trivially on the fiber and therefore
is contained in $\widehat H$.  This implies that
$X=X_1\times X_2=
\widehat G_1/\widehat Q_1\times \widehat G_2/\widehat H_2$
and $M$ splits accordingly.\qed

\subsection {Application to globalization}
Although our goal here is to show that under certain assumptions globalization
is possible, we begin by observing that there can be substantial
problems caused by the presence of certain 3-dimensional
CR-hypersurfaces.

\subsubsection {On the role of the affine quadric}
Here we give several examples to indicate how
the 2-dimensional affine quadric $\SL_2 (\CC)/T$, where $T$
is the subgroup of diagonal matrices, can 
play a role in hindering globalization. This is just the
tangent bundle of $S^2$ where, except for the 0-section,
every $\SO_3(\RR)$ orbit is 1-codimensional and just
a copy of $\SO_3(\RR)$ itself.  Thus the CR-structures
on these orbits are left-invariant hypersurface structures
on the group itself.  The manifolds which cause the difficulty
in globalization are the universal covers of these
manifolds, i.e., left-invariant CR-hypersurface structures
on the group $\SU_2$.  These nonglobalizable structures
can be built into other CR-homogeneous spaces in various
ways in order to construct other nonglobalizatable examples.
Let us discuss several such examples.

\bigskip\noindent
We begin by considering $\SL_2(\CC)$ represented on $\CC^N$
as a direct sum of two irreducible representations, one where
$-\Id$ acts as the identity, i.e., a representation of even 
highest weight, and one where $-\Id$ acts as $-\Id$. Consider
the affine symmetric space $X:=\SL_2(\CC)/D$, where
$D\cong \CC^*$ is the space of diagonal matrices  
in $\SL_2(\CC)$, and equip $\SU_2$ with the left--invariant
CR--structure which comes from one of its real hypersurface
orbits in $X$.  

\bigskip\noindent  
As was noted in the example at the beginning
of section three, with this structure $M=\SU_2$ can not be globalized.
On the other hand, the 2:1 quotient $M/\pm \Id$ is beautifully
globalized in $X$. Here we modify this example slightly to
obtain a 1--codimensional CR--homogeneous space which is much
more difficult to handle.
  
\bigskip\noindent   
For this let $\SL_2(\CC)$ be represented on $\CC^N$
as a direct sum $\rho $ of two irreducible representations, one where
$-\Id$ acts as the identity, i.e., a representation of even 
highest weight, and one where $-\Id$ acts as $-\Id$, and 
let $\widehat G:=\CC^N\rtimes _\rho \SL_2(\CC)$ be the resulting
complex Lie group.  We define $G=\CC^N\rtimes \SU_2$, where
$\SU_2$ is equipped with the CR--structure defined above.
Finally, let $\Gamma $ be a full lattice in $\CC^N$ so that $\CC^N/\Gamma =T$
is a complex torus. Due to the choice of the mixed representation, 
$-\Id\in \SU_2$ does not normalize most $\Gamma $. Hence, we may
choose $\Gamma $ so that the center
$Z=\{\pm \Id\}$ of $\SU_2$ does not act on $M=T\rtimes \SU_2$ 
from the right. In other words, unlike in the example at the 
beginning of this section where it was possible to replace $M$ 
by its quotient with $Z$, there are no simple adjustments
of $M$ which lead to a globalization.

\bigskip\noindent  
One dramatic adjustment is to replace $M$ by its universal cover
$\tilde M=\CC^N\rtimes \SU_2$ and then go down by the $Z$
quotient in $\SU_2$.  This sort of isogeny does allow us to
globalize, but unfortunately the compactness of $M$ is lost.

\bigskip\noindent   
Note that we have the advantage in this construction that
the $\SL_2(\CC)$ which is causing the trouble is a simple
factor of $\widehat G$.  We could have made matters worse
by arranging an affine--rational fibration 
$\widehat G/\widehat J\to \widehat G/\widehat Q$ of 
the base of the anti--canonical bundle where $X$ now
occurs as the fiber $\widehat Q/\widehat J$ and is an orbit of an 
$\SL_2(\CC)$ factor of the Levi--factor
of $\widehat Q$.  
\qed

\bigskip\noindent
The above example shows that whenever the affine quadric $\SL_2(\CC)/D$ 
is involved in the fiber of the affine--rational fibration
we must study the concrete case at hand in order to determine if
$M$ is globalizable.  If $N(D)$ is the normalizer of $D$, then
its 2:1 quotient $\SL_2(\CC)/N(D)$, which is the complement
of a smooth conic in $\PP_2(\CC)$, is equally dangerous.  If
$\widehat G$ is represented as a semisimple group on projective
space and the affine--rational fibration has fiber 
$\widehat Q/\widehat H\cong \CC$, then we can fiber further to 
$\widehat G/\widehat H\to \widehat H/\widehat Q_1$, where
$\widehat Q_1/\widehat Q\cong \PP_1(\CC)$ and $\widehat Q_1/\widehat H$
is again an affine quadric (\cite{AHR}).  So this case must also
be handled on an individual basis. For simplicity of notation
whenever any of these cases occur we say that \emph{the 2--dimensional
affine quadric is involved}.   

\bigskip\noindent
Now the 2-dimensional affine quadric is the affine symmetric space
which is the tangent bundle of the Hermitian symmetric space
$S^2$.  
The phenomenon of nonglobalizability can be traced
to our homotopy condition (see Proposition 3.1) not being satisfied.
Recall that this states that if $G/J$ is a CR-manifold
realized as an orbit in projective space with globalization
$\widehat G/\widehat J$ there, then the inclusion $J\hookrightarrow
\widehat J$ should induce a map $\pi _1(J)\to \pi _1(\widehat J)$ with
at most finite cokernel. The affine quadric is the only
symmetric space of rank one where this does not occur: $J$ is trivial
and $\widehat J\cong \mathbb C^*$.  In the rank two case we find
other examples where the homotopy condition is not satisfied.
Although the affine quadric is a bit hidden in these examples,
it is nevertheless involved. We now indicate how this
happens.

\bigskip\noindent
Let $G$ be semisimple and compact, and let 
$M=G/J$ be an orbit in a projective
space with globalization $\widehat M=\widehat G/\widehat J$ there.
We suppose that $M$ is 2-codimensional in $\widehat M$. As we will
see below, in order to understand the topological obstructions
it is enough to handle the case where $\widehat G/\widehat J$ is 
affine and in that case it is enough to understand the situation
where $\widehat G/\widehat J$ is the tangent bundle of a compact
Hermitian symmetric $G/L$ space of rank two.  There
$G/L$ is a strong deformation retract of $\widehat G/\widehat J$.
Thus our homotopy condition is reduced to showing that
the natural map $\pi_1(J)\to \pi_1(L)$ has finite cokernel.
To show that it is fulfilled it is therefore enough to
show that the semisimple part $L_{ss}$ acts transitively
on $L/J$. 

\bigskip\noindent
One obvious case where $L_{ss}$ does not necessarily act transitively
on $L/J$ is where $G/L$ is a reducible Hermitian symmetric
space, i.e., a product of two Hermitian symmetric spaces
of rank one.  But $L_{ss}$ not acting transitively is equivalent
to one of the factors being $S^2$; in particular the 2-dimensional
affine quadric is involved.

\bigskip\noindent
We need a fact that can, for example, be read from Table E of (\cite {GWZ}) 
and, except for the symmetric space of $E_6$, can be directly verified
by using the standard matrix models of the groups that occur.  
(For the list of Hermitian symmetric spaces see \cite{Helg}.)  
If $G/L$ is an irreducible Hermitian symmetric space of rank two
and $L_{ss}$ does not act transitively on the 2-codimensional
generic orbit $L/J$ in the tangent space of the neutral point,
then $G/L$ is one of the series of symmetric spaces 
$\SO_{p+2}(\mathbb R)/\SO_{p}(\mathbb R)\cdot\SO_{2}(\mathbb R)$.  
The affine quadric is in fact involved in every such space from this 
series and by going to 2:1 coverings one constructs nonglobalizable
examples.  To show this we begin with the simplest example.

\bigskip\noindent 
Let $G = \SO_{4}(\RR)= (\SU_{2}\times \SU_{2})/\Gamma$, 
where $\Gamma$ is the diagonally embedded central subgroup 
of order two.  Note that the diagonal subgroup $\SU_{2}/\Gamma$ can be 
identified with $\SO_{3}(\RR)$.  
In this case the Hermitian symmetric space is 
\[  
        G/L \; = \; \SO_{4}(\RR)/\SO_{2}(\RR)\times \SO_{2}(\RR) \; , 
\]  
where the isotropy group is defined by the standard embedding of 
$\SO_{2}(\RR)\times \SO_{2}(\RR)$ in $\SU_{2}\times \SU_{2}$.  
The generic $G$-orbit $M:=G/J$ in $\widehat G/\widehat J$ 
is just $G$ itself, i.e., $J = \{ e \}$.  
Now let $\SO_{3}(\RR)= \SU_{2}/\Gamma$ be embedded as above, define
$x_1$ to be the base point in $G/L$ and $x_0$ to be the
base point in $\widehat G/\widehat J$, where $M=G.x_0=G/J=G$ and 
such that $x_{0} = g(x_{1})$ with $g\in \SO_{3}(\CC)$.  
It follows that $M_1:=\SO_{3}(\RR) = \SO_{3}(\RR).x_{0}$ is   
a CR-hypersurface in the affine quadric in 
$\SO_{3}(\CC)/\SO_{2}(\CC) = \SO_{3}(\CC).x_{1}$.  
Finally, let 
$$
\widetilde M:=\SU_{2}\times \SU_{2}= \widetilde{G} \to G=M 
$$
be the universal cover.  
Then $\SU_{2} = \widetilde{M_{1}}$ is the universal 
cover of $M_{1}$ equipped with its non-globalizable structure.  
Consequently, $\widetilde{M}$ is also not globalizable.   
  
\bigskip\noindent 
This low-dimensional example turns out to be a special case of a more 
general picture.  
For this, let us return to the affine-rational fibration 
$\widehat{G}/\widehat J \to \widehat G/\widehat Q$ of the 
globalization of the base of the $\mathfrak g$-anticanonical 
fibration.  As was indicated above, it is enough to handle the case   
where the fiber $\widehat Q/\widehat J$ is an 
irreducible rank two affine symmetric space 
which is the tangent bundle of a Hermitian symmetric space.    
The above example fits as follows in this general context.  
Consider the Hermitian symmetric space $\SO_{p+2}/\SO_{p}\cdot \SO_{2}$. 
Then let $\SO_{p-2}(\RR)\times \SO_{4}(\RR)$ be embedded in 
the usual block form so that 
the $\SO_{2}(\RR)$ of the isotropy group 
$\SO_{p}(\RR)\cdot\SO_{2}(\RR)$ is diagonally 
embedded in $\SO_{4}(\RR)$ as above.  
As before, we let $x_{0} = g(x_{1})$, where $x_{1}$ is a base point 
in the Hermitian symmetric space and $g\in \SO_{3}(\CC)$.    
It follows that $M=G.x_{0} = G/J$, where $J = \SO_{p-2}(\RR)\times\Gamma$ 
with $\Gamma$ as above. 
Then analogous to that above, the universal cover  
$\widetilde M = G/J^{\circ}$ contains $\widetilde{M}_{1}= \SU_{2}$ which is 
equipped with a standard nonglobalizable structure.  
Thus $\widetilde M$ itself is nonglobalizable.  

\subsubsection {On the role of real projective space}
\bigskip\noindent
The case where $G$ is represented on the image of the anticanonical
fibration as $\SL_3(\RR)$ with the image of $M$ being 
$G/J=\PP_2(\RR)$ in its projective globalization $\PP_2(\CC)$ 
is one further case where globalization might not be possible.
Similar to the above case, modulo ineffectivity, 
$\pi_1(\widehat J)=\ZZ$ whereas $\pi_1(J)$ is trivial.  In particular,
our homotopy criterion does not guarantee globalization.  On the
other hand, in this case $\widehat G/\widehat J=\PP_2(\CC)$ and so
there is no mixing with fibrations.  Furthermore, if we replace
$G$ by the preimage of the maximal compact subgroup $\SO_3(\RR)$,
then our homotopy condition does guarantee globalization, i.e.,
for the smaller group.

\subsubsection {Conditions for globalization}
Without going into a case-by-case study of exceptions, the following
is the best globalization statement we presently know.
\begin {theorem}
If the affine quadric is not involved and the base $G/J$ of the
anticanonical fibration is not the real projective plane, then
$M$ can be globalized. If $G$ is represented on the base of
the $\lie g$--anticanonical fibration as a compact group, then the 
only exceptions to the existence of globalization are where 
the affine quadric is involved.
\end {theorem}
\begin {proof}
We will apply the homotopy condition of Section \ref{globalization}. 
By Proposition \ref{splitting theorem} it is enough to handle
the individual factors of the base of the anticanonical fibration.
The exceptional case of the projective plane has been discussed
above.  The other case where a noncompact group is involved
is where the factor of $G$ is represented as $\SU(p,q)$ or
$\Sp(2p,2q)$ and the projective CR--homogeneous space of
relevance is the mixed signature quadric. Here the homotopy
condition of Section \ref{globalization} is satisfied. 

\bigskip\noindent   
For example, to see this in the case of $\Sp(2p,2q)$
we consider the decomposition $\CC^n=V_+\oplus V_-$ which
is stabilized by a maximal compact subgroup 
$K\cong \USp(2p)\times \USp(2q)$.  The $\USp(2p)$--isotropy group
of a nonzero vector $v_+\in V_+$ acts on that vector by
a nontrivial character.  Since the same is true of $\USp(2q)$
on $V_-$ and we can choose $v=v_+ + v_-$ to be an isotropic
vector, then a maximal compact subgroup of the 
$\Sp(2p,2q)$--isotropy group at $v$ is just the $K$--isotropy
$K_v$.  The two characters mentioned above give us
a representation of $K_v$ on $S^1\times S^1$ and the appropriately
chosen diagonal will act by the same character on $v$ as
does the $\SL_n(\CC)$--isotropy of the associated point in 
projective space.  In particular, we see that the fundamental
group of the isotropy group of the real form maps onto the
fundamental group of the isotropy group of the complex form.
This shows that the homotopy condition is satisfied in the
case of $\Sp(2p,2q)$.    
The case of $\SU(p,q)$ goes in a similar fashion.

\bigskip\noindent   
Thus it remains to handle the case where $G$ is represented
as a compact group on projective space and the affine quadric
is not involved. For this we consider an affine--rational fibration
$\widehat G/\widehat J\to \widehat G/\widehat Q$ of the
globalization of the base of the anticanonical fibration.
Since $G$ is compact, it acts 
transitively on the base.  If the center $\widehat Z$ of $\widehat Q$ acts 
nontrivially on the fiber, we have an intermediate fibration
$$
\widehat G/\widehat J\to \widehat G/\widehat Z\widehat J\to
\widehat G/\widehat Q\,.
$$
The fiber $\widehat Z\widehat J/\widehat J$ is isomorphic
to $(\CC^*)^n$, $n=1,2$, and the induced fiber of the $G$--orbit
is $(S^1)^n$.  Thus the question of surjectivity of 
homotopy groups reduces to the same question after the
intermediate fibration has been carried out.  Hence,
we may assume that the center of $\widehat Q$ acts trivially
on the fiber.

\bigskip\noindent   
In this case the fiber $\widehat Q/\widehat J$ is either 
$\CC^n$, $n\ge 2$,
or an affine spherical space under the action of a semisimple group. 
Since the affine quadric is not involved, this fiber is either
an affine symmetric space of rank one or two or one of
the two exceptional cases $\SL_{m+1}/\SL_m$ or $\SO_9/\Spin_7$.
The proof is then completed by an analysis of homotopy groups which
is carried out below.
%

\bigskip\noindent  
We consider the long exact homotopy sequences of the fiber
bundles in the rows of the following diagram.
\begin {gather*}
\begin {matrix}
{\widehat J} & \longrightarrow & {\widehat G} & \longrightarrow & {\widehat G/\widehat J} \\
\uparrow     &      &    \uparrow        &      &    \uparrow \\
J      &   \longrightarrow &   G    &   \longrightarrow &  G/J
\end {matrix}
\end {gather*}
Analysis of these sequences shows that in order to guarantee
the surjectivity of $\pi_1(J)\to \pi_1(\widehat J)$ it is
enough for 
\begin {equation}\label{retract}
\pi_2(G/J)\to \pi_2(\widehat G/\widehat J)
\end {equation}
to be surjective.

\bigskip\noindent  
Now the minimal $G$--orbit $G/L$ in $\widehat G/\widehat J$ is
a strong deformation retract of $\widehat G/\widehat J$.  
Thus we must show that the map $G/J\to G/L$ induces a surjective
map at the level of $\pi_2$.  
Thus we look at the homotopy sequence of 
$$
L/J\to G/J\to G/L\,.
$$ 
Since we are really only interested in surjectivity up to
finite cokernel, we see that it is enough to show that
$\pi_1(L/J)$ is finite.  
If $L$ is semisimple, this is immediate.  
So it remains to consider the case where
$Q/L$ is a Hermitian symmetric space and $Q/J$ is 
a generic orbit (of codimension one or two) in its tangent
bundle, i.e., $L/J$ is a generic orbit in its tangent space.
As we showed above, unless the affine quadric is involved,
$L_{ss}$ acts transitively on $L/J$ and thus $\pi_1(L/J)$
is finite.  This then completes the proof of the Theorem.
\end {proof}
 
\section{Fine Classification}  

\subsection{The Levi-nondegenerate Case}   

In this section suppose $M=G/H$ is a compact homogeneous 
CR-manifold of codimension at most two and assume that 
the Levi form of $M$ is nondegenerate.  
Let $G/H \to G/J$ be the $\mathfrak g$-anticanonical fibration of $M$.   
It is then possible to describe the fiber $J/H$ of the 
$\mathfrak g$-anticanonical fibration    
and to give a fine classification for $M$.      
The results in this setting in codimension one were given 
in \cite{AHR}.   

\bigskip\noindent 
We first note that by Proposition \ref{CRnorm} the group $H^{\circ}$ 
is normal in $J=N_{CR}(H)$, since $N_{CR}(H)\subset N(H^{\circ})$.   
Thus $N := J/H^{\circ}$ is a group and $\Gamma := J/H$ is 
a discrete subgroup of $N$.  
We also recall that there is a complexification   
$\widehat N$ of $N$; see Theorem \ref{hatsinm}.     
Without loss of generality one may assume that $\widehat N$ is simply connected 
and let $\widehat N = \widehat R \rtimes\widehat S_{\widehat N}$ be a Levi 
decomposition, where $\widehat R$ denotes the radical of 
$\widehat N$ and $\widehat S_{\widehat N}$ is a maximal semisimple 
subgroup of $\widehat N$.  
Let $\mathfrak n$ denote the Lie algebra of $N$ and 
${\mathfrak m}_{\mathfrak n}$ be its (maximal) complex ideal.  
Further let $\widehat{\mathfrak n}$ denote the Lie algebra 
of the complexification $\widehat N$ of $N$ and assume 
$\widehat{\mathfrak n} = \widehat{ \mathfrak r}\oplus 
\widehat{\mathfrak s}_{\widehat{\mathfrak n}}$ 
is a Levi decomposition of $\widehat{\mathfrak n}$, 
where $\widehat{\mathfrak r}$ denotes 
the radical of $\widehat{\mathfrak n}$ and 
$\widehat{\mathfrak s}_{\widehat{\mathfrak n}}$ a 
maximal semisimple subalgebra of $\widehat{\mathfrak n}$.  
In the proof of Theorem \ref{hatsinm} we observed that  
 \begin{equation}\label{newsinm} 
      \widehat{\mathfrak s}_{\widehat{\mathfrak n}} \; \subset \; 
       {\mathfrak m}_{\mathfrak n}  \; \subset \; \mathfrak n  
\end{equation}    
and we use this fact in the next result.  


\begin{theorem}
Suppose $M=G/H$ is a 2-codimensional compact homogeneous CR-manifold 
whose Levi form is nondegenerate.  
Let $G/H \to G/J$ be the $\mathfrak g$-anticanonical fibration of $M$ and 
let 
\[   
     X \; = \; \widehat G/\widehat H \; \stackrel{\widehat F}{\longrightarrow} \; \widehat G/\widehat J 
\]  
be the corresponding globalization.  
Then the fiber $\widehat F$ is biholomorphic to $(\CC^*)^{k}$ with 
$k \le \codim_{X} M$ and $\widehat N$ is a complex Abelian 
Lie group of dimension at most two.      
Moreover,   
\begin{itemize}  
\item $\dim_{\CC}\widehat F = 0$: $M=G/H \to G/J$ is a (finite) covering 
and a maximal compact subgroup of $G$ acts transitively on $M$.    
\item  $\dim_{\CC}\widehat F = 1$: $\widehat F = \CC^*$, and  the base 
$\Sigma:=G/J$ of the $\mathfrak g$-anticanonical fibration of $M$ is a compact, homogeneous 
CR-hypersurface that is equivariantly embedded in some projective space. 
In particular, $M$ fibers as an $S^{1}$-principal bundle over $\Sigma$.  
\item  $\dim_{\CC}\widehat F = 2$: $\widehat F$ is, up to coverings,    
biholomorphic to $\CC^*\times\CC^*$, and the base 
$G/J=\widehat G/\widehat J = K/L = \widehat K/\widehat P$ 
of the $\mathfrak g$-anticanonical fibration of $M$ is a homogeneous rational manifold.   
One of the following occurs: 
\subitem (i)  $\widehat P$ acts transitively on $\widehat F$ and the complexification 
$\widehat K$ of $K$ acts transitively on $X$ yielding a $\widehat K$-equivariant 
fibration  
\[  
       X \; = \; \widehat K/\widehat I \; \longrightarrow \; \widehat K/\widehat P 
\]  
that is a $\CC^*\times\CC^*$-principal bundle.  
Moreover, $M$ is an $S^1\times S^1$-principal bundle over the homogeneous 
rational manifold $G/J$.  
\subitem (ii)  $\widehat L$ is acting on $\widehat F$ as $\CC^*$ and:  
\subsubitem (a) $\widehat P$ is acting on $\widehat F$ as $\CC^*$ and, up to finite 
coverings, $M$ admits a CR-splitting as a Levi non-degenerate hypersurface times $S^1$  
\subsubitem (b)  $\widehat P$ acts transitively on $\widehat F$ and there is a 
root $\SL_2$ associated to a particular root group $U$ in $\widehat P$.    
The parabolic group $\widehat P$ is a semidirect product $I_{\widehat P}\rtimes A$, 
where $I_{\widehat P}$ is the ineffectivity of the $\widehat P$-action on 
$\widehat F$ and $A$ is a 2-dimensional solvable group that is the 
semidirect product of a 1-dimensional torus with $U$.    
The isotropy $\widehat I$ of the $\widehat K$-action on $X$ is the subgroup 
$I_{\widehat P}$ extended by $U_{\ZZ}$, a subgroup of $U$ isomorphic to $\ZZ$.    
 \end{itemize}   
\end{theorem}     

\begin{proof} 
Since the Levi form of $M$ is nondegenerate,  
${\mathfrak m}_{\mathfrak n} = (0)$.  
Thus by equation (\ref{newsinm}) one has  
$ \widehat{\mathfrak s}_{\widehat{\mathfrak n}} = (0)$.    
So $N$ and $\widehat N$ are solvable.  

\medskip\noindent    
If $\dim_{\CC} \widehat N =0$, then one has a covering 
of one of the manifolds given by Theorem 5.6.  
Since the action of a maximal compact subgroup is transitive and algebraic 
on this projectively embedded manifold, the covering is finite.  
It is straightforward to determine what these are.  

\medskip\noindent  
And, if $\dim_{\CC} \widehat N =1$, the only possibility is  
$N/\Gamma = S^{1} \hookrightarrow \CC^* = \widehat N/\Gamma$,  
because $\Gamma \subset N \subset \widehat N$.   
Then one has a $\CC^*$-bundle over an $\widehat S$-orbit.  
If $\widehat S$ is not transitive on the total space of the bundle, then 
the $\widehat S$-orbits have codimension one in this total 
space and, up to finite coverings, they are sections.  
This implies that a finite covering is biholomorphic to a product.  
Otherwise, $\widehat S$ acts transitively on the total space of the bundle   
$X = \widehat S/\Gamma \to \widehat S/\widehat H$ and its fiber is 
$\CC^*$.    
Thus $M$ is given by an $S^1$-bundle 
over a compact hypersurface orbit of $S$ in $\widehat S/\widehat H$.    
The latter are described in detail in \cite{AHR} and $\widehat H$ is   
acting on the fiber by a character.  
It is straightforward to work out the possibilities in this setting 
and we leave this to the reader.  

\medskip\noindent  
So suppose   $\dim_{\CC} \widehat N >1$.  
By Proposition 7, pp. 50-51 in \cite{HO1} there 
exists a positive dimensional proper closed complex 
subgroup $\widehat I$ of $\widehat N$ that contains $\Gamma$ and 
without loss of generality one may assume that 
$\dim_{\CC}\widehat N/\widehat I = 1$.  
One also has the corresponding fibration $N/\Gamma \to N/I$, where 
$I := N\cap\widehat I$ and $\widehat N/\Gamma\to\widehat N/\widehat I$.  
Now $\widehat N/\widehat I = \CC, \CC^*$, or a complex torus.  
The first case is not possible, since then by the codimension assumption, 
the fiber $I/\Gamma = \widehat I/\Gamma$ would be a compact complex manifold, 
contradicting the assumption that $M$ has a nondegenerate Levi form. 
Let's see that one cannot have a torus as base.  
Here we distinguish two cases depending on $\codim_{\RR}M$.  
Suppose first that $M$ has a codimension two CR-structure.  
By induction the complex fiber $\widehat N/\widehat I$ is either 
$\CC^*\times\CC^*$ or the complex Klein bottle which 
is a two-to-one quotient of $\CC^*\times\CC^*$.  
Hence $\dim_{\CC}\widehat I =2$ and thus $\dim_{\CC}\widehat N =3$.  
Since $M=N/\Gamma$ fibers with real two dimensional fiber, 
$\dim_{\RR} N=4$.   
Thus $\dim_{\CC}  {\mathfrak m}_{\mathfrak n} =1$, 
contradicting the assumption that the Levi form of $M$ 
is nondegenerate.  
So the base $\widehat N/\widehat I$ is biholomorphic to $\CC^*$ 
and the fiber is $\CC^*$, by induction.    
By the additivity of the codimensions of CR-structures, it follows 
that this is an $S^{1}$-bundle over $S^{1}$.  
In particular, $\dim_{\CC} \widehat N =2$ and 
$\widehat N/\Gamma$ fibers as a $\CC^*$-bundle over $\CC^*$ 
which is either biholomorphic to a direct product or a two-to-one covering 
of it is biholomorphic to a product.  

\medskip\noindent 
If the CR-structure of $N/\Gamma$ has  
codimension one, then $\widehat I/\Gamma$ is $\CC^*$ 
by induction and if the base $\widehat N/\widehat I$ would be 
a one dimensional compact, complex torus, then 
$\dim_{\CC}\widehat N=2$, while $\dim_{\RR} N = 3$.  
This again gives the contradiction that ${\mathfrak m}_{\mathfrak n}\not = (0)$.  
It follows that $\widehat N/\Gamma$ must be $\CC^*$.     

\medskip\noindent  
We will now determine the bundle structure    
in the case of a 2-dimensional fiber.  
Here we have the following setup at the level of globalizations:
$$
X=\widehat G/\widehat H\to \widehat G/\widehat J\,,
$$
with fiber $\widehat F=\CC^*\times \CC^*$.  
The base $Z$ is of the form
$$
\widehat S/\widehat Q=S/I=\widehat K/\widehat P=K/L\,,
$$
where $K$ is a maximal compact subgroup of $S$ which is fixed
for the discussion.

\medskip\noindent
Our remarks revolve around the study of the actions of $\widehat P$
and $L$ on the fiber $\widehat F$.  
Since $\widehat F = \CC^*\times\CC^*$, any semisimple subgroup of $L$ 
acts trivially on $\widehat F$.  
Hence $L$ acts on $\widehat F$ as a solvable group and   
since $L$ is compact, $L$ is acting as an Abelian group.
There are two cases that are relatively easy to handle.   

\medskip\noindent   
First we assume that $L$ acts trivially on $F$.  
Then the ineffectivity of the $\widehat P$-action on $\widehat F$ 
is a normal subgroup of $\widehat P$ that contains a 
maximal torus.  
Since the unipotent radical of the parabolic group $\widehat P$ 
is a product of 1-dimensional subgroups and each of these is 
normalized by this maximal torus, it follows 
that this ineffectivity must be all of $\widehat P$.  
In other words, $\widehat P$ acts trivially on $\widehat F$.  
Thus the $\widehat K$-orbits which are the same as the 
$K$-orbits trivialize the bundle.  
However, this case does not occur because the Levi-form 
would be degenerate and we are assuming that this is not so.  
Next if $L$ acts transitively on $F$, then $\widehat P$ 
acts transitively on $\widehat F$ as an algebraic group.  
But then $\widehat K$ acts transitively on $X$ and $K$ 
acts transitively on $M$.  
This yields the principal bundle   
\[  
X=\widehat K/\widehat I\to \widehat K/\widehat P\  
\]    
and the corresponding fibration of $M$ is an 
$S^{1}\times S^{1}$-principal bundle.  

\medskip\noindent
The remaining case occurs when the generic $L$--orbits
are 1-dimensional.  
Throughout the discussion we will use the
fact that $\hat F$ can be equivariantly fibered.  
There are two ways this can happen, depending on whether or not the 
full fiber group $\widehat J$ is acting as an Abelian group.
In both cases, however, we have the intermediate fibration
\begin {equation}
\label{intermediate fibration}
\widehat G/\widehat H\to \widehat G/\widehat I\to \widehat G/\widehat J\,,
\end {equation}
where both fibers are $\CC^*$.  
This induces an $S^1$--bundle fibration of the torus $F$ over $S^1$.  
If $L$ acts transitively on the base $\widehat J/\widehat I$ 
of this fibration, then it is clear that every $L$--orbit in $F$ 
is a copy of $S^1$ and that $L$ is acting as $S^1$.  
Otherwise, $L$ is only acting in the fibers of this fibration.  
But also in that case, linearizing the
$L$--isotropy shows that $L$ acts transitively on all fibers and is
acting as $S^1$.  
So in general $L$ can be regarded as acting
by a free $S^1$-action and $L^\CC$ by $\CC^*$.  

\medskip\noindent
Again there are two cases:
\begin {enumerate}
\item
The group $\widehat P$ is acting as $\CC^*$.
\item
The group $\widehat P$ acts transitively on $\widehat F$.   
\end {enumerate}
In the first case we claim that this only happens if $M$ is CR-isomorphic 
to the product of a Levi non-degenerate hypersurface with $S^{1}$.  
We look at the $\widehat K$-action on the
base $\widehat G/\widehat I$ of the intermediate fibration
in (\ref{intermediate fibration}).  
If $\widehat K$ acts transitively on $\widehat G/\widehat I$, then (up to finite covers) 
the $\widehat K$--orbits in the total space $\widehat G/\widehat H$
are sections of the bundle $\widehat G/\widehat H\to \widehat G/\widehat I$.  
Let's go to such a covering.  
Then for $p\in M$, we consider $N=K.p\hookrightarrow \widehat K.p=Y$.  
The complex tangent space $T_{p}^{CR} M = T_{p}^{CR} N$
and is contained in the tangent space of the section $Y$.  
Now recall that the Levi--form is computed by bracketing 
$(1,0)$-fields with $(0,1)$-fields which are (complexified) tangent fields to $M$.
In this case such fields will be tangent to $N$.  
Thus the values of the Levi-form lie in the tangent space to $Y$.  
Since we are assuming that the Levi-cone is open,
then this is certainly not the case in the situation at hand, 
unless $M$ is a CR-product of a Levi non-degenerate hypersurface 
with $S^{1}$.   
Next assume that $\widehat K$ does not act transitively on the base.  
Hence the base $\widehat G/\widehat I$ is a product $\CC^*\times Z$ and 
we go to its (canonically defined) holomorphic reduction.   
Since $L^\CC$ is acting as $\CC^*$, it follows that $\widehat K$
acts transitively on the fibers of the holomorphic reduction  
and thus we have an exact description of the situation - 
again, $M$ is the product of a Levi non-degenerate CR-hypersurface 
in the $\widehat K$-orbit with an $S^1$.    

\medskip\noindent
Finally, we come to the interesting case where $L^\CC$ is acting as
$\CC^*$ and $\widehat P$ is acting transitively on the fiber $\widehat F$.  
Here we consider the $\widehat P$--action on the fiber of
$\widehat G/\widehat I\to \widehat G/\widehat J$.  
This is a $\widehat K$-homogeneous $\CC^*$--principal bundle 
where $\widehat P$ is acting by a character. 

\medskip\noindent
Since $L^\CC$ is acting as $\CC^*$ on $\widehat F$ it follows that
the unipotent radical $R_u(\widehat P)$ is acting with 1-codimensional 
ineffectivity on the fiber    
$\widehat J/\widehat H$.  
Since this ineffectivity must also be a normal subgroup of $\widehat P$, 
it is a root group. 
Thus its complement in $\widehat R_u(\widehat P)$ is an 
$L^\CC$--invariant root group.  
Hence, we have a complete description
of the $\widehat P$--ineffectivity on $\widehat F$ and it follows
that $\widehat P$ splits as a semidirect product of this ineffectivity
and a 2--dimensional solvable group $A$ which is a semidirect product
of a 1-dimensional torus with a 1-dimensional simple root group $U$.

\medskip\noindent
For the moment let us go to a covering of $X$ so that this ineffectivity
is exactly the $\widehat P$-isotropy.  
In this way we reach a well--understood situation, 
because we can then fiber out the full
reductive part of $\widehat P$ to obtain a $\CC^*$--bundle over
a $\CC$-bundle where the simple root group acts transitively on
the fiber.   
This $\CC$-bundle setting, where $K$ has hypersurface orbits, 
has been classified in detail in \cite{AHR}.    
In particular, we find a further fibration 
$\widehat K/\widehat P\to \widehat K/\widehat P_1$ which is
a $\PP_1$--bundle.  
This is given by the root $\SL_2$ associated to
the simple root group.  
Thus we see that the general example just comes from 
the standard 3-dimensional example $\SL_2(\CC)/U_\ZZ$, where
the $\SL_2(\CC)$ is embedded in the parabolic group $\widehat P_1$
as one of its semisimple factors.
Here the isotropy in the bigger group $\widehat K$ is then the 
ineffectivity described above extended by a copy of $\ZZ$ in the 
simple root group at hand.
\end{proof} 

\noindent  
The basic example where $L$ acts on $F$ with 1--dimensional orbits
and $\widehat P$ acts transitively on $\widehat F$ arises as follows:
$$X=\SL_2(\CC)/U_\ZZ\to \SL_2(\CC)/U=\CC^2\setminus \{0\}
\to \SL_2(\CC)/B=\PP_1(\CC)\,.
$$
If we now consider the $\SL_2(\CC)\times \CC^*$--action on
$X$ which is defined by the right $U$-action, then 
the real subgroup that has 2-codimensional orbits is $\SU_2\times S^1$. 
Note that $\SO_3(\CC)=\SL_2(\CC)/(\pm \Id)$ can be used to make   
essentially the same example.    
In that case we take $\widehat K$ to be the group which is acting 
effectively on the base (in this case $\PP_1$) and therefore 
$L$ does indeed have 1-dimensional orbits.  

\medskip\noindent
This example can be put inside a bigger picture, 
i.e., the one that we describe in the classification proof.  
For this it is enough to consider the following example.

\medskip\noindent{\bf Example:}\   
Let $\widehat G=\SL_n(\CC)$ and $X=\widehat G/\widehat P$ be the 
Grassmannian of 2-dimensional subspaces of $\CC^n$. 
Let $\widehat L\ltimes R_u(\widehat P)$ be the standard Levi 
decomposition of $\widehat P$, where 
$\widehat L=Z(\widehat L)(\widehat L_1\times \widehat L_2)$ 
is the stabilizer of the decomposition 
$$
\CC^n=\Span\{e_1,e_2\}\oplus \Span\{e_3,\ldots e_n\}\,.
$$
Here the center $Z(\widehat L)$ acts diagonally.
Now we attempt to build our basic example into $\widehat P$.
Since $\widehat L_1=\SL_2(\CC)$, we may define
$\widehat H=Z(\widehat L)(U_\ZZ(\widehat L_1)\times \widehat L_2)$.
Since $U(\widehat L_1)$ normalizes $\widehat H$, we have the induced
right--action of $U/U_\ZZ$ and, redefining our group as the
larger group $\widehat G\times U/U_\ZZ$, we have the desired 
example.

\medskip\noindent{\bf Remark:}\   
It should be noted that in general the semisimple factor
$\widehat L_1$ may be $\SO_3(\CC)$ and therefore we must
build in the second basic example mentioned above.  
Furthermore, to put this in the context of our classification result in
the setting where $L$ acts as $S^1$ on the fiber, we must choose
$\widehat K$ to be the group which is acting effectively on the 
base $Z$.  

\subsection{The K\"{a}hlerian Parallelizable Case}  

As we have seen, the only remaining serious difficulty in
determining a fine classification in the general 
2-codimensional case is when the fiber of the
$\lie g$--anticanonical bundle is also 2-codimensional,
i.e., when the base is a compact complex manifold.
In this case we must understand the fiber
in its globalization, and therefore it is enough to
consider the case where the base is just a point.
Thus in this section we assume that $M=G/\Gamma $,
where $\Gamma $ is discrete.    

\medskip\noindent   
By the considerations in \S 2.3 one knows that there exists a complexification 
$(X,\sigma)$ of $G/\Gamma$, where we think of $M$ as already 
embedded in $X$, i.e., we identify $M$ with $\sigma (G/\Gamma)$.  
In this section we assume that $M$ is K\"{a}hlerian in the sense that 
$M$ is contained in a tube $Z$ in $X$ that is K\"{a}hler.  
We also recall that the Lie algebra $\lie g$ of $G$ can be complexified 
to a Lie algebra $\widehat{\lie g}$ that acts locally on some open 
neighborhood of $M$ in $Z$, see \S 2.4.  
Note that without loss of generality   
we may assume that this neighborhood is $Z$ itself by appropriately 
shrinking and renaming, if necessary.    
As usual we let $\lie m:=\lie g\, \cap \, i\lie g$ be
the maximal complex ideal in $\lie g$. 
The purpose of this section is to
show that $\lie m$ is solvable in this general K\"ahler case.  
In addition, if $M$ has a CR-structure of codimension 
at most two, then $\lie g$ is solvable.

\begin {theorem}\label {reduction to solvable}
Let $\Gamma $ be a discrete cocompact subgroup of $G$ and assume
that the compact CR-homogeneous space $M:=G/\Gamma $ is K\"ahlerian
in the sense considered above.  
Then the maximal complex ideal $\lie m$ in $\lie g$ is solvable.
In particular, if the CR-structure of $M$ is at most 2-codimensional,   
then $G$ is solvable.
\end  {theorem}

\noindent
For the proof we apply several of the basic tools for the
theory of quotients via K\"ahlerian reduction (see \cite{HH}).  
Here we give a streamlined argument using the strong assumptions at
hand. 
For this let $L$ be the complex semisimple Lie subgroup 
of $G$ which is associated to 
a maximal semisimple subalgebra $\lie l$ of $\lie m$ and let 
$K$ be a fixed choice of a maximal compact subgroup of $L$. 
Note that we may assume that the tube $Z$ is $K$-invariant.  
By averaging the given K\"ahler 
form $\omega $ over $K$ we may assume that it is $K$-invariant.  
Thus, since $K$ is semisimple, the general theory guarantees the 
existence of a unique ($K$--equivariant) 
moment map $\mu :Z\to \lie k^*$.  
Finally, define $M_0:=\mu^{-1}(0)\,\cap \, M$.   
Note that since the $K$--action on $Z$ is locally free, 
$\mu $ is an open mapping of constant rank which is equal
to $\dim_\RR K$.
\begin {proposition}
The 0-fiber $\mu ^{-1}(0)$ is a submanifold of $Z$ which
is transversal to $L.x$ at every point $x\in M_0$.
In particular, $M_0$ is smooth, nonempty and $M^{ss}:=L.M_0$ is open
in $M$.
\end {proposition}
\begin {proof}
The moment map associated to the pullback of $\omega $
to the orbit $L.x$ of a point $x\in M_0$ is just the
pullback of the moment map.  Since the $K$-action is
locally free, the ranks of these
two maps are the same and the transversality statement
is immediate.  Hence, it only remains to show that
$M_0$ is nonempty. But, since $\mu \vert L.x$ is an
open map for all $x\in M$, if $x_0$ is a point of
$M$ where $\Vert \mu \Vert ^2\vert M$ takes on its minimum, it
follows that $x_0\in M_0$.
\end {proof} 
\noindent
By definition, the complement $M_1:=M\setminus M^{ss}$ is 
a compact $L$--invariant subset of $M$.  Just as in the
above argument, if $M_1$ were nonempty and 
$\Vert \mu \Vert ^2\vert M_1$ takes on its minimum at $x_1$, then 
$x_1$ would be in $M_0$, contrary to assumption.  
\begin {corollary}
$$
M^{ss}=L.M_0=M\,.
$$ 
\end {corollary}

\noindent 
Now, for $x\in M_0$ the orbit $L.x$ is retractible to the 
$K$--orbit $K.x$ which is isotropic in $L.x$.  
Since $L.x$ is Stein, it then follows that $\omega =dd^c\rho $ on $L.x$,
where $\rho $ is a $K$--invariant strictly plurisubharmonic function.  
By the results of Azad and Loeb (\cite {AL}), if
$m_0:=\rho (x_0)$, then $\rho :L_x\to [m_0,\infty)$ is a
proper exhaustion with no critical points except for those
where it takes on its minimum along $K.x_0$. In particular,
$L.x\,\cap \, M_0=K.x$ for all $x\in M_0$.

\medskip\noindent
{\it Proof of Theorem \ref{reduction to solvable}.}
Since for every $x\in M$ the orbit $L.x$ is K\"ahler and the 
isotropy group $L_x$ is just the intersection with $L$ of 
some conjugate of $\Gamma $ in $G$, it follows that $L_x$ is
a discrete subgroup of $L$.  
The fact that $L.x$ is K\"ahler
then implies that $L_x$ is finite (\cite{BeOe}).  
By replacing $\Gamma $ by an appropriate subgroup of finite
index, we may assume that $\Gamma $ is torsion free and therefore
that the $L$-action on $M$ is free.  
Since $L.M_0=M$, it therefore
follows that the action map $L\times M_0\to M$, 
$(\ell ,x)\to \ell (x)$, identifies $M$ with the 
\emph{noncompact} bundle space $L\times _KM_0$. 
Unless $L$ is trivial, this is a contradiction.  
Thus $\lie m$ is solvable.  

\medskip \noindent  
Finally, if the CR-structure of $M$ is at most 2-codimensional, then $\lie m$
is also at most 2-codimensional in $\lie g$ and consequently $\lie g$
is also solvable.
\qed

\medskip\noindent
{\bf Remark.}  
In the hypersurface case finer results are known \cite{R}.  
At present we only know rather simple constructions
of examples using solvable (nonabelian) groups in the 2-codimensional case,
e.g., complexifications of compact 2--dimensional solv--manifolds.
As a Leitfaden for more interesting constructions,
it might be of interest to attempt to build examples $X$ with
$\mathcal O(X)=\CC$.\qed 

\bigskip\noindent 
Authors' addresses:

\noindent 
Dept. of Math. \& Stats.  \hfill Fakult\"{a}t und Institut f\"{u}r Mathematik 

\noindent 
University of Regina \hfill Ruhr-Universit\"{a}t Bochum 

\noindent 
Regina, Canada  S4S 4A5 \hfill   D-44780 Bochum,  Germany  

\noindent  
gilligan@math.uregina.ca  \hfill ahuck@cplx.ruhr-uni-bochum.de 

\end{document}